\documentclass[11pt]{article}

\usepackage{color}
\usepackage{amsfonts}
\usepackage{amssymb,amsmath,amsthm, amsfonts}
\usepackage[numbers,sort&compress]{natbib}
\usepackage{graphicx}
\usepackage{url}

\def\sqr#1#2{{\vcenter{\vbox{\hrule height.#2pt
				\hbox{\vrule width.#2pt height#1pt \kern#1pt \vrule width.#2pt}
				\hrule height.#2pt}}}}
\def\signed #1{{\unskip\nobreak\hfil\penalty50
		\hskip2em\hbox{}\nobreak\hfil#1
		\parfillskip=0pt \finalhyphendemerits=0 \par}}
\def\endpf{\signed {$\sqr69$}}
\def\sqr#1#2{{\vcenter{\vbox{\hrule height.#2pt
				\hbox{\vrule width.#2pt height#1pt \kern#1pt \vrule width.#2pt}
				\hrule height.#2pt}}}}
\def\signed #1{{\unskip\nobreak\hfil\penalty50
		\hskip2em\hbox{}\nobreak\hfil#1
		\parfillskip=0pt \finalhyphendemerits=0 \par}}
\def\endpf{\signed {$\sqr69$}}
\def\3n{\negthinspace \negthinspace \negthinspace }
\def\2n{\negthinspace \negthinspace }
\def\1n{\negthinspace }

\def\={\buildrel \triangle \over =}

\def\limsup{\mathop{\overline{\rm lim}}}
\def\liminf{\mathop{\underline{\rm lim}}}

\def\essinf{\mathop{\rm essinf}}
\def\max{\mathop{\rm max}}
\def\min{\mathop{\rm min}}
\def\exp{\mathop{\rm exp}}
\def\sup{\mathop{\rm sup}}
\def\inf{\mathop{\rm inf}}

\def\({\Big (}
\def\){\Big )}
\def\[{\Big[}
\def\]{\Big]}
\def\be{\begin{equation}}
	
	\def\ee{\end{equation}}

\def\square#1{\vbox{\hrule\hbox{\vrule height#1%
			\kern#1\vrule}\hrule}}
\def\rectangle#1#2{\vbox{\hrule\hbox{\vrule height#1%
			\kern#2\vrule}\hrule}}

\font\tenbb=msbm10 \font\sevenbb=msbm7 \font\fivebb=msbm5

\newfam\bbfam
\scriptscriptfont\bbfam=\fivebb \textfont\bbfam=\tenbb
\scriptfont\bbfam=\sevenbb

\newtheorem{lemma}{Lemma}[section]
\newtheorem{remark}{Remark}[section]

\newtheorem{theorem}{Theorem}[section]
\newtheorem{corollary}{Corollary}[section]

\newtheorem{definition}{Definition}[section]
\newtheorem{proposition}{Proposition}[section]

\newtheorem{assumption}{Assumption}[section]

\makeatletter

\@addtoreset{equation}{section}
\makeatother

\begin{document}
	\title{Optimal control problems with generalized mean-field dynamics\\ and viscosity solution to Master Bellman equation\thanks{The work is supported by the NSF of P.R. China (NO. 12031009), the NSF of Shandong Province (NO. ZR2023ZD35), National Key R and D Program of China (NO. 2018YFA0703900), and NSFC-RS (No. 11661130148; NA150344).}
}
\author{Rainer Buckdahn$^{1,2}$,\,\,Juan Li$^{3,2,}$\footnote{Corresponding authors.}
,\,\, Zhanxin Li$^{3,\dagger}$ \\
{$^1$\small Laboratoire de Math\'ematiques de Bretagne Atlantique, Univ Brest,}\\ {\small UMR CNRS 6205, 6 avenue Le Gorgeu, 29200 Brest, France.}\\
{$^2$\small Research Center for Mathematics and Interdisciplinary Sciences,}\\ {\small Shandong University, Qingdao 266237, P.~R.~China.}\\
 {$^3$\small School of Mathematics and Statistics,}\\ {\small Shandong University, Weihai, Weihai 264209, P.~R.~China.}\\
 {\small{\it E-mails: rainer.buckdahn@univ-brest.fr;\,\ juanli@sdu.edu.cn;\,\ zhanxinli@mail.sdu.edu.cn.}}	}	
  \date{August 15, 2024}
	\maketitle
	\noindent\textbf{Abstract}. We study an optimal control problem of generalized mean-field dynamics with open-loop controls, where the coefficients depend not only on the state processes and controls, but also on the joint law of them.  The value function $V$ defined in a  conventional way, but it does not satisfy the Dynamic Programming Principle (DPP for short).
For this reason we introduce subtly a novel value function $\vartheta$, which is closely related to the original value function $V$, such that, a description of $\vartheta$, as a solution of a partial differential equation (PDE), also characterizes $V$. We establish the DPP for $\vartheta$. By using an intrinsic notion of viscosity solutions, initially introduced in Burzoni, Ignazio, Reppen and Soner \cite{BIRS20} and specifically tailored to our framework, we show that the value function $\vartheta$ is a viscosity solution to a Master Bellman equation on a subset of Wasserstein space of probability measures. The uniqueness of viscosity solution is proved for coefficients which depend on the time and the joint law of the control process and the controlled process. Our approach is inspired by Buckdahn, Li, Peng and Rainer \cite{17AP}, and leads to a generalization of the mean-field PDE in \cite{17AP} to a Master Bellman equation in the case of controls.
	\section{Introduction}
In this paper we study a class of general mean-field stochastic optimal control problems (also known as McKean-Vlasov control problems) with open-loop controls, whose state equations are given as follows: For any $0\leq t\leq s\leq T$, $x\in \mathbb{R}^n$, $\zeta \in L^2(\mathcal{F}_t;\mathbb{R}^n)$,
\begin{equation}\label{SDE1.1}
\begin{aligned}
X_{s}^{t, \zeta, u^2}= \ & \ \zeta + \int_t^s\! b_1(r, (X_{r}^{t, \zeta, u^2},u_r^2),  \mathbb{P}_{(X_{r}^{t, \zeta, u^2},u_r^2)}) d r +\int_t^s\!\sigma_1 (r, (X_{r}^{t, \zeta, u^2},u_r^2), \mathbb{P}_{(X_{r}^{t, \zeta, u^2},u_r^2)}) d B_{r},
\end{aligned}
\end{equation}
\vspace{-2mm}
\begin{equation}\label{SDE1.2}
\begin{aligned}
X_{s}^{t, x, \zeta, u}= \ & \ x + \int_t^s\! b_2(r, (X_{r}^{t,x, \zeta, u},u_r^1),  \mathbb{P}_{(X_{r}^{t, \zeta, u^2},u_r^2)}) d r +\int_t^s\!\sigma_2 (r, (X_{r}^{t,x, \zeta, u},u_r^1), \mathbb{P}_{(X_{r}^{t, \zeta, u^2},u_r^2)}) d B_{r},
\end{aligned}
\end{equation}
where the coefficients $b_1,b_2: [0,T]\times (\mathbb{R}^n\times U)\times\mathcal{P}_2(\mathbb{R}^n\times U) \rightarrow \mathbb{R}^n$, $\sigma_1,\sigma_2: [0,T]\times (\mathbb{R}^n\times U)\times\mathcal{P}_2(\mathbb{R}^n\times U) \rightarrow \mathbb{R}^{n\times d}$ depend not only on the state processes and the controls, but also on the joint law of  $X^{t, \zeta, u^2}$ and $u^2$. The control domain $U\subset\mathbb{R}^n$ is a bounded, closed and non-empty set. In the context of financial market, equation (\ref{SDE1.1}) can be interpreted as describing the average state of all agents in the market (or we say, the state of a representative ``mean-field player"), while equation (\ref{SDE1.2}) characterizes the dynamics of an individual agent who faces the entire market (or ``the mean-field player") (see Remark \ref{Fi_Re} for a detailed interpretation). Our mean-field control problem consists in minimizing over all admissible control processes $u=(u^1,u^2)$ the cost functional:
\begin{equation}\label{1.3}
J (t,x,\zeta,u) :=E \big[ \Phi \big(X_{T}^{t, x, \zeta, u},\mathbb{P}_{X_{T}^{t, \zeta, u^2}}\big)\big|\mathcal{F}_t \big],
\end{equation}
where $\Phi: \mathbb{R}^n\times \mathcal{P}_2(\mathbb{R}^n) \rightarrow \mathbb{R}$, and $\big(X^{t, \zeta, u^2}, X^{t, x, \zeta, u}\big)$ is the solution of (\ref{SDE1.1}) and (\ref{SDE1.2}). The value function $V$ is defined as
\begin{equation}\label{1.4}
\begin{aligned}
V(t, x, \zeta):=\essinf_{u=(u^1,u^2) \in \mathcal{U}_{t,T}} J(t, x, \zeta, u).
\end{aligned}
\end{equation}
However, the value function $V$ does not satisfy the DPP. This is why we introduce a novel value function $\vartheta$: 
\begin{equation}\label{1.5}
\begin{aligned}
\vartheta(t,\theta,\mathbb{P}_{\zeta}):=\inf_{u^2\in \mathcal{U}^0_{t,T}}E \big[\mathop{\essinf}\limits_{u^1\in \mathcal{U}^0_{t,T}}E \big[ \Phi \big(X_{T}^{t, \theta, \zeta, (u^1,u^2)},\mathbb{P}_{X_{T}^{t, \zeta, u^2}}\big)\big|\mathcal{F}_t \big]\big],
\end{aligned}
\end{equation}
where $\theta\in L^2(\mathcal{F}_t;\mathbb{R}^n)$ (for the notations $\mathcal{U}_{t,T}$ and $\mathcal{U}_{t,T}^0$ see Section 2). As we show, $\vartheta(t,x,{\zeta})=V(t,x,{\zeta})$, for all $(t,x,{\zeta})\in [0,T]\times \mathbb{R}^n\times L^2(\mathcal{F}_t;\mathbb{R}^n)$. This allows to characterize $V$ through studying $\vartheta$.
 For this we show that the value function $\vartheta$ satisfies the dynamic programming principle (DPP). This is the key to characterize $\vartheta$ as the  unique viscosity solution of the partial differential equation (\ref{HJB}) on $[0,T)\times \mathcal{M}\times \mathcal{M}$, where $\mathcal{M}\subset\mathcal{P}_2(\mathbb{R})$ is the set of probability measures with $\delta$-exponential moments (for the definition of $\mathcal{M}$, see (\ref{expo})). We emphasize that in the no-control case, under the suitable regularity conditions on the coefficients, for $b_1=b_2$, $\sigma_1=\sigma_2$, the obtained mean field PDE coincides with that obtained in \cite{17AP}, see Remark \ref{Fi_Re-4.2}.

Mean-field optimal control problems have been knowing a growing interest with the emergence of mean-field games, which have been introduced and studied independently by Lasry and Lions \cite{LL07} (the reader is also referred to Lions' lectures at \textit{Coll\`{e}ge de France} \cite{L13}, and their notes redacted by Cardaliaguet \cite{C13}), and by Huang, Caines and Malham\'{e} \cite{HCM06}. Both the mean-field control problems and the mean-field games can be interpreted as describing equilibrium states of large population of agents who interact with each other in a weak and symmetrical manner.
The main difference between these both problems is that they adopt different notions of equilibrium. Mean-field games describe the Nash  non-cooperative equilibrium for a large population, while mean-field control problems mainly describe Pareto optimality, where one wants to describe  an infinite number of players with common behaviour by a representative agent. This agent represents in some sense the average behaviour of the infinite number agents. For a detailed discussion of the both problems and the differences between them, the reader is namely referred to Carmona, Delarue and Lachapelle \cite{CDL13}, Lacker \cite{L17} and the references cited therein.

In the literature there are mainly two approaches to study stochastic optimal control problems. On one hand, many works have been focused on determining 
 necessary conditions satisfied by an optimal control, known as the stochastic maximum principle. For control problems with general McKean-Vlasov dynamics, one can refer to Carmona and Delarue \cite{CD15}, Acciaio et al. \cite{ABC19} and the references therein. They concern the study of Pontryagin's maximum principle under assumptions supposing the convexity of either the control state space or the Hamiltonian. On the other hand, Buckdahn, Li and Ma \cite{BLM16} and Buckdahn, Chen and Li \cite{BCL21} investigated Peng's stochastic maximum principle without any convexity assumptions. Besides the study of the stochastic maximum principle, dynamic programming approaches play an important role in the literature.
They relate the value function of control problems to associated Partial Differential Equations (PDEs for short).
Pham and Wei \cite{PW18} studied the mean-field stochastic control problems where the coefficients of the state equation may depend on the joint law of the state and control. Unlike us, they have used feedback controls and reformulated the problem into a deterministic control problem. Bayraktar, Cosso, and Pham \cite{BCP18} have studied optimal control problems of mean-field type for open-loop controls, but without the dependence of the law on the control. They have proved a so-called Randomized DPP, whose initial idea was introduced in \cite{KMPZ10} in connection with impulse control, and they have got a characterization for value function $V$ through an auxiliary intensity control problem for a Poisson random measure.

In this paper, our main objective is to use a dynamic programming approach to study the mean-field control problem (\ref{SDE1.1})-(\ref{1.4}). 
We focus on the case where the coefficients depend on the joint law of state process, and control process, as well as the controls are open-loop.
 Let us emphasize that, when the coefficients do not depend on the law of control process and we select the same control for both (\ref{SDE1.1}) and (\ref{SDE1.2}), i.e., $u^1=u^2$, the state processes are of the same form as those studied in \cite{BCP18}.

The first difficulty encountered in our study lies in the absence of the DPP for the value function $V$ (as indicated by equation (\ref{DPPNO}) and the subsequent discussion), which makes it impossible to derive the PDE characterization for $V$ by using DPP.
To tackle this problem, based on $V$ we introduce a novel definition of value function $\vartheta$ (see (\ref{1.5})). On one hand, $\vartheta$ satisfies the DPP, and it allows to derive that $\vartheta$ is the unique viscosity solution of the associated Master Bellman equation (a PDE defined on the Wasserstein space of probability measures). On the other hand, $\vartheta$ is constructed subtly to satisfy $\vartheta(t,x,{\zeta})=V(t,x,{\zeta})$, for all $(t,x,{\zeta})\in [0,T]\times \mathbb{R}^n\times L^2(\mathcal{F}_t;\mathbb{R}^n)$, which means that, the  description of $\vartheta$ as a solution of the Master Bellman equation also characterizes $V$.

The second challenge is the uniqueness of the viscosity solution of the second-order PDE on Wasserstein space. Our goal is to characterize the value function $\vartheta$ as the unique viscosity solution of its associated PDE (\ref{HJB}) over a Wasserstein space of probability measures. However, the existing literature has mainly focused on the existence of viscosity solutions, proving that the value function of a control problem is a viscosity solution of the associated PDE. In some special cases, the uniqueness of viscosity solutions has been explored. For example, Wu and Zhang \cite{WZ20} adopted a new definition of viscosity solutions for the path-dependent mean-field control problems, and their maximum/minimum conditions for test functions are stated on compact subsets of the Wasserstein space. Burzoni, Ignazio, Reppen and Soner \cite{BIRS20} have studied viscosity solutions for a particular class of integro-differential Master equations. The uniqueness of viscosity solutions has been proved on Wasserstein spaces of probability measures which have finite exponential moments. For this, the authors have considered deterministic functions of the time as the control processes. Recently, by using refinements of early ideas from the theory of viscosity solutions \cite{L83},
Cosso et al. \cite{CGKPR21} 
  demonstrated the uniqueness of the viscosity solutions on $\mathcal{P}_2(\mathbb{R}^d)$, but only for coefficients which do not depend on the law of the control. These studies use the intrinsic notion of viscosity solutions without lifting (for intrinsic and extrinsic notions of viscosity solutions, see, for example Remark 3.6 of \cite{CGKPR21}). In this paper we adopt the notion of viscosity solution introduced in \cite{BIRS20} and adapt it to our framework of mean-field control problems, in order to derive that value function $\vartheta$ is a viscosity solution of PDE (\ref{HJB}). Under the additional assumption that coefficients only depend on time and joint law of the control process and controlled process, we also  get the uniqueness of the viscosity solution.
	
	The paper is organized as follows. Section 2 is devoted to the formulation of our mean-field optimal control problem.
In Section 3, we study the properties of the value function $V$ and we show that the DPP for $V$ does not hold. To address this issue, we introduce a novel value function, denoted by $\vartheta$, which is closely related to the value function $V$, and we prove the DPP for $\vartheta$. Section 4 is devoted to the introduction of the Master Bellman equation and the definition of viscosity solutions. Some preparations are also made in Section 4, which allow to prove in Section 5 our main result that the value function $\vartheta$ is the unique viscosity solution to the Master Bellman equation.
	
\section{Formulation of the mean-field stochastic control problem}	
\subsection{Notations}
\noindent In this subsection, we introduce the setting in which we study our mean-field stochastic control problem. Throughout our paper we work on the classical Wiener space $(\Omega,\mathcal{F},\mathbb{P})$, and the driving Brownian Motion $B$ is the coordinate process on $\Omega$. To be more precise, for a fixed horizon $T>0$, $\Omega$ is the set of all continuous functions from $[0,T]$ to $\mathbb{R}^d$ starting from $0$ ($\Omega=C_0([0,T]; \mathbb{R}^d)$), $\mathbb{P}$ is Wiener measure on $(\Omega, \mathcal{B}(\Omega))$, $\mathcal{F}=\mathcal{B}(\Omega)\vee \mathcal{N}_{\mathbb{P}} $, where $\mathcal{N}_{\mathbb{P}}$ is the set of all $\mathbb{P}$-null subsets, and the $d$-dimensional coordinate process is given by $B_s(\omega)=\omega_s,\ s\in [0,T],\ \omega\in\Omega$.
By $\mathbb{F}=\left\{\mathcal{F}_{s}, 0 \leq s \leq T\right\}$ we denote the filtration generated by $\{B_s\}_{s\in [0,T]}$ and augmented by all $\mathbb{P}$-null sets, i.e.,
$$
\mathcal{F}_{s}=\sigma\left\{B_{r},0\leq r \leq s\right\} \vee \mathcal{N}_{\mathbb{P}}, \ \ s \in[0, T].
$$
 For any $n\geq 1$, $|z|$ denotes the Euclidean norm of $z\in \mathbb{R}^n$. We denote the law of a random variable $\zeta$ over $(\Omega, \mathcal{F},\mathbb{P})$ by $\mathbb{P}_{\zeta}$, and so the joint distribution of two random variables $\xi$ and $\zeta$ by $\mathbb{P}_{(\xi,\zeta)}$.

 For $k\geq 1$, $\mathcal{P}_{2}(\mathbb{R}^{k})$ denotes the space of the probability measures over $\mathbb{R}^k$ with finite second moment and is endowed with the $2$-Wasserstein metric:
$$\mathcal{W}_{2}(\mu,\nu):=\inf\Big\{\Big(\int_{\mathbb{R}^k\times\mathbb{R}^k}|x-y|^{2}\rho(dxdy)\Big)^{\frac{1}{2}},\ \rho\in \Pi_{\pi,\nu}\Big\},$$
where $\Pi_{\mu,\nu}=\big\{\rho\in\mathcal{P}_{2}(\mathbb{R}^{2k}),\ \text{with} \  \rho(\cdot\times\mathbb{R}^k)=\mu,\ \rho(\mathbb{R}^k\times\cdot)=\nu\big\},$ $\mu,\nu\in \mathcal{P}_{2}(\mathbb{R}^k).$ Recall that $(\mathcal{P}_{2}(\mathbb{R}^{k}),\mathcal{W}_{2})$ is a complete separable space (see, for example Theorem 6.18 of \cite{V09}). Now we introduce the following spaces which will be used frequently in what follows:
	
	\noindent\quad$\bullet\ L^2(\mathcal{F}_t;\mathbb{R}^n)$ is the set of $\mathbb{R}^n$-valued, $\mathcal{F}_t$-measurable random variables $\zeta:\Omega\rightarrow \mathbb{R}^n$ such that $\displaystyle E[|\zeta|^2]<\infty$.

	\noindent\quad$\bullet\ L^2_{\mathbb{F}}([0,T];\mathbb{R}^n)$ is the set of $\mathbb{R}^n$-valued, $\mathbb{F}$\text{-progressively measurable processes} $\phi:\Omega\times [0,T]\rightarrow \mathbb{R}^n$ such that $\displaystyle E\big[\int_0^T |\phi_t|^2 dt\big]<+\infty$.

\noindent\quad$\bullet\ \mathcal{S}^2(0,T;\mathbb{R}^n)$ is the set of $\mathbb{F}$-adapted continuous processes $\phi: \Omega\times[0,T]\rightarrow \mathbb{R}^n$ satisfying
$E\big[\sup\limits_{0\leq s\leq T}|\phi_{s}|^{2}\big]<\infty.$


For simplicity of notations, when $n=1$,
we write $L^2(\mathcal{F}_t):=L^2(\mathcal{F}_t;\mathbb{R})$, and $L^2_{\mathbb{F}}([0,T]):=L^2_{\mathbb{F}}([0,T];\mathbb{R})$.

\subsection{Mean-field stochastic optimal control problems}
Now we introduce our mean-field stochastic control problems. For $0\leq t\leq s\leq T$, the set of admissible control processes on $[t,s]$, denoted by $\mathcal{U}^0_{t,s}$, is the set of all $U$-valued $\mathbb{F}$-adapted stochastic processes $\phi:\Omega\times [t,s]\rightarrow\mathbb{R}^n$, where the control domain $U\subset\mathbb{R}^n$ is a bounded, closed set and $U\neq\varnothing$. We identify two processes $u$ and $\overline{u}$ in $\mathcal{U}^0_{t,s}$ and write $u\equiv \overline{u}$ on $[t,s]$, if $\mathbb{P}\{u=\overline{u}\ \text{a.e.}\ \text{on}\ [t,s]\}=1$. 
For simplicity of notations, we put $\mathcal{U}_{t,s}:=\mathcal{U}^0_{t,s}\times\mathcal{U}^0_{t,s}$.

Our dynamics are the following controlled mean-field SDEs:
\begin{equation}\label{SDE1}
\begin{aligned}
X_{s}^{t, \zeta, u^2}= \ & \ \zeta + \int_t^s b_1(r, (X_{r}^{t, \zeta, u^2},u_r^2),  \mathbb{P}_{(X_{r}^{t, \zeta, u^2},u_r^2)}) d r +\int_t^s\sigma_1 (r, (X_{r}^{t, \zeta, u^2},u_r^2), \mathbb{P}_{(X_{r}^{t, \zeta, u^2},u_r^2)}) d B_{r},
\end{aligned}
\end{equation}
\vspace{-2mm}
\begin{equation}\label{SDE2}
\begin{aligned}
X_{s}^{t, x, \zeta, u}= \ & \ x + \int_t^s\! b_2(r, (X_{r}^{t,x, \zeta, u},u_r^1),  \mathbb{P}_{(X_{r}^{t, \zeta, u^2},u_r^2)}) d r +\int_t^s\!\sigma_2 (r, (X_{r}^{t,x, \zeta, u},u_r^1), \mathbb{P}_{(X_{r}^{t, \zeta, u^2},u_r^2)}) d B_{r},
\end{aligned}
\end{equation}
where\ $t\in [0,T],\ s\in [t,T],\ x\in\mathbb{R}^n,\ \zeta \in L^2(\mathcal{F}_t;\mathbb{R}^n), \ u=(u^1,u^2)\in \mathcal{U}_{t,T}.$
\begin{remark}\label{Fi_Re}
It can be observed that equation (\ref{SDE1}) is a mean-field SDE, and if we substitute the solution of equation (\ref{SDE1}) into equation (\ref{SDE2}), then equation (\ref{SDE2}) becomes a classical SDE.
In financial markets, equation (\ref{SDE1}) can interpreted as the dynamics of a representative agent in the market (referred to as the ``mean-field player"), who plays collectively with his/her ``collective control" $u^2$. It can be understood as describing the average state of all agents in the market.
 On the other hand, equation (\ref{SDE2}) describes the dynamics of an individual agent who faces the ``mean-field player",  and $u^1$ is the ``individual control" played by this individual agent.  Thus, evolution (\ref{SDE2}) is affected not only by its own state process and control $(X^{t,x, \zeta, u},u^1)$, but also the law of $(X^{t, \zeta, u^2},u^2)$.
In other words, (\ref{SDE1}) characterizes the evolution of the law $\displaystyle\mathbb{P}_{(X^{t, \zeta, u^2},u^2)}$, while (\ref{SDE2}) describes the associated trajectories of an individual agent with initial condition $X_{t}^{t, x, \zeta, u}=x$. This interpretation aligns closely with practical observations in financial markets.
%
%
\end{remark}
The coefficients are supposed to satisfy the following assumptions:
\begin{assumption}\label{A1}\ \\
\emph{\textbf{(i)}}
 The mappings $b_1,b_2: [0,T]\times (\mathbb{R}^n\times U)\times\mathcal{P}_2(\mathbb{R}^n\times U) \rightarrow \mathbb{R}^n$ and $\sigma_1,\sigma_2: [0,T]\times (\mathbb{R}^n\times U)\times\mathcal{P}_2(\mathbb{R}^n\times U) \rightarrow \mathbb{R}^{n\times d}$ are continuous and there exists a constant $C_0>0$ such that, for all $t\in [0,T],\ x\in\mathbb{R}^n ,\ u\in U,\ \zeta \in L^2(\mathcal{F}_t;\mathbb{R}^n),\ \eta \in L^2(\mathcal{F}_t;U),$
 $$
 \big|\phi(t, (x, u),\mathbb{P}_{(\zeta,\eta)})\big|
 \leq C_0,\quad \phi=b_1,b_2,\sigma_1,\sigma_2.
 $$
\emph{\textbf{(ii)}}
 There exists a constant $ C>0$, such that
 \vspace{-2mm}
$$
\begin{aligned}
\big|\phi(t, (x, u),\mathbb{P}_{(\zeta,\eta)})-& \phi(t, (x', u), \mathbb{P}_{(\zeta ',\eta)})\big|
 \leq C(|x-x'|+\mathcal{W}_{2}(\mathbb{P}_{(\zeta ,\eta)},\mathbb{P}_{(\zeta ',\eta)})), \\
\end{aligned}
$$
\vspace{-0.5mm}
\hspace{-3mm}for all\ $t\in [0,T],\ x, x'\in\mathbb{R}^n ,\ u\in U,\ \zeta,\zeta' \in L^2(\mathcal{F}_t;\mathbb{R}^n),\ \eta \in L^2(\mathcal{F}_t;U),$\ $\phi=b_1,b_2,\sigma_1,\sigma_2. $
\end{assumption}
Under Assumption \ref{A1}, there exists a unique pair of solutions $\big(X_{s}^{t, \zeta, u^2},X_{s}^{t, x, \zeta, u}\big)_{s\in [t,T]}\in \mathcal{S}^2(0,T;\mathbb{R}^n)\times \mathcal{S}^2(0,T;\mathbb{R}^n)$ to the equations (\ref{SDE1}) and (\ref{SDE2}) (see, e.g., Buckdahn, Li, Peng and Rainer \cite{17AP}, Carmona and Delarue \cite{CD15}). Moreover, for every $p\geq 2$, there exists $C_p\in \mathbb{R}^+$ such that, for all $t\in[0,T],\ u=(u^1,u^2)\in\mathcal{U}_{t,T}$, and $x, x^{\prime} \in \mathbb{R}^n,\ \zeta, \zeta^{\prime} \in L^{2}\left(\mathcal{F}_{t};\mathbb{R}^n\right)$, we have the following estimates:
\begin{equation}\label{Es_SDE}
\begin{gathered}
E\left[\sup _{s \in [t,T]}\left|X_{s}^{t, \zeta, u^{2}}-X_{s}^{t,\zeta^{\prime}, u^{2}}\right|^{p}\right] \leq C_{p}E\left[\left|\zeta-\zeta^{\prime}\right|^{p}\right] ,\\
E\left[\sup _{s \in[t,T]}\left|X_{s}^{t, x, \zeta, u}-X_{s}^{t,x', \zeta^{\prime}, u}\right|^{p}\right] \leq C_p \left(\left|x-x^{\prime}\right|^{p}+E\left[\left|\zeta-\zeta^{\prime}\right|^{p}\right]\right) .
\end{gathered}
\end{equation}
From the uniqueness of the solutions of the both equations, we deduce the following \textit{flow property}: For all $0 \leq t<t+\delta \leq T, \ x \in \mathbb{R}^n,\ \zeta \in L^{2}\left(\mathcal{F}_{t};\mathbb{R}^n\right), u=(u^{1}, u^{2}) \in \mathcal{U}_{t,T}$:
$$\left(X_{s}^{t+\delta, X_{t+\delta}^{t,x, \zeta, u}, X_{t+\delta}^{t,\zeta, u^{2}}, u}, X_{s}^{t+\delta, X_{t+\delta}^{t,\zeta, u^{2}}, u^{2}}\right)=\left(X_{s}^{t, x, \zeta, u}, X_{s}^{t, \zeta, u^{2}}\right),\quad  s \in[t+\delta, T],\ \mathbb{P}\text{-a.s.}$$

\begin{remark}\label{3.2}
As the control $\{u_s=(u^1_s,u^2_s):s\geq t\}$ is $\mathbb{F}$-adapted,  $X^{t,x,\zeta,u}$ is, in general, not independent of $\mathcal{F}_t$, and therefore not independent of $\zeta$.  We note that, although we will not use it in the following discussion,  $X^{t,x,\zeta,u}$ depends on $\zeta$ through the joint law of $(\zeta, \mathring{u}^2)$, where $\mathring{u}^2$ is defined by $u^2$, see the following paragraph for the detailed definition.

For any $u^2\in\mathcal{U}^0_{t,T}$, there exists a measurable nonanticipating functional $f:\ [0,T]\times C_0([0,T];\mathbb{R}^d)\rightarrow U$ such that \begin{equation}
u_s^2=f_s(B_{\cdot\wedge s})=f_s(B_{\cdot\wedge t}+\widehat{B}_{\cdot\wedge s}), \quad s\in[t,T],\ \mathbb{P}\mbox{-a.s.},
\end{equation}
where $\widehat{B}:= B_{\cdot\vee t}-B_t$. Recall that $f:\ [0,T]\times C_0([0,T];\mathbb{R}^d)\rightarrow U$ is nonanticipating means that for all $s\in[t,T]$ and all $\psi,\psi'\in C_0([0,T];\mathbb{R}^d)$ with $\psi(r)=\psi'(r), \ r\in[t,s]$, it holds that $f_r(\psi)=f_r(\psi')$, $r\in[t,s].$
For any $\varphi\in C_0([t,T];\mathbb{R}^d),\ \omega\in\Omega$, we define $\mathring{u}^2_s(\omega,\varphi):=f_s(B_{\cdot\wedge t}(\omega)+\varphi_{\cdot\wedge s})$. Obviously, we know that for any $s\in[t,T],\ \omega\in\Omega$, the mapping $\mathring{u}^2_s(\omega,\cdot): C_0([t,T];\mathbb{R}^d)\rightarrow U$ satisfies $\mathring{u}^2_s(\omega,\widehat{B}_{\cdot\wedge s}(\omega))=u^2_s(\omega)$. Let
$
\mathcal{H}_t:=\Big\{\psi:[t,T]\times C_0([t,T];\mathbb{R}^d)\rightarrow U\mbox{ measurable nonanticipating functional}\ \Big|\ \|\psi\|_{\mathcal{H}_t}^2:=E\big[\int_t^T|\psi_s(\widehat{B}_{\cdot\wedge s})|^2 ds\big]<+\infty \Big\}
$. It is easy to check that $\mathring{u}^2(\omega,\cdot)=\big(\mathring{u}^2_s(\omega,\cdot)\big)_{s\in [t,T]}\in \mathcal{H}_t$.
Moreover, $\mathring{u}^2: (\Omega,\mathcal{F}_t)\rightarrow (\mathcal{H}_t,\mathcal{B}(\mathcal{H}_t))$ is measurable, and $\mathring{u}^2\in L^2(\mathcal{F}_t;\mathcal{H}_t)$\rm{:}
$$
E[\|\mathring{u}^2\|^2_{\mathcal{H}_t}]=E[\int_t^T|u_s^2(B_{\cdot\wedge t}+\widehat{B}_{\cdot\wedge s})|^2 ds]<+\infty.
$$
On the other hand, for any given $\mathring{u}^2\in L^2(\mathcal{F};\mathcal{H}_t)$, we define the random variable $u_s^2:=\mathring{u}_s^2(\cdot,\widehat{B}_{\cdot\wedge s})$, $s\in[t,T]$, then $u^2\in\mathcal{U}^0_{t,T}$. Hence, we can identify $u^2\in\mathcal{U}^0_{t,T}$ with $\mathring{u}^2\in L^2(\mathcal{F}_t;\mathcal{H}_t)$.

Let now $(\zeta, \mathring{u}^2), ({\zeta}',\mathring{u}^{2\prime})\in L^2(\mathcal{F}_t;\mathbb{R}^n)\times L^2(\mathcal{F}_t,\mathcal{H}_t)$ be such that $\mathbb{P}_{(\zeta, \mathring{u}^2)}=\mathbb{P}_{({\zeta}', \mathring{u}^{2\prime})}.$ As $(\zeta, \mathring{u}^2),\ ({\zeta}',\mathring{u}^{2\prime})$ are independent of $\widehat{B}$, we have
\begin{equation*}
\begin{aligned}
&\mathbb{P}_{(\zeta, \mathring{u}^2)}=\mathbb{P}_{({\zeta}', \mathring{u}^{2\prime})}\ \iff \
\mathbb{P}_{(\zeta, \mathring{u}^2, \widehat{B})}=\mathbb{P}_{({\zeta}', \mathring{u}^{2\prime},\widehat{B})} \ \iff \\& \
\mathbb{P}_{(\zeta, \mathring{u}^2(\cdot,\widehat{B}), \widehat{B})}=\mathbb{P}_{({\zeta}', \mathring{u}^{2\prime}(\cdot,\widehat{B}),\widehat{B})}
\ \iff \
\mathbb{P}_{(\zeta, {u}^2, \widehat{B})}=\mathbb{P}_{({\zeta}', {u}^{2\prime},\widehat{B})},
\end{aligned}
\end{equation*}
where $u^2:=\mathring{u}^2(\cdot,\widehat{B}),\ u^{2\prime}:=\mathring{u}^{2\prime}(\cdot,\widehat{B}).$
 As we have the strong uniqueness of SDE (\ref{SDE1}) with control $u^2$, we also have the uniqueness in law. In particular, it follows
$$\mathbb{P}_{X_{s}^{t, \zeta, u^2}}=\mathbb{P}_{X_{s}^{t, \zeta', u^{2\prime}}}, \quad s\in [t,T].$$
For $u^1\in\mathcal{U}^0_{t,T}$, this substituted in SDE (\ref{SDE2}), combined with the strong uniqueness of SDE (\ref{SDE2}) yields:
$$
X_s^{t,x,\zeta,(u^1,u^2)}=X_s^{t,x,\zeta',(u^1,{u^2}')},\quad s\in[t,T],\ \mathbb{P}\text{-a.s.}
$$
This allows to define $X_s^{t,x,u^1,\mathbb{P}_{(\zeta,\mathring{u}^2)}}:=X_s^{t,x,\zeta,(u^1,u^2)},\  s\in [t,T],$ which means $X^{t,x,\zeta,u}$ depends on $\zeta$ through the joint law of $(\zeta, \mathring{u}^2)$, or, equivalently, on the joint law of $(\zeta, u^2, \widehat{B})$.
\end{remark}

\begin{assumption}\label{A2}
Let now be given a Borel measurable function $\Phi: \mathbb{R}^n\times \mathcal{P}_2(\mathbb{R}^n) \rightarrow \mathbb{R}$ that satisfies: There is a constant $C>0$ such that, for all $x, x^{\prime} \in \mathbb{R}^n$, $\mu,\mu'\in\mathcal{P}_2(\mathbb{R}^n),$
\begin{equation}\label{PhiLip}
\left|\Phi(x,\mu)-\Phi\left(x^{\prime},\mu'\right)\right| \leq C\big(\left|x-x^{\prime}\right|+\mathcal{W}_2(\mu,\mu^{\prime})\big).
\end{equation}
\end{assumption}
 Given the control processes $u=\left(u^1, u^{2}\right) \in \mathcal{U}_{t, T}$, we introduce \textbf{the cost functional} of our mean-field stochastic optimal control problem:
\begin{equation}\label{cost1}
J (t,x,\zeta,u) :=E \big[ \Phi \big(X_{T}^{t, x, \zeta, u},\mathbb{P}_{X_{T}^{t, \zeta, u^2}}\big)\big|\mathcal{F}_t \big],\quad (t,x)\in [0,T]\times\mathbb{R}^n,\ \zeta\in L^2(\mathcal{F}_t;\mathbb{R}^n),
\end{equation}
where $\big(X^{t, \zeta, u^2},X^{t, x, \zeta, u}\big)$ are the solutions of (\ref{SDE1}) and (\ref{SDE2}). We first consider \textbf{the common definition of value function}. Supposing that both the ``mean-field player" and the ``individual player" try to minimize the cost, we define for all $(t,x)\in [0,T]\times\mathbb{R}^n,\ \zeta\in L^2(\mathcal{F}_t;\mathbb{R}^n)$,
\begin{equation}\label{value}
\begin{aligned}
V(t, x, \zeta):=&\ \essinf_{u \in \mathcal{U}_{t,T}} J(t, x, \zeta, u)\\
=&\ \essinf_{u^2 \in \mathcal{U}^0_{t,T}}\left(\essinf_{u^1 \in \mathcal{U}^0_{t,T}} J\left(t, x, \zeta,\left(u^{1}, u^{2}\right)\right)\right).
\end{aligned}
\end{equation}
We denote $\displaystyle W\left(t, x, \zeta, u^{2}\right):=\essinf_{u^1 \in \mathcal{U}^0_{t,T}} J\left(t, x, \zeta,\left(u^{1}, u^{2}\right)\right)$, $(t,x)\in [0,T]\times\mathbb{R}^n,\ u^2\in\mathcal{U}^0_{t,T}$. Then
$$
V(t,x,\zeta)=\essinf_{u^2 \in \mathcal{U}^0_{t,T}} W\left(t, x, \zeta, u^{2}\right).
$$
\begin{remark}
As mentioned in Remark \ref{3.2}, $X$ is in general not independent of $\mathcal{F}_t$, and therefore $J(t,x,\zeta,u)$ is an $\mathcal{F}_t$-measurable random variable.
\end{remark}
\section{Dynamic programming principle}
In this section, we will always assume that Assumptions \ref{A1} and \ref{A2} hold true without explicitly stating them. Combining standard arguments from Buckdahn and Li \cite{BL08} with a novel, rather subtle approach, we get the following lemma:
\begin{lemma}\label{W_deter}
For any $(t,x)\in [0,T]\times\mathbb{R}^n$, $\zeta\in L^2(\mathcal{F}_t;\mathbb{R}^n),\ u^2\in\mathcal{U}^0_{t,T}$, $W\left(t, x, \zeta, u^{2}\right)$ is deterministic. Moreover, there is some $C\in\mathbb{R}_+$ such that
$$
|W(t,x,\zeta,u^2)-W(t,x',\zeta',u^2)|\leq C \big(|x-x'|+\big(E[|\zeta-\zeta'|^2]\big)^{\frac{1}{2}}\big),
$$
for all $t\in [0,T]$, $x,x'\in\mathbb{R}^n$, $\zeta,\zeta'\in L^2(\mathcal{F}_t;\mathbb{R}^n)$, $u^2\in\mathcal{U}^0_{t,T}$.
\end{lemma}
\noindent\textit{Proof}. Recall that the space $(\Omega,\mathcal{F},\mathbb{P})$ on which we work is the classical Wiener space, $\Omega=C_0([0,T];\mathbb{R}^d)$, the Brownian motion is the coordinate process $B_t(\omega)=\omega(t)$, $t\in[0,T]$, $\omega\in\Omega$. Let $H$ denote the Cameron-Martin space, $\displaystyle H:=\big\{(\int_0^sh_rdr)_{s\in[0,T]}|\ h\in L^2([0,T];\mathbb{R}^d)\big\}$, and we define $\displaystyle H_t:=\big\{(\int_0^{t\wedge s}h_rdr)_{s\in[0,T]}|\ h\in L^2([0,T];\mathbb{R}^d)\big\}$. Obviously, $H_t\subset H\subset C_0([0,T];\mathbb{R}^d)=\Omega$. For $\omega\in\Omega$, we define the mapping $\displaystyle\tau_h(\omega):=\omega+\int_0^{t\wedge \cdot}h_r dr$, where $\displaystyle\int_0^{t\wedge \cdot}h_r dr\in H_t$. It is easy to see that $\tau_h:\Omega\rightarrow \Omega$ is a bijection, and its law is given by
$$\frac{d\mathbb{P}_{\tau_h}}{d\mathbb{P}}=\exp\Big\{\int_{0}^th_sdB_s-\frac{1}{2}\int_0^t|h_s|^2ds\Big\}.$$
Applying the Girsanov transformation to SDE (\ref{SDE2}), we obtain:
\begin{equation}\label{*}
\begin{aligned}
X_s^{t,x,\zeta,u}\circ\tau_h=&x+\int_t^s b_2\big(r,X_r^{t,x,\zeta,u}\circ\tau_h, u_r^1(\tau_h),\mathbb{P}_{(X_r^{t,\zeta,u^2},u_r^2)}\big)dr\\
&+\int_t^s \sigma_2\big(r,X_r^{t,x,\zeta,u}\circ\tau_h, u_r^1(\tau_h),\mathbb{P}_{(X_r^{t,\zeta,u^2},u_r^2)}\big)dB_r(\tau_h),\  s\in[t,T], \ \mathbb{P}_{\tau_h}(\sim\mathbb{P})\text{-a.s.}
\end{aligned}
\end{equation}
Notice that $dB_s(\tau_h)=dB_s,\ s\in[t,T],$ since $\int_0^{t\wedge s}h_rdr=\int_0^th_rdr,\ s\in[t,T]$. So (\ref{*}) is the same SDE as SDE (\ref{SDE2}), but governed by the Girsanov transformed control $u^1(\tau_h)$. Obviously, also $u^1(\tau_h)$ belongs to $\mathcal{U}^0_{t,s}$, and from the uniqueness of the solution of SDE (\ref{SDE2}) it follows that
\begin{equation}
X_s^{t,x,\zeta,u}\circ \tau_h=X_s^{t,x,\zeta,(u^1(\tau_h),u^2)},\quad s\in[t,T],\ h\in L^2([0,T];\mathbb{R}^d).
\end{equation}
By a similar Girsanov transformation argument applied to the following BSDE
$$Y_t=\Phi\big(X_T^{t,x,\zeta,(u^1,u^2)},\mathbb{P}_{X_T^{t,\zeta,u^2}}\big)-\int_t^TZ_rdB_r,$$
we get from the uniqueness of the solution of the BSDE that
\begin{equation}
E\Big[\Phi\big(X_T^{t,x,\zeta,(u^1,u^2)},\mathbb{P}_{X_T^{t,\zeta,u^2}}\big)\big|\mathcal{F}_t\Big]\circ \tau_h=E\Big[\Phi\big(X_T^{t,x,\zeta,(u^1,u^2)}\circ\tau_h,\mathbb{P}_{X_T^{t,\zeta,u^2}}\big)\big|\mathcal{F}_t\Big].
\end{equation}
Consequently,
\begin{equation}
\begin{aligned}
&J(t,x,\zeta,(u^1,u^2))\circ \tau_h=E\Big[\Phi\big(X_T^{t,x,\zeta,(u^1,u^2)}\circ\tau_h,\mathbb{P}_{X_T^{t,\zeta,u^2}}\big)\big|\mathcal{F}_t\Big]\\
&=E\Big[\Phi\big(X_T^{t,x,\zeta,(u^1(\tau_h),u^2)},\mathbb{P}_{X_T^{t,\zeta,u^2}}\big)\big|\mathcal{F}_t\Big]=J(t,x,\zeta,(u^1(\tau_h),u^2)),\quad \mathbb{P}\text{-a.s.,}
\end{aligned}
\end{equation}
and
\begin{equation}
\begin{aligned}
&W(t,x,\zeta,u^2)\circ\tau_h=\essinf_{u^1\in\mathcal{U}^0_{t,T}}\big(J(t,x,\zeta,(u^1,u^2))\circ\tau_h\big)\\
&=\essinf_{u^1\in\mathcal{U}^0_{t,T}}J(t,x,\zeta,(u^1(\tau_h),u^2))=W(t,x,\zeta,u^2),\quad \mathbb{P}\text{-a.s.,}
\end{aligned}
\end{equation}
as $\mathcal{U}^0_{t,T}=\{u^1(\tau_h)|\ u^1\in\mathcal{U}^0_{t,T}\}.$ Indeed, $u^1=u^1(\tau_{-h})\circ\tau_h$, for all $u^1\in\mathcal{U}^0_{t,T}$; for details the reader is referred to Buckdahn and Li \cite{BL08}. Therefore, $W(t,x,\zeta,u^2)\circ\tau_h=W(t,x,\zeta,u^2)$, $\mathbb{P}$-a.s., for all $h\in L^2([0,T];\mathbb{R}^d)$, and since $W(t,x,\zeta,u^2)$ is $\mathcal{F}_t$-measurable, we even have this relation for all mappings $\tau_h(\omega)=\omega+\int_0^{\cdot}h_rdr,$ where $\int_0^{\cdot}h_rdr\in H$. With the help of Lemma 3.4 in \cite{BL08} it follows that
$$W(t,x,\zeta,u^2)=E[W(t,x,\zeta,u^2)],\quad \mathbb{P}\text{-a.s.}$$
By identifying $W(t,x,\zeta,u^2)$ with its deterministic version $E[W(t,x,\zeta,u^2)]$, we can regard $W(t,\cdot,$ $\cdot,\cdot):\ \mathbb{R}^n\times L^2(\mathcal{F}_t;\mathbb{R}^n)\times \mathcal{U}^0_{t,T}\rightarrow \mathbb{R}$ as a deterministic function.
\endpf
\vspace{3mm}
The above lemma, combined with the definition of the value function in (\ref{value}), yields
$$
V(t,x,\zeta)=\inf_{u^2 \in \mathcal{U}^0_{t,T}} W\left(t, x, \zeta, u^{2}\right)=\inf_{u^2 \in \mathcal{U}^0_{t,T}}\left(\essinf_{u^1 \in \mathcal{U}^0_{t,T}} J\left(t, x, \zeta,\left(u^{1}, u^{2}\right)\right)\right),\quad \mathbb{P}\text{-a.s.}
$$
Since for any $u^2\in\mathcal{U}^0_{t,T}$, $W\left(t, x, \zeta, u^{2}\right)$ is deterministic, we have, obviously, the following statement:
\begin{corollary}
For all $(t ,x) \in [0,T] \times \mathbb{R}^n, \ \zeta \in L^{2}\left(\mathcal{F}_t;\mathbb{R}^n\right),$ the value function $V(t, x, \zeta)$ is deterministic.
\end{corollary}
Next, we want to prove that the value function $V(t,x,\zeta)$ does not depend on $\zeta$ itself but only on its law $\mathbb{P}_{\zeta}$. For this, we have first to study some auxiliary results. We introduce the filtration $\{\mathcal{F}_s^t\}_{s\in [t,T]}$ generated by the increments of Brownian motion $B$ 
on the intervals after $t$, i.e.,
 $$\mathcal{F}_s^t := \sigma\{B_r-B_t,\ r\in[t,s]\},\ s\in[t,T].$$
\begin{lemma}\label{dense}
Consider the following space of elementary control processes:
\begin{equation}\label{U}
\begin{aligned}
\mathcal{U}^e_{t,T} := \Big\{&\ u^2=\sum_{i,j=0}^{N-1}\mathbf{1}_{A_{i,j}}\zeta_{i,j}\mathbf{1}_{(t_i,t_{i+1}]}\Big|\ N\geq 1,\ t=t_0\le\cdots\le t_N=T,\ \zeta_{i,j}\in L^2(\mathcal{F}_{t_i}^t; U),\\& \ A_{i,j}\in\mathcal{F}_t, \ i=0,\cdots,N-1,\ (A_{i,j})_{j=0}^{N-1}\text{is a decomposition of}\ \ \Omega.\Big\}
\end{aligned}
\end{equation}
Then $\mathcal{U}^e_{t,T}$ is dense in $L^2_{\mathbb{F}}([t,T];U)\big(=\mathcal{U}^0_{t,T}\big)$ with respect to the $L^2$-norm over $[t,T]\times\Omega$.
\end{lemma}
For the proof of this lemma the reader is referred to Peng \cite{Peng19}.
Then standard arguments allow to define $V(t,x,\zeta)$ in the following equivalent way:
$$
V(t,x,\zeta)=\inf_{u^2 \in \mathcal{U}^e_{t,T}} W\left(t, x, \zeta, u^{2}\right).
$$
\begin{proposition}\label{tech}
For any\ $\zeta,\zeta'\in L^2(\mathcal{F}_t;\mathbb{R}^n)$ with $\mathbb{P}_{\zeta}=\mathbb{P}_{\zeta'}$, and for any $u^2\in\mathcal{U}^e_{t,T}$, there exists\ ${u^{2\prime}}\in\mathcal{U}^e_{t,T}$ such that $\mathbb{P}_{(\zeta,u^2,B_{t\vee\cdot}-B_t)}=\mathbb{P}_{({\zeta}',{u^{2\prime}},B_{t\vee\cdot}-B_t)}.$
\end{proposition}
The proof of Proposition \ref{tech} is based on a technical result:
\begin{lemma}\label{4.3}
Let $m,k\geq 1$,\ $\eta\in L^2(\mathcal{F}_t;\mathbb{R}^k)$, and $\zeta,\zeta'\in L^2(\mathcal{F}_t;\mathbb{R}^m)$ 
such that $\mathbb{P}_{\zeta}=\mathbb{P}_{\zeta'}$. Then there exists\ $\eta'\in L^2(\mathcal{F}_t;\mathbb{R}^k)$ such that $\mathbb{P}_{(\zeta,\eta)}= \mathbb{P}_{(\zeta ',\eta ')}.$
\end{lemma}
\noindent\textit{Proof}. We prove this result in three steps.\\
\textbf{Step 1}: We first assume that $\eta\in L^2(\mathcal{F}_t;\mathbb{R})$. Let
\begin{equation*}
F_x(y):=\mathbb{P}\{\eta\leq y\mid \zeta=x\},\quad y\in\mathbb{R},\ x\in\mathbb{R}^m
\end{equation*}
 be the regular conditional cumulative distribution function of $\eta$ knowing $\zeta$, which satisfies:\\
\textbf{(i)}\ For all $x\in\mathbb{R}^m$, $F_x(\cdot)$ is a cumulative distribution function;\\
\textbf{(ii)}\ For all $y\in\mathbb{R}$, $\Gamma\in\mathcal{B}(\mathbb{R}^m)$,
$$
\int_{\Gamma}F_x(y)\mathbb{P}_{\zeta}(dx)=\mathbb{P}\{\eta\leq y,\zeta\in\Gamma\};
$$
\textbf{(iii)}\ The function $(x,y)\mapsto F_x(y)$ is Borel-measurable.\\
We denote by
 $$F^*_x(u):=\inf\{y\in\mathbb{R}|\ F_x(y)\geq u\},\quad u\in (0,1),$$
the left-inverse cumulative distribution function of $F_x$. Let $\overline{u}$ be an $\mathcal{F}_t$-measurable random variable, independent of $\zeta'$, $\zeta$ under $\mathbb{P}$, and uniformly distributed on $[0,1]$. We define
\begin{equation*}
\eta':=F^*_{\zeta'}(\overline{u})\big(:=F_x^*(y)\mid_{(x,y)=(\zeta',\overline u)}\big).
\end{equation*}
Then $\eta'$ is $\mathcal{F}_t$-measurable. 
 For all $(y,\Gamma)\in\mathbb{R}\times\mathcal{B}(\mathbb{R})$, we have the following equalities thanks to the independence between $\overline u$ and $\zeta'$:
\begin{equation}
\begin{aligned}
&\mathbb{P}\{\eta'\leq y,\zeta'\in\Gamma\}
=\int_{\Gamma}\mathbb{P}\{\eta'\leq y\mid \zeta'=x\}\mathbb{P}_{\zeta'}(dx)
=\int_{\Gamma}\mathbb{P}\{F_x^*(\overline u)\leq y\mid \zeta'=x\}\mathbb{P}_{\zeta'}(dx)\\
& =\int_{\Gamma}\mathbb{P}\{F_x^*(\overline u)\leq y\}\mathbb{P}_{\zeta}(dx)=\int_{\Gamma}\mathbb{P}\{F_x^*(\overline u)\leq y\mid \zeta=x\}\mathbb{P}_{\zeta}(dx)\\
&=\int_{\Gamma}\mathbb{P}\{\eta\leq y\mid \zeta=x\}\mathbb{P}_{\zeta}(dx)=\mathbb{P}\{\eta\leq y,\zeta\in\Gamma\}.
\end{aligned}
\end{equation}
Consequently, we obtain $\mathbb{P}_{(\zeta',\eta')}=\mathbb{P}_{(\zeta,\eta)}$.\\
\textbf{Step 2}: Let $\eta=(\eta_1,\eta_2)\in L^2(\mathcal{F}_t;\mathbb{R}^2)$. From Step 1, there exists $\eta_1'\in L^2(\mathcal{F}_t;\mathbb{R})$ such that $\mathbb{P}_{(\zeta',\eta_1')}=\mathbb{P}_{(\zeta,\eta_1)}$. This allows to use Step 1 again, but now for $(\zeta,\eta_1)$ and $(\zeta',\eta_1')$ instead of $\zeta$ and $\zeta'$, respectively, and to conclude that with our $\eta_2\in L^2(\mathcal{F}_t;\mathbb{R})$ we can associate $\eta_2'\in L^2(\mathcal{F}_t;\mathbb{R})$ such that $\mathbb{P}_{((\zeta',\eta_1'),\eta_2')}=\mathbb{P}_{((\zeta,\eta_1),\eta_2)}$, i.e., $\mathbb{P}_{(\zeta',\eta')}=\mathbb{P}_{(\zeta,\eta)}$.\\
\textbf{Step 3}: With an argument analogous to Step 2, we show by iteration that the statement of the lemma holds true.
\endpf
\vspace{3mm}
\noindent\textit{Proof of Proposition \ref{tech}}.

For $u^2 \in \mathcal{U}_{t,T}^e$, from the definition of $\mathcal{U}_{t ,T}^e$ in (\ref{U}), we know that $u^2$ is of the form $\displaystyle u^2=\sum_{i, j=0}^{N-1} \mathbf{1}_{A_{i, j}} \zeta_{i, j} \mathbf{1}_{(t_{i}, t_{i+1}]}$,
where $\left(A_{i, j}\right)_{j=0}^{N-1}$ is the $\left(\Omega, \mathcal{F}_t\right)$-partition, $\zeta_{i, j} \in L^2(\mathcal{F}_{t_i}^t; U)$, for all $1 \leq i \leq N$. Let $\eta:=(\mathbf{1}_{A_{i,j}})_{0 \leq i, j \leq N-1}\in L^2(\mathcal{F}_t; \mathbb{R}^{N \times N})$. From the auxiliary result given by Lemma \ref{4.3}, for $\mathbb{P}_{\zeta}=\mathbb{P}_{\zeta'}$, there exists $\eta^{\prime} \in L^2(\mathcal{F}_{t}; \mathbb{R}^{N \times N})$, such that
\begin{equation}\label{6}
\mathbb{P}_{(\zeta,\eta)}=\mathbb{P}_{(\zeta',\eta')}.
\end{equation}
In particular, as (\ref{6}) implies $\mathbb{P}_{\eta}=\mathbb{P}_{\eta'}$, it follows that
$$
\eta^{\prime}=(\mathbf{1}_{A_{i, j}^{\prime}})_{0 \leq i, j \leq N-1} \text {, }
$$
where also $(A_{i, j}^{\prime})_{j=0}^{N-1}$ is an $(\Omega, \mathcal{F}_t)$-partition, for $0 \leq i \leq N-1$. Then, obviously,
\begin{equation}\label{u2'}
u^{2\prime}=\sum_{i, j=0}^{N-1} \mathbf{1}_{A_{i, j}^{\prime}} \zeta_{i, j} \mathbf{1}_{(t_i, t_{i+1}]}\in\mathcal{U}^e_{t,T}\\
\end{equation}
is such that $$\mathbb{P}_{(\zeta^{\prime}, u^{2\prime}, B_{\cdot \vee t}-B_t)}=\mathbb{P}_{(\zeta, u^2, B_{\cdot\vee t}-B_t)}.$$
\endpf
\vspace{3mm}
The following result is a direct consequence of Proposition \ref{tech}.
\begin{corollary}
Given any $\zeta,\zeta^{\prime} \in L^{2}(\mathcal{F}_{t})$ of the same law $\mathbb{P}_{\zeta}= \mathbb{P}_{\zeta' }$, and any $u^{2} \in \mathcal{U}_{t, T}^{e}$ (of form (\ref{U})), we see that, for $u^{2\prime}$ of form (\ref{u2'}),
$$\mathbb{P}\circ \big[(X^{t,{\zeta}',{u^{2\prime}}},{u^{2\prime}})\big]^{-1} =  \mathbb{P}\circ \big[(X^{t,\zeta,u^2},{u^2})\big]^{-1}.$$
\end{corollary}
Thus, replacing $\mathbb{P}_{(X_r^{t,\zeta,{u^{2}}},{u_r^{2}})}$ by $\mathbb{P}_{(X_r^{t,{\zeta}',{u^{2\prime}}},{u_r^{2\prime}})}$ in SDE (\ref{SDE2}), $r \in[t, T]$, we get from the uniqueness of the solution:
$$\quad \quad \quad \quad X^{t,x,{\zeta}',(u^1,{u^{2\prime}})} = X^{t,x,\zeta,(u^1,u^2)},\quad  \mathbb{P}\text{-a.s.},\ 
\ u^1\in\mathcal{U}^0_{t,T}.$$
Consequently, $J\left(t, x, \zeta^{\prime},\left(u^{1}, {u^{2\prime}}\right)\right)=J\big(t, x, \zeta,\big(u^{1}, u^{2}\big)\big),\ \mathbb{P}\text{-a.s.}$, for all $u^{1} \in \mathcal{U}^0_{t,T}$, i.e., we have shown the following result:
\begin{lemma}
 For every $\zeta,\zeta^{\prime} \in L^{2}\left(\mathcal{F}_{t};\mathbb{R}^n\right)$ 
 with $\mathbb{P}_{\zeta}=\mathbb{P}_{\zeta'}$, and every $u^{2} \in \mathcal{U}^e_{t, T}$, there exists ${u^{2\prime}} \in \mathcal{U}^e_{t, T}$ such that
$$J\big(t, x, \zeta^{\prime},\big(u^{1}, {u^{2\prime}}\big)\big)=J\left(t, x, \zeta,\left(u^{1}, u^{2}\right)\right),\quad \mathbb{P}\text{-a.s., for all } u^{1} \in \mathcal{U}^0_{t,T},\ x\in\mathbb{R}^n,$$
and, in particular,
\begin{equation}\label{8}
W\big(t, x, \zeta ', {u^{2\prime}}\big):=\essinf_{u^1 \in \mathcal{U}^0_{t,T}} J\big(t, x, \zeta ',\big(u^{1}, {u^{2\prime}}\big)\big)=\essinf_{u^1 \in \mathcal{U}^0_{t,T}} J\big(t, x, \zeta,\big(u^{1}, u^{2}\big)\big)=W\big(t, x, \zeta, u^{2}\big).
\end{equation}
\end{lemma}
From (\ref{8}) in the above lemma we deduce now easily:
\begin{proposition}
Let\ $(t,x)\in [0,T]\times \mathbb{R}^n$, $\zeta\in L^2(\mathcal{F}_t;\mathbb{R}^n)$. Then, for all\ ${\zeta}'\in L^2(\mathcal{F}_t;\mathbb{R}^n)$\ with\ $\mathbb{P}_{\zeta}=\mathbb{P}_{\zeta '}$, we have $$V(t,x,\zeta)=V(t,x,\zeta '),$$ i.e., $V$ depends on $\zeta$ only through $\mathbb{P}_{\zeta}$. We write:
$$V(t,x,\mathbb{P}_{\zeta}):=V(t,x,\zeta),$$
where\ $V: [0,T]\times\mathbb{R}^n\times\mathcal{P}_2(\mathbb{R}^n)\rightarrow \mathbb{R}$. Moreover, by standard estimates, there exists\ $C\in\mathbb{R}_{+}$, such that
\begin{equation}\label{estimate_V}
\big|V(t,x,\mathbb{P}_{\zeta})-V(t,\widetilde{x},\mathbb{P}_{\widetilde{\zeta}})\big|\leq C\big(|x-\widetilde{x}|+\mathcal{W}_2(\mathbb{P}_{\zeta},\mathbb{P}_{\widetilde{\zeta}})\big),
\end{equation}
for all\ $t\in[0,T],\ x,\widetilde{x}\in \mathbb{R}^n,\ \zeta,\widetilde{{\zeta}'}\in L^2(\mathcal{F}_t;\mathbb{R}^n).$
\end{proposition}

\vspace{3mm}
Our next objective is to prove the dynamic programming principle (DPP) for $W$, and 
 to use this to study the properties of the value function $V$. Recall that $W(t,\cdot,\cdot,\cdot):\ \mathbb{R}\times L^2(\mathcal{F}_t;\mathbb{R}^n)\times\mathcal{U}^0_{t,T}\rightarrow \mathbb{R}$, is defined by
$$ W\left(t, x, \zeta, u^{2}\right)=\essinf_{u^1 \in \mathcal{U}^0_{t,T}} J\left(t, x, \zeta,\left(u^{1}, u^{2}\right)\right)
=\essinf_{u^1 \in \mathcal{U}^0_{t,T}}E \big[ \Phi \big(X_{T}^{t, x, \zeta, (u^{1}, u^{2})},\mathbb{P}_{X_{T}^{t, \zeta, u^2}}\big)\big|\mathcal{F}_t \big].$$
Lemma \ref{W_deter} shows that $W(t,x,\zeta,u^2)$, $u^2\in \mathcal{U}^0_{t,T}$, is deterministic and there is some constant $C\in\mathbb{R}_+$ such that
\begin{equation}\label{W_est}
|W(t,x,\zeta,u^2)-W(t,x',\zeta',u^{2})|\leq C \left(\left|x-x^{\prime}\right|+\big(E[\left|\zeta-\zeta^{\prime}\right|^2]\big)^{\frac{1}{2}}\right),
\end{equation}
for all $t\in [0,T]$, $x,x'\in\mathbb{R}^n$, $\zeta,\zeta'\in L^2(\mathcal{F}_t;\mathbb{R}^n)$, $u^2\in\mathcal{U}^0_{t,T}$.
Using the above properties of $W$, we obtain the following DPP for $W$.
\begin{theorem}[DPP for $W$]\label{DPP_W}
For all\hspace{1mm} $0\leq t< t+\delta\leq T,\ u^2\in\mathcal{U}^0_{t,T},\ x\in\mathbb{R}^n,\ \zeta\in L^2(\mathcal{F}_t;\mathbb{R}^n),$
$$W(t,x,\zeta,u^2)=\essinf_{u^1\in \mathcal{U}^0_{t,t+\delta}}E\big[W\big(t+\delta,X_{t+\delta}^{t,x,\zeta,(u^1,u^2)},X_{t+\delta}^{t,\zeta,u^2},u^2\big)\big|\mathcal{F}_t\big].$$
\end{theorem}
\noindent\textit{Proof}. We prove this by using a standard argument from Buckdahn and Li  \cite{BL08}. For $(t, x) \in[0, T] \times \mathbb{R}^n, \ \zeta \in L^2\left(\mathcal{F}_t\right)$, we fix $u^2 \in \mathcal{U}_{t+\delta , T}^0$. Then there exists a sequence $\left\{u^{1, k}\right\}_{k \geq 1} \subset \mathcal{U}_{t, t+\delta}^0$, such that
\begin{equation}
\begin{aligned}
\widehat{W}_\delta\left(t, x, \zeta, u^2\right):=&\essinf_{u^{1} \in \mathcal{U}_{t,t+\delta}^0}E\big[W\big(t+\delta,X_{t+\delta}^{t,x,\zeta,(u^1,u^2)},X_{t+\delta}^{t,\zeta,u^2},u^2\big)\big|\mathcal{F}_t\big]\\
=&\inf_{k\geq 1}E\big[W\big(t+\delta,X_{t+\delta}^{t,x,\zeta,(u^{1,k},u^2)},X_{t+\delta}^{t,\zeta,u^2},u^2\big)\big|\mathcal{F}_t\big], \quad \mathbb{P}\text{-a.s.}
\end{aligned}
\end{equation}
Throughout the proof we put
$$
I_\delta(t, x, \zeta,(u^1, u^2)):=E\big[W\big(t+\delta,X_{t+\delta}^{t,x,\zeta,(u^1,u^2)},X_{t+\delta}^{t,\zeta,u^2},u^2\big)\big|\mathcal{F}_t\big].
$$
For any $\varepsilon>0$, we set $\Gamma_k:=\{I_\delta(t, x, \zeta,(u^{1,k}, u^2))\leq \widehat{W}_\delta\left(t, x, \zeta, u^2\right)+\varepsilon \}\in\mathcal{F}_t$, $k\geq 1$. Then $\displaystyle\mathbb{P}(\cup_{k\geq 1}\Gamma_k)=1$, the family of sets $\displaystyle\Delta_1:=\Gamma_1,\ \Delta_k:=\Gamma_k\backslash (\cup_{\ell\leq k-1}\Gamma_{\ell})\in\mathcal{F}_t,\ k\geq 1$ forms an $(\Omega,\mathcal{F}_t)$-decomposition, and $u^{1,\varepsilon}:=\sum_{k \geq 1} \mathbf{1}_{\Delta_k} u^{1,k}$ belongs obviously to $\mathcal{U}^0_{t, t+\delta}$. Moreover, from the uniqueness of the solution of (\ref{SDE2}), we deduce that
\begin{equation*}
X_{s}^{t,x,\zeta,(u^{1,\varepsilon},u^2)}=\sum_{k\geq 1}\mathbf{1}_{\Delta_k}X_{s}^{t,x,\zeta,(u^{1,k},u^2)},\quad s\in[t,T],\ \mathbb{P}\text{-a.s.}
\end{equation*}
Therefore,
\begin{equation}\label{11}
\begin{aligned}
&\widehat{W}_{\delta}(t,x,\zeta,u^2)+\varepsilon \\&\geq  \sum_{k\geq 1}\mathbf{1}_{\Delta_k}I_\delta(t, x, \zeta,(u^{1,k}, u^2))=E\big[\sum_{k\geq 1}\mathbf{1}_{\Delta_k}W\big(t+\delta,X_{t+\delta}^{t,x,\zeta,(u^{1,k},u^2)},X_{t+\delta}^{t,\zeta,u^2},u^2\big)\big|\mathcal{F}_t\big]\\&=E\big[W\big(t+\delta,X_{t+\delta}^{t,x,\zeta,(u^{1,\varepsilon},u^2)},X_{t+\delta}^{t,\zeta,u^2},u^2\big)\big|\mathcal{F}_t\big]=I_\delta(t, x, \zeta,(u^{1,\varepsilon}, u^2)), \quad \mathbb{P}\text{-a.s.}
\end{aligned}
\end{equation}
Let $\left\{O_i\right\}_{i \geq 1} \subset \mathcal{B}\left(\mathbb{R}^n\right)$ be a partition of $\mathbb{R}^n$ such that $\sum_{i \geq 1} O_i=\mathbb{R}^n$ and $\operatorname{diam}\left(O_i\right)$ $\leq \varepsilon,\ i \geq 1$. Let $y_i$ be an arbitrarily fixed element of $O_i,\ i \geq 1$.
As \begin{equation*}\begin{aligned}
W\big(t+\delta,y_i,X_{t+\delta}^{t,\zeta,u^2},u^2\big)&=\essinf_{\widetilde{u}^1\in\mathcal{U}^0_{t+\delta,T}}
E\big[\Phi\big(X_{T}^{t+\delta,y_i,X_{t+\delta}^{t,\zeta,u^2},(\widetilde{u}^1,u^2)},\mathbb{P}_{X_{T}^{t+\delta,X_{t+\delta}^{t,\zeta,u^2},u^2}}\big)\big|\mathcal{F}_{t+\delta}\big]
\\&=\essinf_{\widetilde{u}^1\in\mathcal{U}^0_{t+\delta,T}}
E\big[\Phi\big(X_{T}^{t+\delta,y_i,X_{t+\delta}^{t,\zeta,u^2},(\widetilde{u}^1,u^2)},\mathbb{P}_{X_{T}^{t,\zeta,u^2}}\big)\big|\mathcal{F}_{t+\delta}\big]
,\end{aligned}\end{equation*}
we see in analogy to the first part of the proof that, for every $i\geq 1$, there exist a sequence $\left\{\widetilde{u}^{1, i}\right\}_{i \geq 1} \subset \mathcal{U}_{t+\delta,T}^0$, such that
\begin{equation*}
W\big(t+\delta,y_i,X_{t+\delta}^{t,\zeta,u^2},u^2\big)\geq E\big[\Phi\big(X_{T}^{t+\delta,y_i,X_{t+\delta}^{t,\zeta,u^2},(\widetilde{u}^{1,i},u^2)},\mathbb{P}_{X_{T}^{t,\zeta,u^2}}\big)\big|\mathcal{F}_{t+\delta}\big]-\varepsilon,
\end{equation*}
and, thus,
\begin{equation}\label{4.19}
\begin{aligned}
W\big(t+\delta,y_i,X_{t+\delta}^{t,\zeta,u^2},u^2\big)&\geq E\big[\Phi\big(X_{T}^{t+\delta,X_{t+\delta}^{t,x,\zeta,(u^{1,\varepsilon},u^2)},X_{t+\delta}^{t,\zeta,u^2},(\widetilde{u}^{1,i},u^2)},\mathbb{P}_{X_{T}^{t,\zeta,u^2}}\big)\big|\mathcal{F}_{t+\delta}\big]
\\&\quad -C\big|X_{t+\delta}^{t,x,\zeta,(u^{1,\varepsilon},u^2)}-y_i\big|-\varepsilon,
\end{aligned}
\end{equation}
where the latter estimate follows from (\ref{Es_SDE}). Hence, for $$\widehat{u}^{1,\varepsilon}:=u^{1,\varepsilon}\mathbf{1}_{[t,t+\delta)}+\sum_{i\geq 1}\mathbf{1}_{\big\{X_{t+\delta}^{t,x,\zeta,(u^{1,\varepsilon},u^2)}\in O_i\big\}}\widetilde{u}^{1,i}\mathbf{1}_{[t+\delta,T]}\in\mathcal{U}_{t,T}^0,$$ by (\ref{W_est}) and (\ref{4.19}) we have
\begin{equation}
\begin{aligned}
&W\big(t+\delta,X_{t+\delta}^{t,x,\zeta,(u^{1,\varepsilon},u^2)},X_{t+\delta}^{t,\zeta,u^2},u^2\big)\\&\geq \sum_{i\geq 1} W\big(t+\delta,y_i,X_{t+\delta}^{t,\zeta,u^2},u^2\big)\mathbf{1}_{\big\{X_{t+\delta}^{t,x,\zeta,(u^{1,\varepsilon},u^2)}\in O_i\big\}}-C\varepsilon\\
&\geq \sum_{i\geq 1}\Big( E\big[\Phi\big(X_{T}^{t+\delta,X_{t+\delta}^{t,x,\zeta,(u^{1,\varepsilon},u^2)},X_{t+\delta}^{t,\zeta,u^2},(\widetilde{u}^{1,i},u^2)},\mathbb{P}_{X_{T}^{t,\zeta,u^2}}\big)\big|\mathcal{F}_{t+\delta}\big]\\
&\hspace{1.5cm} -(C+1)\varepsilon\Big)\mathbf{1}_{\big\{X_{t+\delta}^{t,x,\zeta,(u^{1,\varepsilon},u^2)}\in O_i\big\}}-C\varepsilon
\\&=E\big[\Phi\big(X_{T}^{t,x,\zeta,(\widehat{u}^{1,\varepsilon},u^2)},\mathbb{P}_{X_{T}^{t,\zeta,u^2}}\big)\big|\mathcal{F}_{t+\delta}\big]-(2C+1)\varepsilon
,\quad \mathbb{P}\text{-a.s.}\end{aligned}
\end{equation}
Thus,
\begin{equation}
\begin{aligned}
&E\big[W\big(t+\delta,X_{t+\delta}^{t,x,\zeta,(u^{1,\varepsilon},u^2)},X_{t+\delta}^{t,\zeta,u^2},u^2\big)\big|\mathcal{F}_{t}\big]\geq E\big[\Phi\big(X_{T}^{t,x,\zeta,(\widehat{u}^{1,\varepsilon},u^2)},\mathbb{P}_{X_{T}^{t,\zeta,u^2}}\big)\big|\mathcal{F}_{t}\big]-(2C+1)\varepsilon
\\&\geq
\essinf_{{u}^1\in\mathcal{U}^0_{t,T}}E\big[\Phi\big(X_{T}^{t,x,\zeta,(u^{1},u^2)},\mathbb{P}_{X_{T}^{t,\zeta,u^2}}\big)\big|\mathcal{F}_{t}\big]-(2C+1)\varepsilon
=W(t,x,\zeta,u^2)-(2C+1)\varepsilon
,\quad \mathbb{P}\text{-a.s.}\end{aligned}
\end{equation}
Consequently, combined with (\ref{11}) we have
\begin{equation}\label{12}
\widehat{W}_{\delta}(t,x,\zeta,u^2)\geq W(t,x,\zeta,u^2)-(2C+1)\varepsilon-\varepsilon, \quad \varepsilon>0.
\end{equation}
Letting $\varepsilon\rightarrow0$, this yields
\begin{equation}\label{13}
\widehat{W}_{\delta}(t,x,\zeta,u^2)\geq W(t,x,\zeta,u^2).
\end{equation}

On the other hand, for all $u^1\in\mathcal{U}^0_{t,T}$:
\begin{equation}\label{3.23}
\begin{aligned}
&\widehat{W}_\delta\left(t, x, \zeta, u^2\right)\\&\leq E\big[W\big(t+\delta,X_{t+\delta}^{t,x,\zeta,(u^1,u^2)},X_{t+\delta}^{t,\zeta,u^2},u^2\big)\big|\mathcal{F}_t\big]=E\big[W\big(t+\delta,y,X_{t+\delta}^{t,\zeta,u^2},u^2\big)\big|_{y=X_{t+\delta}^{t,x,\zeta,(u^1,u^2)}}\big|\mathcal{F}_t\big]
\\&\leq E\big[E\big[\Phi\big(X_{T}^{t+\delta,y,X_{t+\delta}^{t,\zeta,u^2},({u}^{1},u^2)},\mathbb{P}_{X_{T}^{t+\delta,X_{t+\delta}^{t,\zeta,u^2},u^2}}\big)\big|\mathcal{F}_{t+\delta}\big]\big|_{y=X_{t+\delta}^{t,x,\zeta,(u^1,u^2)}}\big|\mathcal{F}_{t}\big]
\\&=E\big[E\big[\Phi\big(X_{T}^{t,x,\zeta,({u}^{1},u^2)},\mathbb{P}_{X_{T}^{t+\delta,X_{t+\delta}^{t,\zeta,u^2},u^2}}\big)\big|\mathcal{F}_{t+\delta}\big]\big|\mathcal{F}_{t}\big]
=E\big[\Phi\big(X_{T}^{t,x,\zeta,({u}^{1},u^2)},\mathbb{P}_{X_{T}^{t,\zeta,u^2}}\big)\big|\mathcal{F}_{t}\big]
, \quad \mathbb{P}\text{-a.s.,}
\end{aligned}
\end{equation}
and hence,
\begin{equation}
\widehat{W}_{\delta}(t,x,\zeta,u^2)\leq W(t,x,\zeta,u^2).
\end{equation}

Together with (\ref{12}), this proves the DPP for $W(\cdot,\cdot,\cdot,u^2),$ i.e., for all $u^2\in\mathcal{U}^0_{t,T}$:
\begin{equation}
\widehat{W}_{\delta}(t,x,\zeta,u^2)= W(t,x,\zeta,u^2).
\end{equation}
\endpf
\vspace{3mm}
\begin{remark}
For all $u^2 \in \mathcal{U}_{t, T}^0,\ \varepsilon>0$, from the inequalities (\ref{11}) and (\ref{13})-(\ref{3.23}) we obtain there exists $u^{1,\varepsilon}\in\mathcal{U}^0_{t,t+\delta}$ such that
 $$W(t, x, \zeta, u^2)=E\left[\widehat{W}_\delta(t, x, \zeta, u^2)\right] \geq E\left[W\big(t+\delta,X_{t+\delta}^{t,x,\zeta,(u^{1,\varepsilon},u^2)},X_{t+\delta}^{t,\zeta,u^2},u^2\big)\right]-\varepsilon.$$
On the other hand,
$$W(t, x, \zeta, u^2)\leq E\left[W\big(t+\delta,X_{t+\delta}^{t,x,\zeta,(u^{1},u^2)},X_{t+\delta}^{t,\zeta,u^2},u^2\big)\right] ,\quad
\text{for all}\ u^1 \in \mathcal{U}_{t, t+\delta}^0.$$
Consequently,
\begin{equation}\label{W_remark}
W(t, x, \zeta, u^2)=\inf_{{u}^1\in\mathcal{U}^0_{t,t+\delta}}E\left[W\big(t+\delta,X_{t+\delta}^{t,x,\zeta,(u^{1},u^2)},X_{t+\delta}^{t,\zeta,u^2},u^2\big)\right],\  t<t+\delta \leq T.
\end{equation}
In particular, for $t+\delta=T$,
$$W(t, x, \zeta, u^2)=\inf_{{u}^1\in\mathcal{U}^0_{t,T}}E\left[\Phi\big(X_{T}^{t,x,\zeta,(u^{1},u^2)},\mathbb{P}_{X_{T}^{t,\zeta,u^2}}\big)\right],$$
and therefore, for the value function $V(t,x,\mathbb{P}_{\zeta})=V(t,x,\zeta)$\rm{:}
\begin{equation}
\begin{aligned}
&V(t,x,\mathbb{P}_{\zeta})=V(t,x,\zeta)=\essinf_{({u}^1,u^2)\in\mathcal{U}_{t,T}}E\big[\Phi\big(X_{T}^{t,x,\zeta,(u^{1},u^2)},\mathbb{P}_{X_{T}^{t,\zeta,u^2}}\big)\big|\mathcal{F}_{t}\big]
\\&=\inf_{u^2\in\mathcal{U}_{t,T}^0}W(t,x,\zeta,u^2)=\inf_{({u}^1,u^2)\in\mathcal{U}_{t,T}}E\big[\Phi\big(X_{T}^{t,x,\zeta,(u^{1},u^2)},\mathbb{P}_{X_{T}^{t,\zeta,u^2}}\big)\big].
\end{aligned}
\end{equation}
\end{remark}

Now we use the properties of $W$ we have proved above to study $V$. As the following inequality shows, we can only get the one-sided DPP for $V$: For $0\leq t<t+\delta\leq T$, $x\in\mathbb{R}^n$, $\zeta\in L^2(\mathcal{F}_t;\mathbb{R}^n)$, thanks to equality (\ref{W_remark}) and the fact that for any $u^2\in\mathcal{U}^0_{t,T}$, $W(t,x,\zeta,u^2)$ is deterministic, we get
\begin{equation}\label{DPPNO}
\begin{aligned}
V(t,x,\mathbb{P}_{\zeta})&=\inf_{u^2\in\mathcal{U}_{t,T}^0}W(t,x,\zeta,u^2)\\
&=\inf_{u^2\in\mathcal{U}_{t,T}^0}\Big(\inf_{u^1\in\mathcal{U}_{t,t+\delta}^0}E\big[W\big(t+\delta,X_{t+\delta}^{t,x,\zeta,(u^{1},u^2)},X_{t+\delta}^{t,\zeta,u^2},u^2\big)\big]\Big)
\\&\geq \inf_{u^2\in\mathcal{U}_{t,T}^0}\Big(\inf_{u^1\in\mathcal{U}_{t,t+\delta}^0}E\big[V\big(t+\delta,X_{t+\delta}^{t,x,\zeta,(u^{1},u^2)},X_{t+\delta}^{t,\zeta,u^2}\big)\big]\Big)
\\&=\inf_{(u^1,u^2)\in\mathcal{U}_{t,t+\delta}}E\big[V\big(t+\delta,X_{t+\delta}^{t,x,\zeta,(u^{1},u^2)},\mathbb{P}_{X_{t+\delta}^{t,\zeta,u^2}}\big)\big].
\end{aligned}
\end{equation}
But we can not get the above inequality in the opposite direction.

To deal with this problem, we introduce a novel definition of the value function, denoted by $\vartheta$, which is closely related to the original value function $V$. As we will see below, the DPP for $\vartheta$ holds, and we can thus derive the description of $\vartheta$ as a solution of a PDE. Moreover, Remark \ref{Characterize} shows that, this PDE can also characterize the original value function $V$.\\

\noindent\textbf{Novel definition of value function\ $\vartheta$.}\quad For $\theta,\zeta\in L^2(\mathcal{F}_t;\mathbb{R}^n)$, we introduce a new definition of the value function:
\begin{equation}\label{NewValue}
\begin{aligned}
\vartheta(t,\theta,\mathbb{P}_{\zeta}):=&\ \inf_{u^2\in \mathcal{U}^0_{t,T}}E \big[W(t,\theta,\zeta,u^2)\big]\\
=&\ \inf_{u^2\in \mathcal{U}^0_{t,T}}E \big[\mathop{\essinf}\limits_{u^1\in \mathcal{U}^0_{t,T}}E \big[ \Phi \big(X_{T}^{t, \theta, \zeta, (u^1,u^2)},\mathbb{P}_{X_{T}^{t, \zeta, u^2}}\big)\big|\mathcal{F}_t \big]\big].
\end{aligned}
\end{equation}
\noindent The function $\vartheta$ is obviously deterministic, and we can prove that $\vartheta(t,\theta,\mathbb{P}_{\zeta})$ depends on $\theta$ only through the law $\mathbb{P}_{\theta}$. Indeed, since $W(t,x,\zeta,u^2)$ is deterministic, also $\vartheta(t,\theta,\mathbb{P}_{\zeta})$ is deterministic and
 $$\vartheta(t,\theta,\mathbb{P}_{\zeta})=\inf_{u^2\in \mathcal{U}^0_{t,T}}\int_{\mathbb{R}^n}W(t,x,\zeta,u^2)\mathbb{P}_{\theta}(dx).$$
 So $\vartheta(t,\theta,\mathbb{P}_{\zeta})$ depends on $\theta$ only through its law, which allows to write
  $$\vartheta(t, \mathbb{P}_{\theta},\mathbb{P}_{\zeta}) := \vartheta(t,\theta,\zeta),$$
  and to see $\vartheta$ as a function over $[0,T]\times\mathcal{P}_2(\mathbb{R}^n)\times\mathcal{P}_2(\mathbb{R}^n)$. From Lemma \ref{W_deter} and a   standard argument we see that there exists a constant $C>0$ such that
 \begin{equation}\label{vartheta_con}
 |\vartheta(t, \mathbb{P}_{\theta},\mathbb{P}_{\zeta})-\vartheta(t, \mathbb{P}_{\theta'},\mathbb{P}_{\zeta'})|\leq C\big(\mathcal{W}_2(\mathbb{P}_{\theta},\mathbb{P}_{\theta'})+\mathcal{W}_2(\mathbb{P}_{\zeta},\mathbb{P}_{\zeta'})\big),
 \end{equation}
 for all $t\in [0,T]$, $\theta,\theta',\zeta,\zeta'\in L^2(\mathcal{F}_t;\mathbb{R}^n)$.
 \begin{remark}\label{Characterize}
  Notice that
$$
\vartheta(t,x,\mathbb{P}_{\zeta}) (= \vartheta(t,\delta_x,\mathbb{P}_{\zeta})) = V(t,x,\mathbb{P}_{\zeta}),\ \ (t,x)\in [0,T]\times\mathbb{R}^n, \ \zeta\in L^2(\mathcal{F}_t;\mathbb{R}^n),
$$
\noindent where $\delta_x$ denotes the Dirac measure at $x$, which means that a description of $\vartheta(t,\mathbb{P}_{\theta},\mathbb{P}_{\zeta}),\ (\theta,\zeta)\in L^2(\mathcal{F}_t;\mathbb{R}^n)\times L^2(\mathcal{F}_t;\mathbb{R}^n)$, as a solution of a PDE also characterizes $V(t,x,\mathbb{P}_{\zeta})$.
\end{remark}
Let us now study the new value function $\vartheta$. To begin with, we prove that $\vartheta$ obeys the following $\mathrm{DPP}$:
\begin{theorem}[DPP for $\vartheta$]\label{DPP2}
For any\ $0\leq t< t+\delta\leq T,\ \theta,\zeta\in L^2(\mathcal{F}_t;\mathbb{R}^n),$
\begin{equation}\label{DPP_2_Eq}
\vartheta(t,\mathbb{P}_{\theta},\mathbb{P}_{\zeta})= \inf_{u\in \mathcal{U}_{t,t+\delta}} \vartheta \big(t+\delta, \mathbb{P}_{X_{t+\delta}^{t, \theta, \zeta, u}},\mathbb{P}_{X_{t+\delta}^{t, \zeta, u^2}}\big).
\end{equation}
\end{theorem}
\noindent \textit{Proof}.
Recall from Lemma \ref{W_deter} and Theorem \ref{DPP_W} that, for $x \in \mathbb{R}^n,\ \delta>0,\ u^2\in\mathcal{U}^0_{t,T},$ $W$ satisfies the DPP
\begin{equation}\label{14}
W(t,x,\zeta,u^2)=\essinf_{u^1\in \mathcal{U}^0_{t,t+\delta}}E\big[W\big(t+\delta,X_{t+\delta}^{t,x,\zeta,(u^1,u^2)},X_{t+\delta}^{t,\zeta,u^2},u^2\big)\big|\mathcal{F}_t\big],\end{equation}
and is deterministic and Lipschitz with respect to $x \in \mathbb{R}^n$, uniformly in $(t,\zeta,u^2)$.
So standard arguments allow to substitute in (\ref{14}) for $x$ the random variable $\theta \in L^2(\mathcal{F} _t)$ :
\begin{equation}\label{15}
W(t,\theta,\zeta,u^2)=\essinf_{u^1\in \mathcal{U}^0_{t,t+\delta}}E\big[W\big(t+\delta,X_{t+\delta}^{t,\theta,\zeta,(u^1,u^2)},X_{t+\delta}^{t,\zeta,u^2},u^2\big)\big|\mathcal{F}_t\big], \quad \mathbb{P}\text{-a.s.}\end{equation}
Again, from standard arguments (see, for instance, the proof of Theorem \ref{DPP_W} or Buckdahn and Li \cite{BL08}), we know that, for $\varepsilon>0$, there exists $u^{1,\varepsilon}\in\mathcal{U}^0_{t,t+\delta}$ (depending on $u^2$) such that
\begin{equation}\begin{aligned}
W(t, \theta, \zeta, u^2)\geq E\left[W\big(t+\delta,X_{t+\delta}^{t,\theta,\zeta,(u^{1,\varepsilon},u^2)},X_{t+\delta}^{t,\zeta,u^2},u^2\big)\big|\mathcal{F}_t\right]-\varepsilon,\quad \mathbb{P}\text{-a.s.}
\end{aligned}\end{equation}
Hence, $\text{for all }u^2\in\mathcal{U}_{t,T}^0$,
\begin{equation}\label{*4}\begin{aligned}
E\big[W(t, \theta, \zeta, u^2)\big]&\geq E\left[W\big(t+\delta,X_{t+\delta}^{t,\theta,\zeta,(u^{1,\varepsilon},u^2)},X_{t+\delta}^{t,\zeta,u^2},u^2\big)\right]-\varepsilon\\
&\geq \inf_{u^1\in \mathcal{U}^0_{t,t+\delta}}E\left[W\big(t+\delta,X_{t+\delta}^{t,\theta,\zeta,(u^{1},u^2)},X_{t+\delta}^{t,\zeta,u^2},u^2\big)\right]-\varepsilon.
\end{aligned}\end{equation}

On one hand, identifying $u^2$ with $(u^{2\prime},u^{2\prime\prime})\in\mathcal{U}^0_{t,t+\delta}\times\mathcal{U}^0_{t+\delta,T}$, where $u^{2\prime}:=u^2\big|_{[t,t+\delta]}$, $u^{2\prime\prime}:=u^2\big|_{[t+\delta,T]}$,
letting $\varepsilon\rightarrow 0$ in (\ref{*4}) and taking infimum over $u^2\in \mathcal{U}^0_{t,T}$ yields
$$\begin{aligned}\vartheta\left(t, \mathbb{P}_\theta, \mathbb{P}_{\zeta}\right)\left(=\vartheta\left(t,\theta, \mathbb{P}_{\zeta}\right)\right)&=\inf_{u^2\in \mathcal{U}^0_{t,T}} E\left[W\left(t, \theta, \zeta, u^2\right)\right]\\&\geq\inf_{u=\left(u^1, u^{2\prime}\right) \in \mathcal{U}_{t, t+\delta}} \Big(\inf_{u^{2\prime\prime}\in \mathcal{U}^0_{t,T}}E\big[W\big(t+\delta,X_{t+\delta}^{t,\theta,\zeta,(u^{1},u^{2\prime})},X_{t+\delta}^{t,\zeta,u^{2\prime}},u^{2\prime\prime}\big)\big]\Big)\\&=\inf_{u=\left(u^1, u^{2\prime}\right) \in \mathcal{U}_{t, t+\delta}} \vartheta\big(t+\delta,\mathbb{P}_{X_{t+\delta}^{t,\theta,\zeta,(u^{1},u^{2\prime})}},\mathbb{P}_{X_{t+\delta}^{t,\zeta,u^{2\prime}}}\big).\end{aligned}$$
Hence,
$$\begin{aligned}\vartheta\left(t, \mathbb{P}_\theta, \mathbb{P}_{\zeta}\right)\geq\inf_{u=\left(u^1, u^{2}\right) \in \mathcal{U}_{t, t+\delta}} \vartheta\big(t+\delta,\mathbb{P}_{X_{t+\delta}^{t,\theta,\zeta,(u^{1},u^2)}},\mathbb{P}_{X_{t+\delta}^{t,\zeta,u^{2}}}\big).\end{aligned}$$

On the other hand, for all $u=(u^{1},u^2)\in\mathcal{U}_{t,T}$, it follows from (\ref{15}) that
\begin{equation*}
W(t,\theta,\zeta,u^2)\leq E\big[W\big(t+\delta,X_{t+\delta}^{t,\theta,\zeta,(u^1,u^2)},X_{t+\delta}^{t,\zeta,u^2},u^2\big)\big|\mathcal{F}_t\big],\ \ \mathbb{P}\text{-a.s.,}
\end{equation*}
and, so
\begin{equation*}
E\big[W(t,\theta,\zeta,u^2)\big]\leq E\big[W\big(t+\delta,X_{t+\delta}^{t,\theta,\zeta,(u^1,u^2)},X_{t+\delta}^{t,\zeta,u^2},u^2\big)\big].
\end{equation*}
Thus, we have
$$\begin{aligned}\vartheta\left(t, \mathbb{P}_\theta, \mathbb{P}_{\zeta}\right)&=\inf_{u^2\in \mathcal{U}^0_{t,T}} E\left[W\left(t, \theta, \zeta, u^2\right)\right]
\\&\leq\inf_{u=\left(u^1, u^{2}\right) \in \mathcal{U}_{t, t+\delta}} \Big(\inf_{\widetilde{u}^{2}\in \mathcal{U}_{t+\delta,T}}E\big[W\big(t+\delta,X_{t+\delta}^{t,\theta,\zeta,(u^{1},u^2)},X_{t+\delta}^{t,\zeta,u^{2}},\widetilde{u}^{2}\big)\big]\Big)
\\&=\inf_{u=\left(u^1, u^{2}\right) \in \mathcal{U}_{t, t+\delta}} \vartheta\big(t+\delta,\mathbb{P}_{X_{t+\delta}^{t,\theta,\zeta,u}},\mathbb{P}_{X_{t+\delta}^{t,\zeta,u^{2}}}\big).\end{aligned}$$

Consequently, we obtain the $\mathrm{DPP}$ for $\vartheta$: For all $0\leq t\leq t+\delta\leq T$, $\theta,\zeta\in L^2(\mathcal{F}_t;\mathbb{R}^n)$,
\begin{equation}\label{DPP_2}\begin{aligned}\vartheta\left(t, \mathbb{P}_\theta, \mathbb{P}_{\zeta}\right)=\inf_{u=\left(u^1, u^{2}\right) \in \mathcal{U}_{t, t+\delta}} \vartheta \big(t+\delta,\mathbb{P}_{X_{t+\delta}^{t,\theta,\zeta,u}},\mathbb{P}_{X_{t+\delta}^{t,\zeta,u^{2}}}\big).\end{aligned}\end{equation}
\endpf
 With the help of Theorem \ref{DPP2}, we can now also study the continuity property of the value function $\vartheta(t,\mathbb{P}_{\theta},\mathbb{P}_{\zeta})$ with respect to $t$.

\begin{proposition}
The value function $\vartheta$ is $\frac{1}{2}$-H\"{o}lder continuous in $t$: There exists a constant $C$ such that, for every $t, t^{\prime} \in[0, T]$, $\theta,\zeta\in L^2(\mathcal{F}_t;\mathbb{R}^n),$
\begin{equation}\label{vartheta_con_t}
\left|\vartheta(t, \mathbb{P}_{\theta},\mathbb{P}_{\zeta})-\vartheta\left(t^{\prime}, \mathbb{P}_{\theta},\mathbb{P}_{\zeta}\right)\right| \leq C\big(1+(E[|\zeta|^2])^{\frac{1}{2}}+(E[|\theta|^2])^{\frac{1}{2}}\big)\left|t-t^{\prime}\right|^{\frac{1}{2}} .
\end{equation}
\end{proposition}
\noindent\textit{Proof}. Let $(t, x) \in[0, T] \times \mathbb{R}^n$ and $\delta>0$ be arbitrarily given such that $0<\delta \leq$ $T-t$. We note that due to (\ref{DPP_2_Eq}), for an arbitrarily small $\varepsilon>0$, there exists a control $u=(u^1,u^2)\in\mathcal{U}_{t,t+\delta}$ such that
$$
 |\vartheta \big(t+\delta, \mathbb{P}_{X_{t+\delta}^{t, \theta, \zeta, u}},\mathbb{P}_{X_{t+\delta}^{t, \zeta, u^2}}\big)-\vartheta(t, \mathbb{P}_{\theta},\mathbb{P}_{\zeta})|\leq \varepsilon.
$$
Therefore, with the help of the estimate (\ref{vartheta_con}),
$$
\begin{aligned}
&|\vartheta(t+\delta, \mathbb{P}_{\theta},\mathbb{P}_{\zeta})-\vartheta(t, \mathbb{P}_{\theta},\mathbb{P}_{\zeta})|\\
&\leq\ |\vartheta(t+\delta, \mathbb{P}_{\theta},\mathbb{P}_{\zeta})-\vartheta \big(t+\delta, \mathbb{P}_{X_{t+\delta}^{t, \theta, \zeta, u}},\mathbb{P}_{X_{t+\delta}^{t, \zeta, u^2}}\big)|+|\vartheta \big(t+\delta, \mathbb{P}_{X_{t+\delta}^{t, \theta, \zeta, u}},\mathbb{P}_{X_{t+\delta}^{t, \zeta, u^2}}\big)-\vartheta(t, \mathbb{P}_{\theta},\mathbb{P}_{\zeta})|\\
&\leq |\vartheta(t+\delta, \mathbb{P}_{\theta},\mathbb{P}_{\zeta})-\vartheta \big(t+\delta, \mathbb{P}_{X_{t+\delta}^{t, \theta, \zeta, u}},\mathbb{P}_{X_{t+\delta}^{t, \zeta, u^2}}\big)|+\varepsilon\\
&\leq\  C\big(\mathcal{W}_2(\mathbb{P}_{\theta},\mathbb{P}_{X_{t+\delta}^{t, \theta, \zeta, u}})+\mathcal{W}_2(\mathbb{P}_{\zeta},\mathbb{P}_{X_{t+\delta}^{t, \zeta, u^2}})\big)+\varepsilon\\
&\leq\  C\big(E\big[|X_{t+\delta}^{t, \theta, \zeta, u}-\theta|^2\big]\big)^{\frac{1}{2}}+ C\big(E\big[|X_{t+\delta}^{t, \zeta, u^2}-\zeta|^2\big]\big)^{\frac{1}{2}}+\varepsilon.
\end{aligned}
$$
From Proposition 3.1 of \cite{HL16}, we see that there exists a constant $C>0$,
$$
E\big[|X_{t+\delta}^{t, \theta, \zeta, u}-\theta|^2\big]+E\big[|X_{t+\delta}^{t, \zeta, u^2}-\zeta|^2\big]\leq C(1+E[|\theta|^2]+E[|\zeta|^2])\delta.
$$
Therefore,
$$
\begin{aligned}
|\vartheta(t+\delta, \mathbb{P}_{\theta},\mathbb{P}_{\zeta})-\vartheta(t, \mathbb{P}_{\theta},\mathbb{P}_{\zeta})|\leq\  C\big(1+(E[|\theta|^2])^{\frac{1}{2}}+(E[|\zeta|^2])^{\frac{1}{2}}\big)\delta^{\frac{1}{2}}+\varepsilon,
\end{aligned}
$$
and letting $\varepsilon\downarrow 0$ we get (\ref{vartheta_con_t}). The proof is complete.
\endpf
\vspace{3mm}
\section{Master Bellman equation and viscosity solution}
This section is devoted to the Master Bellman equation related to our value function $\vartheta$. For this we give the definition of viscosity solution, introduced in Burzoni et al. \cite{BIRS20} and adapted to our framework. Some preparation is also made for proving  in Section 6 that the value function $\vartheta$ is the unique viscosity solution of this Master Bellman equation. We recall first some necessary notions to state the Master Bellman equation.

Let us begin with the notion of derivative w.r.t. the measure over the Wasserstein space.
\begin{definition} A function $\varphi: \mathcal{P}_2\left(\mathbb{R}^k\right) \rightarrow \mathbb{R}$ is said to have a linear functional derivative if there exists a function
$$
\frac{\delta \varphi}{\delta \mu}: \mathcal{P}_2\left(\mathbb{R}^k\right) \times \mathbb{R}^k \ni(\mu, x) \mapsto \frac{\delta \varphi}{\delta \mu}(\mu)(x) \in \mathbb{R},
$$
which is continuous with respect to the product topology ($\mathcal{P}_2\left(\mathbb{R}^k\right)$ is equipped with the 2-Wasserstein distance) such that, for any bounded subset $\mathcal{K} \subset \mathcal{P}_2\left(\mathbb{R}^k\right)$, the function $\mathbb{R}^d \ni x \mapsto$ $[\delta \varphi / \delta \mu](\mu)(x)$ is at most of quadratic growth in $x$ uniformly in $\mu$, for $\mu \in \mathcal{K}$, and for all $\mu$ and $\mu^{\prime}$ in $\mathcal{P}_2\left(\mathbb{R}^k\right)$, it holds:
$$
\varphi\left(\mu^{\prime}\right)-\varphi(\mu)=\int_0^1 \int_{\mathbb{R}^d} \frac{\delta \varphi}{\delta \mu}\left(r \mu^{\prime}+(1-r) \mu\right)(x) d(\mu^{\prime}-\mu)(x) d r .
$$
\end{definition}
This leads us to the definition of the L-derivative of $\varphi$ (``L" stands for Lions, the derivative introduced by P.L. Lions \cite{L13}).
\begin{definition} If $\displaystyle\frac{\delta \varphi}{\delta \mu}$ is of class $\mathcal{C}^1$ with respect to the second variable, the L-derivative $\partial_{\mu} \varphi: \mathcal{P}\left(\mathbb{R}^k\right) \times \mathbb{R}^k \rightarrow \mathbb{R}^k$ is defined by
$$
\partial_{\mu} \varphi(\mu, x):=\partial_{x} \big(\frac{\delta \varphi}{\delta \mu} \big)(\mu, x),\quad (\mu,x)\in\mathcal{P}_2(\mathbb{R}^k)\times \mathbb{R}^k.
$$
\end{definition}

For any $\mu\in\mathcal{P}_2(\mathbb{R}^k)$ and any $\mu$-integrable function $f:\mathbb{R}^k\rightarrow\mathbb{R}$, we use the notation
$$
\langle\mu,f\rangle:= \int_{\mathbb{R}^k}f(x)\mu(dx).
$$
For simplicity of notation, we only restrict ourselves to dimension $n=d=k=1$ in what follows. We need to impose some additional assumptions:
\begin{assumption}\label{a2}\ \\
\emph{\textbf{(i)}}  There exists a constant $\kappa_{0}>0$ and a finite set $\mathcal{I} \subset \mathbb{N}$ such that for all $s,s'\in[t,T],\ u\in U,$ $\gamma, \gamma^{\prime} \in \mathcal{P}_2(\mathbb{R}\times U)$ satisfying $\gamma(\mathbb{R}\times\cdot)=\gamma^{\prime}(\mathbb{R}\times\cdot)$, and for every $x\in\mathbb{R}$,
$$
\big|\phi(s, (x, u),\gamma)- \phi(s', (x, u), \gamma')\big|
 \leq \kappa_{0}(|s-s'|+\sum_{i \in \mathcal{I}}|\langle{\mu}-{\mu}', x^{i}\rangle|), \\
$$
where ${\mu}=\gamma(\cdot\times U),{\mu}'=\gamma'(\cdot\times U)$, and $\phi=b_1,b_2,\sigma_1,\sigma_2. $\\
\emph{\textbf{(ii)}} The function $\Phi: \mathbb{R}\times \mathcal{P}_2(\mathbb{R}) \rightarrow \mathbb{R}$ has $\delta$-exponential growth with respect to $x$\ (for some $\delta>0$), uniformly in $\mu$, i.e., there exists a constant $C$, such that for every $x\in\mathbb{R}$,
\begin{equation*}
\sup_{\mu\in\mathcal{P}_2(\mathbb{R})}\big|\Phi(x,\mu)\big|\leq C e^{\delta|x|}.
\end{equation*}
Note that any polynomial has $\delta$-exponential growth.

\end{assumption}

In this section, we study our optimal control problem on $\mathcal{O}:=[0,T)\times \mathcal{M}\times \mathcal{M}$, where $\mathcal{M}\subset \mathcal{P}_2(\mathbb{R})$ is the set of probability measures with $\delta$-exponential moments, i.e.,
\begin{equation}\label{expo}
\mathcal{M}:=\big\{ \mu\in\mathcal{P}_2(\mathbb{R}):\langle\mu, \exp (\delta|\cdot|)\rangle=\int_{\mathbb{R}} \exp (\delta|x|) \mu(\mathrm{d} x)<\infty\big\},
\end{equation}
where $\delta$ is as in Assumption \ref{a2}-(ii). We endow $\mathcal{M}$ with the topology of weak convergence $\sigma(\mathcal{M},C_b(\mathbb{R}))$, where $C_b(\mathbb{R})$ is the space of continuous and bounded functions on $\mathbb{R}$.

\begin{remark}
Condition {\rm(i)} in Assumption \ref{a2} is a form of Lipschitz continuity with respect to cylindrical functions of the measure arguments.
This condition implies the Lipschitz continuity with respect to a distance-like function $d$ ($d(\hat{\mu},\hat{\mu}'):=\sum_{j=1}^{\infty}c_j|\langle \hat{\mu}-\hat{\mu}',f_j\rangle|$, where $\{f_j\}_{j\in\mathbb{N}}$ contains all monomials $f_j(x)=x^j,\ j\geq 0$). As it is shown in Lemma A.1 of Burzoni et al. \cite{BIRS20}, the function $d$, restricted to a chosen compact set, is a metric compatible with weak convergence.
\end{remark}

Now we consider the following second-order partial differential equation on $\mathcal{O}$, the so-called
\noindent{\textbf{Master Bellman equation:}}\\
For any $(t,\mu_1,\mu_2)\in[0,T)\times \mathcal{M}\times \mathcal{M}$,	
\begin{equation}\label{HJB}
\left\{\begin{aligned}&
-\partial_t \vartheta(t, \mu_1,\mu_2)-\mathcal{H}\big(t,\mu_1,\mu_2,\partial_{\mu_1}\vartheta (t,\mu_1,\mu_2;\cdot),\partial_{\mu_2}\vartheta (t,\mu_1,\mu_2;\cdot)\big)=0,
\\& \vartheta(T,\mu_1,\mu_2)= \langle\mu_1,\Phi(\cdot,\mu_2)\rangle,
\end{aligned}\right.
\end{equation}
where, with the notation
$$
\Pi_\mu=\big\{\gamma\in\mathcal{P}_2(\mathbb{R}\times U):\gamma (\cdot\times U)=\mu\big\},\ \mu\in\mathcal{P}_2(\mathbb{R}),
$$
 the Hamiltonian $\mathcal{H}$ is defined by
\begin{equation}\label{Hamilton}
 \begin{aligned}
 \mathcal{H}\big(t,\mu_1,\mu_2,p_1,p_2\big)=
 \inf\big\{
 &\langle \gamma_1,\mathcal{L}^{\mu_1,\mu_2,\gamma_2}_t[p_1]\rangle+
\langle \gamma_2,\overline{\mathcal{L}}^{\mu_1,\mu_2,\gamma_2}_t[p_2]\rangle:\ \gamma_i\in\Pi_{\mu_i},
\ i=1,2 \big\},
\end{aligned}
\end{equation}
$p_1,p_2\in C^1(\mathbb{R})$, with, for $(y,v)\in\mathbb{R}\times U$,
\begin{equation}\label{L}
\begin{gathered}
\mathcal{L}^{\mu_1,\mu_2,\gamma_2}_t[p_1](y,v)=p_1(y)b_2(t,(y,v),\gamma_2)+\frac{1}{2}\partial_y p_1(y)(\sigma_2(t,(y,v),\gamma_2))^2,\\
\overline{\mathcal{L}}^{\mu_1,\mu_2,\gamma_2}_t[p_2](y,v)=p_2(y)b_1(t,(y,v),\gamma_2)+\frac{1}{2}\partial_y p_2(y)(\sigma_1(t,(y,v),\gamma_2))^2.
\end{gathered}
\end{equation}
\begin{remark}\label{Fi_Re-4.2}
When the coefficients $b_1$ and $\sigma_1$ do not depend on the control $u^2$, from PDE (\ref{HJB}) we get the following PDE related with $V(t, x, \mu)$: For $(t,x,\mu)\in[0,T)\times \mathbb{R}\times \mathcal{M},$
\begin{equation}\label{HJB-1}
\left\{\begin{aligned}&
-\partial_t V(t, x,\mu)-\inf_{u\in U}\big\{\partial_{x}V(t,x,\mu)b_2(t, (x, u), \mu)+\frac{1}{2}\partial^2_{x}V(t,x,\mu)(\sigma_2(t, (x, u), \mu))^2\big\}\\&
-\int_{\mathbb{R}}\partial_{\mu}V (t,x,\mu;y)b_1(t,y,\mu)\mu(dy)-\int_{\mathbb{R}}\frac{1}{2}\partial_{y}\partial_{\mu}V (t,x,\mu;y)(\sigma_1(t,y,\mu))^2\mu(dy)=0, \\&
V(T,x,\mu)= \Phi(x,\mu).
\end{aligned}\right.
\end{equation}
Obviously, PDE (\ref{HJB-1}) is just the classical HJB equation if $b_1=\sigma_1=0$, and $b_2, \sigma_2, \Phi$ do not depend on the law. On the other hand, if $b_1=b_2$, and $\sigma_1=\sigma_2$, i.e., all the coefficients are free of controls, then PDE (\ref{HJB}) is just the mean-field PDE obtained  in \cite{17AP}.

\end{remark}

The space $\mathcal{O}$ has a suitable $\sigma$-compact structure, which allows to establish uniqueness in the next section. This form of $\mathcal{O}$ is crucial to obtain uniform integrability of the viscosity test functions as well as some continuity properties of the Hamiltonian. We continue with some preliminaries.

We consider the function
 \begin{equation}\label{edelta}
e_{\delta}(x)
:=\exp \left(\delta(\sqrt{x^{2}+1}-1)\right), \quad x \in \mathbb{R},
\end{equation}
and we note that $e_{\delta}$ is twice continuously differentiable, and
$$
\exp (\delta(|x|-1)) \leq e_{\delta}(x) \leq \exp (\delta|x|) \leq e^{\delta} e_{\delta}(x), \quad  x \in \mathbb{R}.
$$
For $N\in\mathbb{N}$ and $C_0,\delta$ as in Assumption \ref{a2}, let
$$\mathcal{O}_{N}:=\left\{(t, \mu_1,\mu_2) \in[0, T) \times \mathcal{P}_2(\mathbb{R})\times\mathcal{P}_2(\mathbb{R}) \mid\left\langle\mu_i, e_{\delta}\right\rangle \leq N e^{K^{*} t},\ i=1,2\right\},$$
$$\overline{\mathcal{O}}_{N}:=\left\{(t, \mu_1,\mu_2) \in[0, T] \times\mathcal{P}_2(\mathbb{R})\times\mathcal{P}_2(\mathbb{R}) \mid\left\langle\mu_i, e_{\delta}\right\rangle \leq N e^{K^{*} t},\ i=1,2\right\},$$

\noindent where
\begin{equation}\label{K}
K^{*}=K^{*}(C_0,\delta)=C_0\delta+\frac{1}{2}C_0^2(\delta+{\delta}^2).
\end{equation}
The definition of $K^*$ is derived from Lemma \ref{invariant} further down, in order to ensure that $\mathcal{O}$ is invariant for the dynamics (\ref{SDE1})-(\ref{SDE2}).
Let $\overline{\mathcal{O}}:=[0, T] \times \mathcal{M}\times \mathcal{M}$. Note that $\displaystyle \mathcal{O}=[0, T) \times \mathcal{M}\times \mathcal{M}=\cup_{N=1}^{\infty} \mathcal{O}_{N}$ and $\displaystyle\overline{\mathcal{O}}=\cup_{N=1}^{\infty} \overline{\mathcal{O}}_{N}$. For a constant $b$, we put
$$\mathcal{M}_{b}:=\left\{\mu \in \mathcal{P}_2(\mathbb{R}) \mid\left\langle\mu, e_{\delta}\right\rangle \leq b\right\}.$$

Now we prove that, for every $N,\ \mathcal{O}_{N}$, and thus $\mathcal{O}$, remains invariant with respect to the controlled dynamics (\ref{SDE1})-(\ref{SDE2}). This implies that for any given initial law $(t,\mu_1,\mu_2) \in \mathcal{O}_{N}$, we may restrict the Master Bellman equation (\ref{HJB}) to $\mathcal{O}_{N}$.
\begin{lemma}\label{invariant}
Under the Assumption \ref{A1}
, for all $N\in\mathbb{N}$, the set $\mathcal{O}_N$ is invariant for the SDEs (\ref{SDE1})-(\ref{SDE2}), namely,
$$
\big(t,\mathbb{P}_{\theta},\mathbb{P}_{\zeta}\big) \in \mathcal{O}_{N} \Longrightarrow\Big(s,\mathbb{P}_{X_{s}^{t, \theta,\zeta, u}},\mathbb{P}_{X_{s}^{t,\zeta, u^2}}\Big) \in \mathcal{O}_{N}, $$
for all $t\in [0,T],\ s\in[t,T],\ \theta ,\zeta\in L^2(\mathcal{F}_t)$, and $u=(u^1,u^2)\in \mathcal{U}_{t,T},$ where $\big(X_{s}^{t,\zeta, u^2},X_{s}^{t, \theta,\zeta, u}\big)_{s \in[t, T]}$ is the solution to (\ref{SDE1})-(\ref{SDE2}) with initial condition $\big(X_{t}^{t,\zeta, u^2},X_{t}^{t, \theta,\zeta, u}\big)=(\zeta,\theta)$.
\end{lemma}
\noindent\textit{Proof}.
It is sufficient to prove $\big\langle \mathbb{P}_{X_s^{t,\zeta,u^2}},e_{\delta}\big\rangle \leq N e^{K^*s}$ and $\langle \mathbb{P}_{X_s^{t,\theta,\zeta,u}},e_{\delta}\rangle \leq N e^{K^*s}$ for all $s\in[t,T]$. 
 Let $\varphi(x):=\sqrt{x^2+1}-1$, for $x\in\mathbb{R}$. Notice that
\begin{equation}\label{*2}
\begin{aligned}
&0\leq|\varphi' (x)|=\frac{|x|}{\sqrt{1+{x^2}}}<1,\quad 0<\varphi'' (x)=(x^2+1)^{-\frac{3}{2}}\leq 1,\quad x\in\mathbb{R}.
\end{aligned}
\end{equation}
For any $(t,\mathbb{P}_{\theta},\mathbb{P}_{\zeta})\in\mathcal{O}_N$, $(u^1,u^2)\in \mathcal{U}_{t,T},$ we set
$$Y_s^1:=e_{\delta}(X_s^{t,\zeta,u^2})=e^{\delta \varphi(X_s^{t,\zeta,u^2})},\quad s\in[t,T].$$
 By It\^{o}'s formula, for $s\in[t,T]$,
\begin{equation}\label{**1}
\begin{aligned}
Y_s^1=&\ Y_t^1+\int_t^s \delta \varphi'(X_r^{t,\zeta,u^2})Y_r^1b_1\big(r,(X_r^{t,\zeta,u^2},u_r^2),\mathbb{P}_{(X_r^{t,\zeta,u^2},u_r^2)}\big)dr\\
&+\int_t^s \delta \varphi'(X_r^{t,\zeta,u^2})Y_r^1\sigma_1\big(r,(X_r^{t,\zeta,u^2},u_r^2),\mathbb{P}_{(X_r^{t,\zeta,u^2},u_r^2)}\big)dB_r\\
&+\frac{1}{2}\int_t^s Y_r^1\Big(\delta^2 \big(\varphi'(X_r^{t,\zeta,u^2})\big)^2+\delta \varphi^{\prime\prime}(X_r^{t,\zeta,u^2})\Big)\Big(\sigma_1\big(r,(X_r^{t,\zeta,u^2},u_r^2),\mathbb{P}_{(X_r^{t,\zeta,u^2},u_r^2)}\big)\Big)^2dr.
\end{aligned}
\end{equation}
From (\ref{*2}), (\ref{**1}) and Assumption \ref{A1} it follows that
$$
\begin{aligned}
E\left[Y_s^1\right] &\leq E\left[e^{\delta \varphi(\zeta)}\right]+\big(C_0 \delta+\frac{1}{2} C_0^2(\delta^2+\delta)\big) \int_t^s E[Y_r^1] d r\\
&= E\left[e^{\delta \varphi(\zeta)}\right]+K^* \int_t^s E[Y_r^1] d r,\quad s\in [t,T],
\end{aligned}
$$
where $C_0$ is as in Assumption \ref{A1}. Using Gronwall's inequality and the condition $\mathbb{P}_{\zeta}\in\mathcal{M}_{Ne^{K^*t}}$, for $s\in [t,T]$,
$$
\begin{aligned}
&\big\langle\mathbb{P}_{X_s^{t, \zeta, u^2}}, e_\delta\big\rangle=E\left[Y_s^{1}\right] \leq E\left[e^{\delta \varphi(\zeta)}\right]  e^{K^*(s-t)}
&=\left\langle\mathbb{P}_{\zeta}, e_\delta\right\rangle e^{K^*(s-t)}\leq Ne^{K^*t} e^{K^*(s-t)}=Ne^{K^*s}.
\end{aligned}
$$

Let $Y_s^2:= e^{\delta\varphi (X_s^{t,\theta,\zeta,u})},\ s\in [t,T]$. In the same manner as above we can prove that
$$
\big\langle\mathbb{P}_{X_s^{t, \theta, \zeta, u}}, e_\delta\big\rangle \leq N e^{K^*s},\quad s \in [t,T],
$$
which completes the proof.
\endpf
\vspace{3mm}
\begin{lemma}\label{compact}
For $N\in\mathbb{N}$, $\overline{\mathcal{O}}_N$ is a compact subset of $[0,T]\times \mathcal{P}_2(\mathbb{R})\times \mathcal{P}_2(\mathbb{R})$.
\end{lemma}
For the proof of the above lemma, the reader is referred to Lemma 5.2 of Burzoni et al. \cite{BIRS20}.
Now we give the definition of a test function and that of a viscosity solution to the Master Bellman equation (\ref{HJB}), which were first introduced in \cite{BIRS20}. We adapt them here to our framework.

\begin{definition}
A cylindrical function is a mapping of the form $(t, \mu_1,\mu_2) \mapsto F(t,\langle\mu_1, f_1\rangle,$ $\langle\mu_2, f_2\rangle)$ for some functions $f_1,f_2: \mathbb{R} \rightarrow \mathbb{R}$ and $F:[0, T] \times \mathbb{R} \times \mathbb{R}\rightarrow \mathbb{R}$. This function is called cylindrical polynomial, if $f_1,f_2$ are polynomials and $F$ is continuously differentiable.
\end{definition}

We extend the above class to its linear span. For any polynomial $f,$ we denote\ the degree of $f$ by $\operatorname{deg}(f)$.

\begin{definition}[Test Functions]
 For $E \subset \overline{\mathcal{O}}, a$ viscosity test function on $E$ is a function of the form
$$
\varphi(t, \mu_1,\mu_2)=\sum_{j=1}^{\infty} \varphi_{j}(t, \mu_1,\mu_2), \quad(t, \mu_1,\mu_2) \in E,
$$
where $\left\{\varphi_{j}\right\}_{j\in\mathbb{N}}$ is a sequence of cylindrical polynomials that is absolutely convergent at every $(t, \mu_1,\mu_2)$ and, for $i =1,2$,

\begin{equation}\label{testcondition}
\lim _{M \rightarrow \infty} \sum_{j=M}^{\infty} \sup_{\left(t, \mu_1, \mu_2\right) \in \mathcal{O}_N} \sum_{i=1,2}\Big(\langle\mu_i,\left|\partial_{\mu_i} \varphi_j(t, \mu_1, \mu_2;\cdot)\right|\rangle+\langle\mu_i,\left|\partial_y\partial_{\mu_i} \varphi_j(t, \mu_1, \mu_2;\cdot)\right|\rangle\Big)=0 .
\end{equation}
Let $\Gamma_{E}$ be the set of all viscosity test functions on $E$.
\end{definition}
\begin{remark}
According to the definition of the derivative w.r.t. the measure, for a linear function $g(\mu)=\langle\mu,f\rangle$ with $f:\mathbb{R}\rightarrow\mathbb{R}$ differentiable, we have $\partial_{\mu}g(\mu;y)=f'(y).$\
Then, for cylindrical functions with measure arguments of the form $g(\mu)=F(\langle\mu,f\rangle)$, where $F,f: \mathbb{R} \rightarrow \mathbb{R}$ are differentiable, it immediately follows that $\partial_{\mu}g(\mu;y)=F'(\langle\mu,f\rangle) f'(y).$\ For more details the reader is referred, e.g., to \cite{17AP}.
\end{remark}
\begin{definition} For $E \subseteq \overline{\mathcal{O}}$ and $(t, \mu_1,\mu_2) \in E$ with $t<T$, the superjet of a function $u:E\rightarrow\mathbb{R}$ at $(t, \mu_1,\mu_2)$ is given by
$$
\begin{aligned}
&J_{E}^{+} u(t, \mu_1,\mu_2):=\\&\left\{\Big(\partial_{t} \varphi(t, \mu_1,\mu_2), \partial_{\mu_1}\varphi(t, \mu_1,\mu_2; \cdot),\partial_{\mu_2}\varphi(t, \mu_1,\mu_2; \cdot)\Big) \Big| \varphi \in \Gamma_{E},(u-\varphi)(t, \mu_1,\mu_2)=\max _{E}(u-\varphi)\right\}.
\end{aligned}
$$
The subjet of $u$ at $(t, \mu_1,\mu_2)$ is defined as $J_{E}^{-} u(t, \mu_1,\mu_2):=-J_{E}^{+}(-u(t, \mu_1,\mu_2))$.
\end{definition}

\begin{definition}[Viscosity Solution]\label{viscosity} We say that a continuous function $u: \mathcal{O}_{N} \rightarrow \mathbb{R}$ with $u(T,\mu_1,\mu_2)= \langle\mu_1,\Phi(\cdot,\mu_2)\rangle$ is a viscosity subsolution of (\ref{HJB}) on $\mathcal{O}_{N}$ if, for every $(t, \mu_1,\mu_2) \in \mathcal{O}_{N}$,
$$
-\pi_{t}-\mathcal{H}\left(t, \mu_1,\mu_2, \pi_{\mu_1},\pi_{\mu_2}\right)\leq 0, \quad \left(\pi_{t}, \pi_{\mu_1},\pi_{\mu_2}\right) \in J_{\mathcal{O}_{N}}^{+} u(t, \mu_1,\mu_2) .
$$
We say that a continuous function $u: \mathcal{O}_{N} \rightarrow \mathbb{R}$ with $u(T,\mu_1,\mu_2)= \langle\mu_1,\Phi(\cdot,\mu_2)\rangle$ is a viscosity supersolution of (\ref{HJB}) on $\mathcal{O}_{N}$ if for every $(t, \mu_1,\mu_2) \in \mathcal{O}_{N}$,
$$
-\pi_{t}-\mathcal{H}\left(t, \mu_1,\mu_2, \pi_{\mu_1},\pi_{\mu_2}\right) \geq 0, \quad \left(\pi_{t}, \pi_{\mu_1},\pi_{\mu_2}\right) \in J_{\mathcal{O}_{N}}^{-} u(t, \mu_1,\mu_2) .
$$
A viscosity solution of (\ref{HJB}) is a function on $\mathcal{O}$ that is both a subsolution and a supersolution of (\ref{HJB}) on $\mathcal{O}_{N}$, for every $N \in \mathbb{N}$.
\end{definition}
In what follows we show that $(t, \mu_1,\mu_2) \mapsto \mathcal{H}\big(t,\mu_1,\mu_2,\partial_{\mu_1}\varphi,\partial_{\mu_2}\varphi\big)$ is continuous on $\mathcal{O}_{N}$, for any $\varphi\in\Gamma_{\mathcal{O}_N}$. For this purpose, we first state two lemmas.
\begin{lemma}\label{BIR}
Let $\delta>0$ be as in Assumption \ref{a2}. For any constant $b>0$ and any continuous function $g$ with $\delta$-exponential growth,
$$
\sup _{\mu \in \mathcal{M}_{b}}\langle\mu,|g|\rangle<\infty \quad \text { and } \quad \lim _{R \rightarrow \infty} \sup _{\mu \in \mathcal{M}_{b}} \int_{|x| \geq R}|g(x)| \mu(\mathrm{d} x)=0 .
$$
\end{lemma}
For the proof of this lemma, the reader is referred to Lemma 6.8 of \cite{BIRS20}.
\begin{lemma}\label{contin}
Let $\delta>0$ be as in Assumption \ref{a2}. If $\mu_n\rightarrow\mu$ weakly in $\mathcal{M}_b$ as $n\rightarrow \infty$, then for any continuous function $g$ with $\delta$-exponential growth, we have
$$
\langle\mu_n,g\rangle\longrightarrow \langle\mu,g\rangle,\quad n\rightarrow \infty.
$$
\end{lemma}
Lemma \ref{contin} is a straightforward conclusion of Lemma 5.1.12 of Ambrosio, Gigli and Savar\'{e} \cite{AGS05} and the above Lemma \ref{BIR}.

\begin{proposition}\label{4.1}
Under Assumption \ref{A1}, for every $\varphi\in\Gamma_{\mathcal{O}_N}$, the mapping $$(t, \mu_1,\mu_2) \mapsto \mathcal{H}\big(t,\mu_1,\mu_2,\partial_{\mu_1}\varphi(t,\mu_1,\mu_2;\cdot),\partial_{\mu_2}\varphi(t,\mu_1,\mu_2;\cdot)\big)$$ is continuous on $\mathcal{O}_{N}$.
\end{proposition}
\noindent\textit{Proof.}
We first consider the case where $\varphi\in \Gamma_{\mathcal{O}_N}$ is a cylindrical polynomial, namely, $\varphi(t, \mu_1, \mu_2)$ $:=F(t,\langle\mu_1, f_1\rangle,\langle\mu_2, f_2\rangle)$,
where $f_1, f_2$ are polynomials and $F$ is continuously differentiable.
For $\left(t, \mu_1, \mu_2\right) \in \mathcal{O}_N$, $\partial \mu_i \varphi\left(t_1, \mu_1, \mu_2; y\right)=\left(\partial_{z_i}F \right)\left(t_1,\langle\mu_1, f_1\rangle, \langle\mu_2, f_2\rangle\right) f_i^{\prime}(y)$,\ $i=1,2$, and recalling the definition of the Hamiltonian in (\ref{Hamilton}), we have
\begin{equation*}
 \begin{aligned}
 &\mathcal{H}\left(t, \mu_1, \mu_2, \partial_ {\mu_1} \varphi(t,\mu_1,\mu_2;\cdot), \partial_{ \mu_2} \varphi(t,\mu_1,\mu_2;\cdot)\right)\\
& =\inf \Big\{\left\langle\gamma_1, \mathcal{L}_t^{\mu_1, \mu_2, \gamma_2}\left[\partial_{\mu_1} \varphi(t,\mu_1,\mu_2;\cdot)\right]\right\rangle+\big\langle\gamma_2,\overline{\mathcal{L}}_t^{ \mu_1, \mu_2, \gamma_2}\left[\partial_{\mu_2} \varphi(t,\mu_1,\mu_2;\cdot)\right]\big\rangle:
 \gamma_{i} \in \Pi_{\mu_i},\ i=1,2\Big\},
\end{aligned}
\end{equation*}
where
\begin{equation*}
 \begin{aligned}
& \left\langle\gamma_1, \mathcal{L}_t^{\mu_1, \mu_2, \gamma_2}[\partial_{\mu_1} \varphi(t,\mu_1,\mu_2;\cdot)]\right\rangle
\\&=\left\langle\gamma_1(dydv),\left(\partial_{\mu_1} \varphi\right)\left(t, \mu_1, \mu_2 ; y\right) b_2(t,(y, v), \gamma_2)\right\rangle
\\&\qquad +\frac{1}{2}\big\langle\gamma_1(dydv),\left(\partial_y\partial_{\mu_1} \varphi\right)\left(t, \mu_1, \mu_2 ; y\right) \big(\sigma_2(t,(y, v), \gamma_2)\big)^2\big\rangle
\\&=(\partial_{z_1}F)\left(t,\big\langle\mu_1, f_1\rangle, \langle\mu_2, f_2\rangle\right)\\&\qquad\cdot\Big\{\langle\gamma_1(dydv),f_1'(y)b_2(t,(y,v),\gamma_2)\big\rangle+\frac{1}{2}\big\langle\gamma_1(dydv),f_1^{\prime\prime}(y)\big(\sigma_2(t,(y,v),\gamma_2)\big)^2\big\rangle\Big\},
\end{aligned}
\end{equation*}
and
\begin{equation*}
 \begin{aligned}
& \left\langle\gamma_2, \overline{\mathcal{L}}_t^{\mu_1, \mu_2, \gamma_2}[\partial_{\mu_2} \varphi]\right\rangle
\\&=\left\langle\gamma_2(dydv),\left(\partial_{\mu_2} \varphi\right)\left(t, \mu_1, \mu_2 ; y\right) b_1(t,(y, v), \gamma_2)\right\rangle
\\&\qquad +\frac{1}{2}\big\langle\gamma_2(dydv),\left(\partial_y\partial_{\mu_2} \varphi\right)\left(t, \mu_1, \mu_2 ; y\right) \big(\sigma_1(t,(y, v), \gamma_2)\big)^2\big\rangle
\\&=(\partial_{z_2}F)\left(t,\langle\mu_1, f_1\rangle, \langle\mu_2, f_2\rangle\right)\\&\qquad\cdot\Big\{\big\langle\gamma_2(dydv),f_2'(y)b_1(t,(y,v),\gamma_2)\big\rangle+\frac{1}{2}\big\langle\gamma_2(dydv),f_2^{\prime\prime}(y)\big(\sigma_1(t,(y,v),\gamma_2)\big)^2\big\rangle\Big\}.
\end{aligned}
\end{equation*}
It follows from Lemma \ref{contin} that for any continuous $g$ with $\delta$-exponential growth, $\mu\mapsto\langle\mu, g\rangle$ is continuous on $\mathcal{M}_b,\ \text{for all } b>0$.
 Let us first remark that, thanks to our assumptions, namely that $b_1,b_2,\sigma_1,\sigma_2$ are bounded and continuous, all derivatives of $f_1,f_2$ are continuous and of $\delta$-exponential growth, and $F\in C^1$, we have the continuity of the mappings
\begin{equation}\label{a}
\big(t,(\mu_1,\mu_2),(\gamma_1,\gamma_2)\big)\rightarrow \big\langle \gamma_1,\mathcal{L}_t^{\mu_1,\mu_2,\gamma_2}[(\partial_{\mu_1}\varphi)(t,\mu_1,\mu_2;\cdot)]\big\rangle
\end{equation}
and
\begin{equation}\label{b}
\big(t,(\mu_1,\mu_2),\gamma_2\big)\rightarrow \big\langle \gamma_2,\overline{\mathcal{L}}_t^{\mu_1,\mu_2,\gamma_2}[(\partial_{\mu_2}\varphi)(t,\mu_1,\mu_2;\cdot)]\big\rangle
\end{equation}

\noindent over $\mathcal{O}_N\times \big(\mathcal{P}_2(\mathbb{R}\times U)\big)^2$ and $\mathcal{O}_N\times \mathcal{P}_2(\mathbb{R}\times U)$, respectively. Let now $t_n\rightarrow t$ in $[0,T]$ and $\mu_i,\ \mu_i^n\in\mathcal{M}_b,\ n\geq 1$, be such that $\mu_i^n\rightarrow \mu_i,\ n\rightarrow \infty$ weakly, $i=1,2$. Then, as due to the definition of $\mathcal{M}_{b}$, $\displaystyle\sup_{n\geq 1}\int_{\mathbb{R}}e^{\delta|x|}\mu_i^n(dx)\leq  e^\delta b <\infty$, we even have $\mathcal{W}_2(\mu_i^n,\mu_i)\rightarrow 0,\ n\rightarrow \infty$, $i=1,2$. Let now $\gamma_i^n\in\Pi_{\mu_i^n},\ n\geq 1,\ i=1,2,$ be such that
\begin{equation}\label{c}
\begin{aligned}
&\mathcal{H}\big(t_n,\mu_1^n,\mu_2^n,(\partial_{\mu_1}\varphi)(t_n,\mu_1^n,\mu_2^n;\cdot),(\partial_{\mu_2}\varphi)(t_n,\mu_1^n,\mu_2^n;\cdot)\big)+\frac{1}{n}\\
& \geq \big\langle \gamma_1^n,\mathcal{L}_t^{\mu^n_1,\mu^n_2,\gamma^n_2}[(\partial_{\mu_1}\varphi)(t_n,\mu^n_1,\mu^n_2;\cdot)]\big\rangle
 +\big\langle \gamma^n_2,\overline{\mathcal{L}}_t^{\mu^n_1,\mu^n_2,\gamma^n_2}[(\partial_{\mu_2}\varphi)(t_n,\mu^n_1,\mu^n_2;\cdot)]\big\rangle.
\end{aligned}
\end{equation}
Since $\gamma_i^n(\cdot\times U)=\mu_i^n\rightarrow \mu_i,\ n\rightarrow \infty$ in $\mathcal{P}_2(\mathbb{R})$, and, on the other hand, $U$ is compact, we get the tightness of $(\gamma_i^n)_{n\geq 1}$ over $C_b(\mathbb{R}\times U)$. Moreover, as $\displaystyle\sup_{n\geq 1}\int_{ \mathbb{R}\times U}e^{\delta|y|}\gamma_i^n(dydv)=\sup_{n\geq 1}\int_{\mathbb{R}}e^{\delta|y|}\mu_i^n(dy)\leq e^\delta b<\infty$, we can conclude that there exists some $\gamma_i\in\mathcal{P}_2(\mathbb{R}\times U)$ and a subsequence $(\gamma_i^{n_\ell})_{\ell\geq 1}\subset (\gamma_i^{n})_{n\geq 1}$ such that $\mathcal{W}_2(\gamma_i^{n_{\ell}},\gamma_i)\rightarrow 0$, as $\ell\rightarrow \infty$. Consequently, by a standard argument, in $\mathcal{W}_2(\cdot,\cdot)$-sense, as $\gamma_i^{n_{\ell}}\in\Pi_{\mu^{n_\ell}_i},$
$$\gamma_i(\cdot\times U)=\lim_{\ell\rightarrow \infty}\gamma_i^{n_\ell}(\cdot\times U)=\lim_{\ell\rightarrow\infty}\mu_i^{n_\ell}=\mu_i,\ i=1,2,\quad \text{i.e.,}\ \gamma_i\in\Pi_{\mu_i},\ i=1,2.$$
Then, from the continuity of the mapping (\ref{a}) and (\ref{b}), from (\ref{c}) and the definition of the Hamiltonian $\mathcal{H}$,
\begin{equation}\label{d}
\begin{aligned}
& \liminf_{\ell\rightarrow \infty}\mathcal{H}\big(t_{n_\ell},\mu_1^{n_\ell},\mu_2^{n_\ell},(\partial_{\mu_1}\varphi)(t_{n_\ell},\mu_1^{n_\ell},\mu_2^{n_\ell};\cdot),(\partial_{\mu_2}\varphi)(t_{n_\ell},\mu_1^{n_\ell},\mu_2^{n_\ell};\cdot)\big)\\
&\geq \big\langle \gamma_1,\mathcal{L}_t^{\mu_1,\mu_2,\gamma_2}[(\partial_{\mu_1}\varphi)(t,\mu_1,\mu_2;\cdot)]\big\rangle
 +\big\langle \gamma_2,\overline{\mathcal{L}}_t^{\mu_1,\mu_2,\gamma_2}[(\partial_{\mu_2}\varphi)(t,\mu_1,\mu_2;\cdot)]\big\rangle\\
 & \geq \mathcal{H}\big(t,\mu_1,\mu_2,(\partial_{\mu_1}\varphi)(t,\mu_1,\mu_2;\cdot),(\partial_{\mu_2}\varphi)(t,\mu_1,\mu_2;\cdot)\big).
\end{aligned}
\end{equation}
On the other hand, from the definition of the definition of the Hamiltonian $\mathcal{H}$, for any $\varepsilon>0$, there exists $\gamma_i\in\Pi_{\mu_i},\ i=1,2,$ such that
\begin{equation}\label{e}
\begin{aligned}
&\mathcal{H}\big(t,\mu_1,\mu_2,(\partial_{\mu_1}\varphi)(t,\mu_1,\mu_2;\cdot),(\partial_{\mu_2}\varphi)(t,\mu_1,\mu_2;\cdot)\big)+\varepsilon\\
&\geq \big\langle \gamma_1,\mathcal{L}_t^{\mu_1,\mu_2,\gamma_2}[(\partial_{\mu_1}\varphi)(t,\mu_1,\mu_2;\cdot)]\big\rangle
 +\big\langle \gamma_2,\overline{\mathcal{L}}_t^{\mu_1,\mu_2,\gamma_2}[(\partial_{\mu_2}\varphi)(t,\mu_1,\mu_2;\cdot)]\big\rangle.
\end{aligned}
\end{equation}
We choose any $(\xi_i,\theta_i)\in L^2(\mathcal{F};\mathbb{R}\times U)$ such that $\mathbb{P}_{(\xi_i,\theta_i)}=\gamma_i,\ i=1,2,$ and we let $\xi_i^n\in L^2(\mathcal{F};\mathbb{R})$ be such that $\mathbb{P}_{\xi_i^n}=\mu^n_i$, and
\begin{equation*}
\begin{aligned}
E\big[|\xi_i^n-\xi|^2\big]\leq \big(\mathcal{W}_2(\mathbb{P}_{\xi_i^n},\mathbb{P}_{\xi_i})\big)^2+\frac{1}{n}\ \Big(=\big(\mathcal{W}_2(\mu_i^n,\mu_i)\big)^2+\frac{1}{n}\rightarrow 0,\ \text{as }n\rightarrow \infty\Big).
\end{aligned}
\end{equation*}
Then, $\xi_i^n\rightarrow \xi_i$, $n\rightarrow \infty$ in $L^2(\mathcal{F};\mathbb{R})$, $i=1,2$. Let us consider $\gamma_i^n:=\mathbb{P}_{(\xi_i^n,\theta_i)},\ n\geq 1$, and observe that $\gamma_i^n\in\Pi_{\mu_i^n},\ n\geq 1,\ i=1,2.$ Moreover, as
\begin{equation*}
\begin{aligned}
\big(\mathcal{W}_2(\gamma_i^n,\gamma_i)\big)^2\leq E\big[|\xi_i^n-\xi_i|^2\big]\leq \big(\mathcal{W}_2(\mathbb{P}_{\xi_i^n},\mathbb{P}_{\xi_i})\big)^2+\frac{1}{n},
\ n\geq 1,
\end{aligned}
\end{equation*}
we have $\mathcal{W}_2(\gamma_i^n,\gamma_i)\rightarrow 0,\ n\rightarrow \infty,\ i=1,2.$ Consequently, because of the continuity of the mappings (\ref{a}) and (\ref{b}), and considering (\ref{e}) we get
\begin{equation*}
\begin{aligned}
&\mathcal{H}\big(t,\mu_1,\mu_2,(\partial_{\mu_1}\varphi)(t,\mu_1,\mu_2;\cdot),(\partial_{\mu_2}\varphi)(t,\mu_1,\mu_2;\cdot)\big)+\varepsilon\\
&\geq \lim_{n\rightarrow \infty}\Big(\big\langle \gamma^n_1,\mathcal{L}_t^{\mu^n_1,\mu^n_2,\gamma^n_2}[(\partial_{\mu_1}\varphi)(t_n,\mu_1^n,\mu_2^n;\cdot)]\big\rangle
 +\big\langle \gamma^n_2,\overline{\mathcal{L}}_t^{\mu^n_1,\mu^n_2,\gamma^n_2}[(\partial_{\mu_2}\varphi)(t_n,\mu^n_1,\mu^n_2;\cdot)]\big\rangle\Big)\\
 &\geq \limsup_{n\rightarrow\infty}\mathcal{H}\big(t_n,\mu_1^n,\mu_2^n,(\partial_{\mu_1}\varphi)(t_n,\mu_1^n,\mu_2^n;\cdot),(\partial_{\mu_2}\varphi)(t_n,\mu_1^n,\mu_2^n;\cdot)\big).
\end{aligned}
\end{equation*}
Let now $\varepsilon$ tend to zero. Together with (\ref{d}) this yields that, for any converging sequence $(t_n,\mu_1^n,\mu_2^n)\rightarrow (t,\mu_1,\mu_2)$ in $\mathcal{O}_N$, there exists a subsequence $\big(t_{n_\ell},\mu_1^{n_\ell},\mu_2^{n_\ell}\big)_{\ell\geq 1}\subset \big(t_{n},\mu_1^{n},\mu_2^{n}\big)_{n\geq 1}$ such that
\begin{equation*}
\begin{aligned}
&\lim_{\ell\rightarrow \infty}\mathcal{H}\big(t_{n_\ell},\mu_1^{n_\ell},\mu_2^{n_\ell},(\partial_{\mu_1}\varphi)(t_{n_\ell},\mu_1^{n_\ell},\mu_2^{n_\ell};\cdot),(\partial_{\mu_2}\varphi)(t_{n_\ell},\mu_1^{n_\ell},\mu_2^{n_\ell};\cdot)\big)\\
&=\mathcal{H}\big(t,\mu_1,\mu_2,(\partial_{\mu_1}\varphi)(t,\mu_1,\mu_2;\cdot),(\partial_{\mu_2}\varphi)(t,\mu_1,\mu_2;\cdot)\big).
\end{aligned}
\end{equation*}
Consequently, we even have this latter convergence for the whole sequence $\big(t_n,\mu_1^n,\mu_2^n\big)_{n\geq 1}$, i.e., for $n_\ell=\ell,\ \ell\geq 1$, which implies that $\left(t, \mu_1, \mu_2\right) \mapsto \mathcal{H}\left(t, \mu_1, \mu_2, \partial_{\mu_1} \varphi, \partial_{\mu_2} \varphi\right)$ is continuous on $\mathcal{O}_N$. This continuity extends directly to all functions of the type $\displaystyle \varphi^M\left(t, \mu_1, \mu_2\right):=\sum_{j=1}^M \varphi_j\left(t, \mu_1, \mu_2\right)$,
where for any $j,\ \varphi_j\left(t, \mu_1, \mu_2\right):=F_j\left(t,\left\langle\mu_1, f_{1 j}\right\rangle, \left\langle \mu_2,f_{2 j}\right\rangle\right)$ is a cylindrical polynomial.

Now consider a general test function
\begin{equation}\label{generaltest}
\varphi:=\sum_{j=1}^{\infty} \varphi_j \in \Gamma_{\mathcal{O}_N},\
\end{equation}
where, for all $j\in\mathbb{N}$, $\varphi_j$ is a cylindrical polynomial. For $N \in \mathbb{N}$ and the cylindrical polynomials $\left\{\varphi_j\right\}_{j \in \mathbb{N}}\subset \Gamma_{\mathcal{O}_{N}}$, since $b_1,b_2,\sigma_1,\sigma_2$ are bounded, there exists a constant $C>0$ such that
\begin{equation*}
\label{uniform}
\begin{aligned}
&\sup _{\left(t, \mu_1, \mu_2\right) \in \mathcal{O}_N}\left|\mathcal{H}\left(t, \mu_1, \mu_2, \partial_{\mu_1}\varphi_j(t,\mu_1,\mu_2;\cdot), \partial_{\mu_2}\varphi_j(t,\mu_1,\mu_2;\cdot)\right)\right|
\\&\leq C\sup_{\left(t, \mu_1, \mu_2\right) \in \mathcal{O}_N} \sum_{i=1,2}\Big(\langle\mu_i,\left|\partial_{\mu_i} \varphi_j(t, \mu_1, \mu_2;\cdot)\right|\rangle+\langle\mu_i,\left|\partial_y\partial_{\mu_i} \varphi_j(t, \mu_1, \mu_2;\cdot)\right|\rangle\Big).
\end{aligned}
\end{equation*}
This together with (\ref{testcondition}) implies that, for $N\geq M\geq 1$,
$$
\begin{aligned}
& \lim _{M \rightarrow \infty} \sup _{\left(t, \mu_1, \mu_2\right) \in \mathcal{O}_N}\big|\mathcal{H}\left(t, \mu_1, \mu_2, \partial_{\mu_1} \varphi^N(t,\mu_1,\mu_2;\cdot), \partial_{\mu_2} \varphi^N(t,\mu_1,\mu_2;\cdot)\right)\\
&\hspace{3.2cm} -\mathcal{H}\left(t, \mu_1, \mu_2, \partial_{\mu_1} \varphi^M(t,\mu_1,\mu_2;\cdot), \partial_{\mu_2} \varphi^M(t,\mu_1,\mu_2;\cdot)\right)\big| \\
&\leq C \lim _{M \rightarrow \infty} \sum_{j=M}^{\infty}\sup_{\left(t, \mu_1, \mu_2\right) \in \mathcal{O}_N} \sum_{i=1,2}\Big(\langle\mu_i,\left|(\partial_{\mu_i} \varphi_j)(t, \mu_1, \mu_2;\cdot)\right|\rangle+\langle\mu_i,\left|(\partial_y\partial_{\mu_i} \varphi_j)(t, \mu_1, \mu_2;\cdot)\right|\rangle\Big)\\
&= 0,
\end{aligned}
$$
which means that $\mathcal{H}\left(t, \mu_1, \mu_2, \partial_{\mu_1} \varphi^M,\partial_{\mu_2} \varphi^M\right)$ converges uniformly to $\mathcal{H}\left(t, \mu_1, \mu_2, \partial_{\mu_1} \varphi, \partial_{\mu_2} \varphi\right)$ as $M$ tends to infinity. Hence, also $\left(t, \mu_1, \mu_2\right) \mapsto \mathcal{H}\left(t, \mu_1, \mu_2, \partial_{\mu_1} \varphi, \partial_{\mu_2} \varphi\right)$ is continuous on $\mathcal{O}_N$.
\endpf
\vspace{2mm}
\begin{lemma}
Under the Assumptions \ref{A1}, \ref{A2} and \ref{a2}, for all $N$, the function $\vartheta$ is bounded on $\mathcal{O}_N$.
\end{lemma}
\noindent\textit{Proof}. For any $(t, \mathbb{P}_{\theta},\mathbb{P}_{\zeta}) \in \mathcal{O}_{N} ,\ u=(u^1,u^2)\in\mathcal{U}_{t,T},$
$$
\begin{aligned}
\vartheta(t,\mathbb{P}_{\theta},\mathbb{P}_{\zeta})&= \inf_{u^2\in \mathcal{U}^0_{t,T}}E \big[\mathop{\essinf}\limits_{u^1\in \mathcal{U}^0_{t,T}} E \big[ \Phi \big(X_{T}^{t, \theta, \zeta, u},\mathbb{P}_{X_{T}^{t, \zeta, u^2}}\big)\big|\mathcal{F}_t \big]\big]\\
&\leq E\big[ \Phi \big(X_{T}^{t, \theta, \zeta, u},\mathbb{P}_{X_{T}^{t, \zeta, u^2}}\big)\big]=\big\langle \mathbb{P}_{X_{T}^{t, \theta,\zeta, u}},\Phi(\cdot,\mathbb{P}_{X_{T}^{t,\zeta, u^2}})\big\rangle.
\end{aligned}$$

By Assumption \ref{a2}-(ii) and (\ref{PhiLip}) we know that $\displaystyle\sup_{\mu\in\mathcal{M}_b}|\Phi(\cdot,\mu)|$ is continuous and has $\delta$-exponential growth. Then Lemma \ref{BIR} implies that for any constant $b>0$,
$$\sup_{\mu_1\in\mathcal{M}_b}\big\langle\mu_1,\sup_{\mu_2\in\mathcal{M}_b}|\Phi(\cdot,\mu_2)|\big\rangle<\infty.$$
Moreover, thanks to Lemma \ref{invariant}, $\big(T, \mathbb{P}_{X_{T}^{t, \theta,\zeta, u}},\mathbb{P}_{X_{T}^{t,\zeta, u^2}}\big) \in \mathcal{O}_{N}$. Hence, for all $(t,\mathbb{P}_{\theta},\mathbb{P}_{\zeta})\in\mathcal{O}_N$,
$$\vartheta(t,\mathbb{P}_{\theta},\mathbb{P}_{\zeta})\leq\big\langle \mathbb{P}_{X_{T}^{t, \theta,\zeta, u}},\Phi(\cdot,\mathbb{P}_{X_{T}^{t,\zeta, u^2}})\big\rangle\leq\sup_{\mu_1\in\mathcal{M}_{Ne^{K^*T}}}\big\langle\mu_1,\sup_{\mu_2\in\mathcal{M}_{Ne^{K^*T}}}|\Phi(\cdot,\mu_2)|\big\rangle<\infty.$$
\endpf
\section{Main results}
This section is devoted to the both main results of this paper: One states that the value function $\vartheta$ is a viscosity solution of (\ref{HJB}) on $\mathcal{O}$, and the other shows that the comparison theorem for (\ref{HJB}) holds.
\begin{theorem}\label{verify}
 Let the Assumptions \ref{A1} and \ref{A2}
 hold true. Then, for all $N \in \mathbb{N}$, the value function $\vartheta$ is both a viscosity sub- and a supersolution to (\ref{HJB}) on $\mathcal{O}_{N}$.
\end{theorem}
\noindent\textit{Proof}. The relation $\vartheta(T, \mu_1,\mu_2)=\langle \mu_1 , \Phi (\cdot,\mu_2)\rangle,\ \mu_1,\mu_2\in\mathcal{M}_{N e^{K^{*} T}}$, can be obtained directly from the definition of $\vartheta$. In the following we show that $\vartheta$ is a viscosity solution to (\ref{HJB}) on $\mathcal{O}_{N}$.\\
\textbf{Step 1}: We first prove $\vartheta$ is a viscosity subsolution for $t<T$. Suppose that for $\varphi \in \Gamma_{\mathcal{O}_N}$ and $\left(t, \mu_{1}, \mu_{2}\right) \in \mathcal{O}_{N}$,
$$
0=\left(\vartheta-\varphi\right)\left(t, \mu_{1}, \mu_{2}\right)=\max _{\mathcal{O}_{N}}\left(\vartheta-\varphi\right).
$$
This implies that $\vartheta\leq\varphi$ on $\mathcal{O}_{N}$. We choose $\theta, \zeta\in L^2(\mathcal{F}_t)$ such that $(\mathbb{P}_{\theta},\mathbb{P}_{\zeta})=(\mu_1,\mu_2)$.

We fix any control process $u=\left(u^{1}, u^{2}\right) \in
\mathcal{U}_{t, \beta}$, where $\beta:=t+h$, and $0<h<T-t$, and we denote by $\big( X^{t, \zeta, u^{2}},X^{t, \theta, \zeta, u}\big)$ the solution of the SDEs (\ref{SDE1})-(\ref{SDE2}) with initial value $(t,\zeta,\theta)$ and control $u$. Then we apply the dynamic programming principle in Theorem \ref{DPP2} from $t$ to $\beta$ and get
$$
\vartheta\left(t, \mu_{1}, \mu_{2}\right) \leq \vartheta\big(\beta,\mathbb{P}_{X_{\beta}^{t, \theta, \zeta,u}},\mathbb{P}_{X_{\beta}^{t,\zeta, u^{2}}}\big)
\leq \varphi\big(\beta,\mathbb{P}_{X_{\beta}^{t, \theta, \zeta,u}},\mathbb{P}_{X_{\beta}^{t,\zeta, u^{2}}}\big).
$$
By applying It\^{o}'s formula to $\varphi\big(\cdot, \mathbb{P}_{X_{\cdot}^{t, \theta, \zeta,u}},\mathbb{P}_{X_{\cdot}^{t,\zeta, u^{2}}}\big)$ (see, e.g., Theorem 7.1 in \cite{17AP}) between $t$ and $\beta$, and using the above inequality
we obtain
$$
\begin{aligned}
&0\geq\big(\vartheta\left(t, \mu_{1}, \mu_{2}\right)-\varphi(t, \mu_{1}, \mu_{2})\big)-\int_{t}^{\beta}\Big\{\partial_{t} \varphi\big(r, \mathbb{P}_{X_{r}^{t, \theta, \zeta,u}},\mathbb{P}_{X_{r}^{t,\zeta, u^{2}}}\big)
\\ & +{E}\Big[\partial_{\mu_1}\varphi\big(r,\mathbb{P}_{X_{r}^{t, \theta, \zeta,u}},\mathbb{P}_{X_{r}^{t,\zeta, u^{2}}};{X}_{r}^{t, \theta, \zeta,u}\big) b_2\Big(r,({X}_{r}^{t, \theta,
 \zeta,u},{u}_r^1),\mathbb{P}_{(X_{r}^{t,\zeta, u^{2}},u_r^2)}\Big)\Big]\\&
+ {E}\Big[\partial_{\mu_2}\varphi\big(r,\mathbb{P}_{X_{r}^{t, \theta, \zeta,u}},\mathbb{P}_{X_{r}^{t,\zeta, u^{2}}};{X}_{r}^{t,  \zeta,u^2}\big) b_1\Big(r,({X}_{r}^{t,\zeta,u^2},{u}_r^2),\mathbb{P}_{(X_{r}^{t,\zeta, u^{2}},u_r^2)}\Big)\Big]\\&
 + \frac{1}{2}{E}\Big[\partial_y\partial_{\mu_1}\varphi\big(r,\mathbb{P}_{X_{r}^{t, \theta, \zeta,u}},\mathbb{P}_{X_{r}^{t,\zeta, u^{2}}};{X}_{r}^{t, \theta, \zeta,u}\big)
 \Big(\sigma_2\Big(r,({X}_{r}^{t, \theta, \zeta,u},{u}_r^1),\mathbb{P}_{(X_{r}^{t,\zeta, u^{2}},u_r^2)}\Big)\Big)^2\Big]\\&
  + \frac{1}{2}{E}\Big[\partial_y\partial_{\mu_2}\varphi\big(r,\mathbb{P}_{X_{r}^{t, \theta, \zeta,u}},\mathbb{P}_{X_{r}^{t,\zeta, u^{2}}};{X}_{r}^{t,  \zeta,u^2}\big)
 \Big(\sigma_1\Big(r,({X}_{r}^{t,  \zeta,u^2},u_r^2),\mathbb{P}_{(X_{r}^{t,\zeta, u^{2}},u_r^2)}\Big)\Big)^2\Big]\Big\}dr.
\end{aligned}
$$

\noindent This holds for all $h>0$. By dividing by $h >0$ and letting $h\downarrow 0$, we use the continuity of the intergrand in $r=t$. This yields
$$
\begin{aligned}
& -\partial_{t} \varphi\big(t, \mathbb{P}_{\theta},\mathbb{P}_{\zeta}\big)
-{E}\Big[\partial_{\mu_1}\varphi\big(t, \mathbb{P}_{\theta},\mathbb{P}_{\zeta};{\theta}\big) b_2\big(t,({\theta},{u}_t^1),\mathbb{P}_{(\zeta,{u}_t^2)}\big)\Big]\\&
-{E}\Big[\partial_{\mu_2}\varphi\big(t, \mathbb{P}_{\theta},\mathbb{P}_{\zeta};{\zeta}\big) b_1\big(t,({\zeta},{u}_t^2),\mathbb{P}_{(\zeta,{u}_t^2)}\big)\Big]\\&
-\frac{1}{2}{E}\Big[\partial_y\partial_{\mu_1}\varphi\big(t, \mathbb{P}_{\theta},\mathbb{P}_{\zeta};{\theta} \big) \big(\sigma_2\big(t,({\theta},{u}_t^1),\mathbb{P}_{(\zeta,{u}_t^2)}\big)\big)^2\Big]\\&
-\frac{1}{2}{E}\Big[\partial_y\partial_{\mu_2}\varphi\big(t, \mathbb{P}_{\theta},\mathbb{P}_{\zeta};{\zeta} \big) \big(\sigma_1\big(t,({\zeta},{u}_t^2),\mathbb{P}_{(\zeta,{u}_t^2)}\big)\big)^2\Big]\leq 0,
\end{aligned}
$$
for all $u\in\mathcal{U}_{t,t+h}$, and  $\theta, \zeta\in L^2(\mathcal{F}_t)$ with $(\mathbb{P}_{\theta},\mathbb{P}_{\zeta})=(\mu_1,\mu_2)$.

Therefore, we have
\begin{equation}
\begin{aligned}
& -\partial_{t} \varphi\big(t, \mu_1,\mu_2\big)-\inf_{u_t^1,u_t^2\in L^2(\mathcal{F}_t;U)}\Big\{{E}\Big[\partial_{\mu_1}\varphi\big(t, \mu_1,\mu_2;{\theta}\big) b_2\big(t,({\theta},{u}_t^1),\mathbb{P}_{(\zeta,{u}_t^2)}\big)\Big]\\&
+{E}\Big[\partial_{\mu_2}\varphi\big(t,  \mu_1,\mu_2;{\zeta}\big) b_1\big(t,({\zeta},{u}_t^2),\mathbb{P}_{(\zeta,{u}_t^2)}\big)\Big]\\&
+\frac{1}{2}{E}\Big[\partial_y\partial_{\mu_1}\varphi\big(t,  \mu_1,\mu_2;{\theta} \big) \big(\sigma_2\big(t,({\theta},{u}_t^1),\mathbb{P}_{(\zeta,{u}_t^2)}\big)\big)^2\Big]\\&
+\frac{1}{2}{E}\Big[\partial_y\partial_{\mu_2}\varphi\big(t,  \mu_1,\mu_2;{\zeta} \big) \big(\sigma_1\big(t,({\zeta},{u}_t^2),\mathbb{P}_{(\zeta,{u}_t^2)}\big)\big)^2\Big]\Big\}\leq 0,
\end{aligned}
\end{equation}
which is equivalent to the following inequality:
\begin{equation}
\begin{aligned}
& -\partial_{t} \varphi\big(t, \mu_1,\mu_2\big)-
 \inf_{\gamma_i\in\Pi_{\mu_i},i=1,2}\Big\{\left\langle\gamma_1(dydv),\left(\partial_{\mu_1} \varphi\right)\left(t, \mu_1, \mu_2 ; y\right) b_2(t,(y, v), \gamma_2)\right\rangle\\&
+\left\langle\gamma_2(dydv),\left(\partial_{\mu_2} \varphi\right)\left(t, \mu_1, \mu_2 ; y\right) b_1(t,(y, v), \gamma_2)\right\rangle\\&+\frac{1}{2}\big\langle\gamma_1(dydv),\left(\partial_y\partial_{\mu_1} \varphi\right)\left(t, \mu_1, \mu_2 ; y\right) \big(\sigma_2(t,(y, v), \gamma_2)\big)^2\big\rangle\\
&+\frac{1}{2}\big\langle\gamma_2(dydv),\left(\partial_y\partial_{\mu_2} \varphi\right)\left(t, \mu_1, \mu_2 ; y\right) \big(\sigma_1(t,(y, v), \gamma_2)\big)^2\big\rangle
\Big\}\leq 0.
\end{aligned}
\end{equation}
Hence,
\begin{equation}
\begin{aligned}
-\partial_t \varphi(t, \mu_1,\mu_2)- \mathcal{H}\big(t,\mu_1,\mu_2,\partial_{\mu_1}\varphi(t,\mu_1,\mu_2;\cdot),\partial_{\mu_2}\varphi(t,\mu_1,\mu_2;\cdot)\big)\leq 0,
\end{aligned}
\end{equation}
where $\mathcal{H}$ is just the Hamiltonian $\mathcal{H}$ introduced in (\ref{Hamilton}).\\
\textbf{Step 2}: We prove $\vartheta$ is a viscosity supersolution for $t<T$. Suppose that, for $\varphi \in \Gamma_{\mathcal{O}_N}$ and $\left(t, \mu_{1}, \mu_{2}\right) \in \mathcal{O}_{N}$,
$$
0=\left(\vartheta-\varphi\right)\left(t, \mu_{1}, \mu_{2}\right)=\min _{\mathcal{O}_{N}}\left(\vartheta-\varphi\right).
$$
This implies that $\vartheta\geq\varphi$ on $\mathcal{O}_{N}$. In view of Lemma 7.1 in \cite{BIRS20}, we may assume without loss of generality that above minimum is strict. We prove the result by contradiction. We assume on the contrary that
\begin{equation}\label{contrary}
-\partial_t \varphi(t, \mu_1,\mu_2)- \mathcal{H}\big(t,\mu_1,\mu_2,\partial_{\mu_1}\varphi(t,\mu_1,\mu_2;\cdot),\partial_{\mu_2}\varphi(t,\mu_1,\mu_2;\cdot)\big)< 0.
\end{equation}
Then there is some $\zeta>0$ such that
\begin{equation*}\label{contrary}
-\partial_t \varphi(t, \mu_1,\mu_2)- \mathcal{H}\big(t,\mu_1,\mu_2,\partial_{\mu_1}\varphi(t,\mu_1,\mu_2;\cdot),\partial_{\mu_2}\varphi(t,\mu_1,\mu_2;\cdot)\big)<-\zeta< 0.
\end{equation*}
By the continuity of the function $\mathcal{H}$ on $\mathcal{O}_N$ stated in Proposition \ref{4.1}, there exists a neighbourhood $B$ of $(t,\mu_1,\mu_2)$ such that
\begin{equation}\label{6.4}
-\partial_t \varphi(t', \mu_1',\mu_2')
-\mathcal{H}\big(t',\mu_1',\mu_2',\partial_{\mu_1}\varphi(t',\mu_1',\mu_2';\cdot),\partial_{\mu_2}\varphi(t',\mu_1',\mu_2';\cdot)\big)
\leq -\zeta,
\end{equation}
for all $(t',\mu_1',\mu_2')\in B_N:=B\cap \mathcal{O}_N$. We choose $\theta, \zeta\in L^2(\mathcal{F}_t)$ again such that $(\mathbb{P}_{\theta},\mathbb{P}_{\zeta})=(\mu_1,\mu_2)$, and we fix an arbitrary control $u=(u^1,u^2)\in \mathcal{U}_{t,\beta}$, and put
$$\mu_r^{t,\theta,\zeta,u}:=\mathbb{P}_{X_{r}^{t, \theta, \zeta,u}},\quad \mu_r^{t,\zeta,u^2}:=\mathbb{P}_{X_{r}^{t,\zeta, u^{2}}}.$$
 Let us consider the deterministic time
$$
\beta:=\inf \left\{s \geq t:\Big(s,\mu_{s}^{t,\theta,\zeta,u}, \mu_{s}^{t,\zeta,u^2}\Big) \notin B_{N}\right\} \wedge T \text {. }
$$
Note that $\beta>t$. This implies that $\Big(s,\mu_{s}^{t,\theta,\zeta,u}, \mu_{s}^{t,\zeta,u^2}\Big) \in B_{N}$, for all $t\leq s<\beta$. 
 Applying It\^{o}'s formula to $\varphi\big(\cdot, \mu_{\cdot}^{t,\theta,\zeta,u}, \mu_{\cdot}^{t,\zeta,u^2}\big)$ between $t$ and $\beta$, 
we obtain from (\ref{6.4}) that
$$
\begin{aligned}
&\varphi\big(\beta, \mu_{\beta}^{t,\theta,\zeta,u}, \mu_{\beta}^{t,\zeta,u^2}\big)=\varphi\left(t, \mu_{1}, \mu_{2}\right)+\int_{t}^{\beta}\Big\{\partial_{t} \varphi\big(r, \mu_r^{t,\theta,\zeta,u}, \mu_r^{t,\zeta,u^2}\big)
\\ &\hspace{4cm} +\big\langle \gamma_r^{t,\theta,\zeta,u},\mathcal{L}^{\mu_r^{t,\theta,\zeta,u}, \mu_r^{t,\zeta,u^2},\gamma_r^{t,\zeta,u^2}}_r[\partial_{\mu_1}\varphi(r,\mu^{t,\theta,\zeta,u}_r,\mu^{t,\zeta,u^2}_r;\cdot)]\big\rangle\\
&\hspace{4cm}+
\big\langle \gamma_r^{t,\zeta,u^2},\overline{\mathcal{L}}^{\mu_r^{t,\theta,\zeta,u}, \mu_r^{t,\zeta,u^2},\gamma_r^{t,\zeta,u^2}}_r[\partial_{\mu_2}\varphi(r,\mu^{t,\theta,\zeta,u}_r,\mu^{t,\zeta,u^2}_r;\cdot)]\big\rangle\Big\}dr\\
&\hspace{3.2cm}\geq  \varphi(t,\mu_1,\mu_2)+(\beta-t)\zeta,
\end{aligned}
$$
where $$\gamma_r^{t,\theta,\zeta,u}=\mathbb{P}_{\big(X_{r}^{t, \theta, \zeta,u},u_r^1\big)},\quad \gamma_r^{t,\zeta,u^2}=\mathbb{P}_{\big(X_{r}^{t, \zeta,u^2},u_r^2\big)}.
$$
Consequently,
$$
\zeta(\beta-t)+\varphi\left(t, \mu_{1}, \mu_{2}\right)\leq\varphi\big(\beta, \mu_{\beta}^{t,\theta,\zeta,u}, \mu_{\beta}^{t,\zeta,u^2}\big).
$$
Note that $\overline{\mathcal{O}}_{N} \backslash B_{N}=\overline{\mathcal{O}}_{N} \backslash B$ is compact. On the other hand, we have supposed  that $\vartheta-\varphi$ has a strict minimum $0$ at $(t, \mu_1,\mu_2)$ on $\overline{\mathcal{O}}_{N}$, then
$$
\zeta(\beta-t)+\varphi\left(t, \mu_{1}, \mu_{2}\right)\leq\varphi\big(\beta, \mu_{\beta}^{t,\theta,\zeta,u}, \mu_{\beta}^{t,\zeta,u^2}\big) \leq \vartheta\big(\beta, \mu_{\beta}^{t,\theta,\zeta,u}, \mu_{\beta}^{t,\zeta,u^2}\big).
$$
Since the $(\varphi-\vartheta)\left(t, \mu_1,\mu_2\right) 
=0$, we have
$$
\zeta(\beta-t)+\vartheta\left(t, \mu_{1}, \mu_{2}\right)\leq\vartheta\big(\beta, \mu_{\beta}^{t,\theta,\zeta,u}, \mu_{\beta}^{t,\zeta,u^2}\big).
$$
As $\zeta(\beta-t)>0$ and the above inequality holds for all $u\in \mathcal{U}_{t,\beta}$, it is in contradiction with Theorem \ref{DPP2}. Hence, $\vartheta$ is a viscosity supersolution to (\ref{HJB}).
\endpf
\vspace{2mm}
The remaining part of this section is devoted to a comparison theorem for the value function $\vartheta$. For this, we first need to introduce some definitions and notations. \begin{definition}
 We say that a set of polynomials $\mathcal{X}$ has the $(*)$-property, if it satisfies
$$\text{for all }g \in \mathcal{X}, \ g^{(i)} \in \mathcal{X},\ 0\leq i\leq \operatorname{deg}(g),$$
where $g^{(i)}$ is the $i$-th derivative of $g$, and $\operatorname{deg}(g)$ denotes the degree of polynomial $g$. Let $\sum$ be the collection of all sets of polynomials that have the $(*)$-property.
\end{definition}
For $f$ a real polynomial we set
$$
\mathcal{X}(f):=\bigcap_{\mathcal{X} \in \Sigma, f \in \mathcal{X}} \mathcal{X} \text {. }
$$
We can easily check the following properties of $\mathcal{X}(f)$.
\begin{lemma} \label{property}
For every polynomial $f$, we have:\\
$\mathrm{(a)}$ $\mathcal{X}(f)$ is the smallest set of polynomials with the $(*)$-property that includes $f$;\\
$\mathrm{(b)}$ For every $g \in \mathcal{X}(f),\ \mathcal{X}(g) \subset \mathcal{X}(f)$;\\
$\mathrm{(c)}$ $\mathcal{X}(f)$ is finite.
\end{lemma}

Let $\displaystyle\Theta:=\bigcup_{j=1}^{\infty} \mathcal{X}\left(\psi_j\right)$, where $\psi_j(x)=x^j,\ x\in\mathbb{R}$. Then $\Theta$ is countable, $\{\psi_j\}_{j=1}^{\infty} \subset \Theta$ and for any $f \in \Theta,\ \mathcal{X}(f) \subset \Theta$. Let $\left\{f_j\right\}_{j=1}^{\infty}$ be an enumeration of $\Theta$. We define the finite index set $I_j$ by putting
$$\mathcal{X}\left(f_j\right)=\big\{f_i \mid i \in I_j\big\},\ j\geq 1.$$
Then, for all $i \in I_j$, we have $\mathcal{X}\left(f_i\right) \subset \mathcal{X}\left(f_j\right)$ and, therefore, $I_i \subset I_j$. Moreover, we define for $b>0$ and $j\in\mathbb{N}$
$$
c_j(b):=\big(\sum_{k\in I_j}2^k\big)^{-1}\big(\sum_{k\in I_j}s_k(b)\big)^{-2},
$$
where $\displaystyle s_j(b):=1+\sup _{\mu \in \mathcal{M}_b}\left\langle\mu, f_j\right\rangle$. We observe that it follows from Lemma \ref{contin} that $1\leq s_j(b)<\infty$, for all $j\in\mathbb{N}$. Since $f_j \in \mathcal{X}\left(f_j\right)$, we have $j \in I_j$, and, thus, $c_j(b) \leq 2^{-j}$. Hence, $\sum_{j=1}^{\infty} c_j(b) \leq 1$.
 Moreover, for $i \in I_j$, from $I_i \subset I_j$ we get $c_j(b) \leq c_i(b)$. Finally, we observe that, by the definition of $s_j(b)$ and $c_j(b)$,
\begin{equation}\label{c_j}
\sum_{j=1}^{\infty} c_j(b)\left\langle\mu, f_j\right\rangle^2 \leq 1, \quad  \mu \in \mathcal{M}_b.
\end{equation}
\begin{theorem}[Comparison Theorem]\label{comparison}
 In addition to our standing assumptions, we suppose that
 \begin{equation}
 b_j(t,(y,v),\gamma)=b_j(t,\gamma),\quad \sigma_j(t,(y,v),\gamma)=\sigma(t,\gamma),
 \end{equation}
$(t,(y,v),\gamma)\in [0,T]\times(\mathbb{R}\times U)\times \mathcal{P}_2(\mathbb{R}\times U),\ j=1,2,$ i.e., the coefficients $b_j,\sigma_j$ are independent of $(y,v)$.
 Let $u\in C(\mathcal{O}_N)$ be a viscosity subsolution to HJB equation (\ref{HJB}) on $\mathcal{O}_{N}$ and $v\in C(\mathcal{O}_N)$ be a viscosity supersolution to HJB equation (\ref{HJB}) on $\mathcal{O}_{N}$, satisfying $u(T, \mu_1,\mu_2) \leq v(T, \mu_1,\mu_2)$, for any $(T, \mu_1,\mu_2) \in$ $\overline{\mathcal{O}}_{N}$. Then $u \leq v$ on $\overline{\mathcal{O}}_{N}$.
\end{theorem}
\begin{remark}
 Burzoni et al. \cite{BIRS20} consider coefficients $(b,\sigma)(t,\mu,v),$\ $(t,\mu,v)\in [0,T]\times \mathcal{P}(\mathbb{R})\times U$.  While they use only deterministic control processes, we overcome this difficulty  by considering our stochastic control in the law.
\end{remark}
\noindent\textit{Proof of Theorem \ref{comparison}.}
Fix $N \in \mathbb{N}$ and let $c_j:=c_j\left(N e^{K^* T}\right)$ (recall that $K^*$ is defined in (\ref{K})). Then, for all $\left(t, \mu_1, \mu_2\right) \in \overline{\mathcal{O}}_N$, $\mu_1, \mu_2 \in \mathcal{M}_{N e^{K^* t}} \subset \mathcal{M}_{N e^{K^* T}}$, it follows from (\ref{c_j}) that
$$
\sup_{\left(t, \mu_1, \mu_2\right) \in \overline{\mathcal{O}}_N} \sum_{j=1}^{\infty} c_j\langle\mu_i, f_j\rangle^2 \leq 1, \quad i=1,2 \text {. }
$$
We suppose that
$$
\sup _{\overline{\mathcal{O}}_N}(u-v)>0 \text {, }
$$
and we prove that this leads to a contradiction. Since $u-v$ is continuous and $\overline{\mathcal{O}}_N$ is compact, the maximum
$$
\ell:=\max _{\left(t, \mu_1, \mu_2\right) \in \overline{\mathcal{O}}_N}\big((u-v)(t, \mu_1, \mu_2)-2\eta (T-t)\big)
$$
can be achieved and there exists sufficiently small $\eta_0$, such that, for all $\eta \in\left(0, \eta_0\right]$, we have $\ell>0$.

Now we use the standard argument of doubling variables to construct test functions for $u$ and $v$.\smallskip\\
\textbf{Step\ 1.\ Doubling of variables.}

For $\varepsilon>0$ and $\eta \in\left(0, \eta_0\right]$, we define
\begin{equation}\label{phi}
\begin{aligned}
\phi_{\varepsilon}\left(t, \mu_1, \mu_2, s, \nu_1, \nu_2\right):=&u\left(t, \mu_1, \mu_2\right)-v\left(s, \nu_1, \nu_2\right)
\\&-\frac{1}{\varepsilon}\big((t-s)^2+d(\mu_1, \nu_1)+d(\mu_2, \nu_2)\big)-\eta(T-t+T-s),
\end{aligned}
\end{equation}
where $d$ is a distance-like function defined for $\mu, \nu \in \mathcal{M}_{N e^{K^* T}}$ by the relation
\begin{equation}\label{dist}
d(\mu, \nu):=\sum_{j=1}^{\infty} c_j\langle\mu-\nu, f_j\rangle^2;
\end{equation}
recall that $\left\{f_j\right\}_{j=1}^{\infty}=\Theta=\bigcup_{j=1}^{\infty} \mathcal{X}\left(x^j\right)$ (see Lemma \ref{property}). 
We observe that $d(\cdot,\cdot)$ is compatible with the weak convergence in $\mathcal{M}_{Ne^{K^* T}}$: Both generate the same topology.
  Since $\phi_{\varepsilon}$ is continuous and $\overline{\mathcal{O}}_N$ is compact, $\phi_{\varepsilon}$ admits a maximizer $(t^*_\varepsilon, \mu^*_{1,\varepsilon}, \mu^*_{2,\varepsilon}, s^*_\varepsilon, \nu^*_{1,\varepsilon}, \nu^*_{2,\varepsilon})$ such that
\begin{equation}\label{max}
\phi_\varepsilon^*:= \max_{(t,\mu_1,\mu_2)\in\overline{\mathcal{O}}_N,(s,\nu_1,\nu_2)\in\overline{\mathcal{O}}_N}\phi_\varepsilon\left(t, \mu_1, \mu_2, s, \nu_1, \nu_2\right)=\phi_\varepsilon (t^*_\varepsilon, \mu^*_{1,\varepsilon}, \mu^*_{2,\varepsilon}, s^*_\varepsilon, \nu^*_{1,\varepsilon}, \nu^*_{2,\varepsilon}).
\end{equation}
Notice that if we take $\left(s, \nu_1, \nu_2\right)=\left(t, \mu_1, \mu_2\right)$ in (\ref{phi}),  $$\phi_{\varepsilon}\left(t, \mu_1, \mu_2, t, \mu_1, \mu_2\right)=(u-v)\left(t, \mu_1, \mu_2\right)-2 \eta(T-t),$$ then we obtain
\begin{equation}\label{>0}
\phi_\varepsilon^*\geq\max_{(t,\mu_1,\mu_2)\in\overline{\mathcal{O}}_N}\phi_\varepsilon\left(t, \mu_1, \mu_2, t, \mu_1, \mu_2\right)=\ell>0.
\end{equation}
Set
\begin{equation}\label{varphi_u}
\varphi_u\left(t,\mu_1, \mu_2\right):=\frac{1}{\varepsilon} \big((t-s_\varepsilon^*)^2 + d\left(\mu_1, \nu^*_{1,\varepsilon}\right)+d\left(\mu_2, \nu^*_{2,\varepsilon}\right)\big)+\eta\left(T-t+T-s_{\varepsilon}^*\right).
\end{equation}
We can check that $\varphi_u$ is a test function for the viscosity subsolution $u$ at $\left(t_{\varepsilon}^*, \mu_{1, \varepsilon}^*, \mu_{2, \varepsilon}^*\right)$. Actually, from the form of $\varphi_u$ ((\ref{dist}) and (\ref{varphi_u})) and (\ref{max}), we observe that $\varphi_u \in \Phi_{\overline{\mathcal{O}}_N}$ and $$u\left(t, \mu_1, \mu_2\right)-\varphi_u\left(t, \mu_1, \mu_2\right)=\phi_{\varepsilon}\left(t, \mu_1, \mu_2, s_{\varepsilon}^*, \nu_{1, \varepsilon}^*, \nu_{2 ,\varepsilon}^*\right)+v\left(s_{\varepsilon}^*, \nu_{1, \varepsilon}^*, \nu_{2, \varepsilon}^*\right)$$ achieves its maximum $\phi_{\varepsilon}^*+v\left(s_{\varepsilon}^*, \nu_{1, \varepsilon}^*, \nu_{2, \varepsilon}^*\right)$ at $\left(t^*_\varepsilon, \mu_{1,\varepsilon}^*,\mu_{2, \varepsilon}^*\right)$.
Similarly, by setting
\begin{equation}\label{varphi_v}
\varphi_v\left(s, \nu_1, \nu_2\right):=-\frac{1}{\varepsilon}\big(\left(t_{\varepsilon}^*-s\right)^2+d\left(\mu_{1, \varepsilon}^*, \nu_1\right)+ d\left(\mu_{2,\varepsilon}^*,\nu_2\right)\big)-\eta\left(T-t_{\varepsilon}^*+T-s\right),
\end{equation}
we get a test function $\varphi_v$ for the viscosity supersolution $v$ at $\left(s_{\varepsilon}^*, \nu_{1, \varepsilon}^*, \nu_{2, \varepsilon}^*\right)$. Indeed, $\varphi_v \in \Phi_{\overline{\mathcal{O}}_N}$ and
$$
\begin{aligned}
& v\left(s, \nu_1, \nu_2\right)-\varphi_v\left(s, \nu_1, \nu_2\right) =u\left(t^*_\varepsilon, \mu_{1,\varepsilon}^*,\mu_{2, \varepsilon}^*\right)-\phi_{\varepsilon}\left(t^*_\varepsilon, \mu_{1,\varepsilon}^*,\mu_{2, \varepsilon}^*,s,\nu_1,\nu_2\right)
\end{aligned}
$$
achieves its minimum $u\left(t^*_\varepsilon, \mu_{1,\varepsilon}^*,\mu_{2, \varepsilon}^*\right)-\phi_{\varepsilon}^*$ at $\left(s_{\varepsilon}^*, \nu_{1, \varepsilon}^*, \nu_{2, \varepsilon}^*\right)$.\smallskip\\
\textbf{Step 2.\ Some preparations.}

Recalling that $u$ and $v$ are continuous, and $\overline{\mathcal{O}}_N$ is compact, we put $\displaystyle M:=\max_{\overline{\mathcal{O}}_N}u \in \mathbb{R}$, and  $\displaystyle m:=\min_{\overline{\mathcal{O}}_N} v \in \mathbb{R}$. Let us define $\zeta_{\varepsilon}:=d\left(\mu_{1, \varepsilon}^*, \nu_{1, \varepsilon}^*\right)+d\left(\mu_{2, \varepsilon}^*, \nu_{2, \varepsilon}^*\right)+\left(t^*_{\varepsilon}-s_\varepsilon^*\right)^2$. From (\ref{max}) and (\ref{>0}), it follows 
\begin{equation}\label{jing}
\begin{aligned}
0 \leq \frac{\zeta_{\varepsilon}}{\varepsilon} \leq u\left(t_\varepsilon^*, \mu_{1, \varepsilon}^*, \mu^* _{2,\varepsilon}\right)-v\left(s_{\varepsilon}^*, \nu_{1, \varepsilon}^*, \nu_{2, \varepsilon}^*\right)-\eta(T-t_{\varepsilon}^*+T-s_{\varepsilon}^*)-\ell\leq M-m-\ell<\infty .
\end{aligned}
\end{equation}

Since $\overline{\mathcal{O}}_N$ is compact, for any sequence $0<\varepsilon_i \downarrow 0\ (\text{as } i\rightarrow \infty)$, there exists subsequence of $\left\{\left(t_{\varepsilon_i}^*, \mu_{1, \varepsilon_i}^*, \mu_{2, \varepsilon_i}^*,s_{\varepsilon_i}^*,\nu_{1, \varepsilon_i}^*, \nu_{2, \varepsilon_i}^*\right)\right\}_{i \in \mathbb{N}}$ (defined by (\ref{max})), still denoted as it, which converges, i.e., there exist $(t^*,\mu_1^*,\mu_2^*),\ (s^*,\nu_1^*,\nu_2^*)\in\overline{\mathcal{O}}_N$, such that
$$
\lim_{i\rightarrow \infty}(t_{\varepsilon_i}^*, \mu_{1, \varepsilon_i}^*, \mu_{2, \varepsilon_i}^*)=(t^*,\mu_1^*,\mu_2^*),\quad \lim_{i\rightarrow \infty}(s_{\varepsilon_i}^*,\nu_{1, \varepsilon_i}^*, \nu_{2, \varepsilon_i}^*)=(s^*,\nu_1^*,\nu_2^*).
$$
From (\ref{jing}),
$$
\frac{1}{\varepsilon}\big(d\left(\mu_{1, \varepsilon}^*, \nu_{1, \varepsilon}^*\right)+d\left(\mu_{2, \varepsilon}^*, \nu_{2, \varepsilon}^*\right)+\left(t^*_{\varepsilon}-s_\varepsilon^*\right)^2\big)=\frac{\zeta_\varepsilon}{\varepsilon}\leq M-m,\quad \varepsilon>0,
$$
it follows that $\zeta_\varepsilon\rightarrow 0$ ($\varepsilon\downarrow 0$), and, hence, $s^*=t^*$.
Next we prove $\left(\mu_1^*, \mu_2^*\right)=\left(\nu_1^*, \nu_2^*\right),\ t^*<T$ and $\displaystyle\limsup _{i \rightarrow \infty} \frac{\zeta_{\varepsilon_i}}{\varepsilon_i}=0$.\smallskip\\
\textbf{(1)} Let us first show $\left(\mu_1^*, \mu_2^*\right)=\left(\nu_1^*, \nu_2^*\right)$:
Since $\displaystyle\lim _{\varepsilon\rightarrow 0} \zeta_{\varepsilon}=0$, it holds that $\displaystyle\lim_{\varepsilon \rightarrow 0} \langle\mu_{i, \varepsilon}^*-\nu_{i,\varepsilon}^*, f_j\rangle=0$, for all $j \geq 1$, $i=1,2$.
As $\left\{x^j\right\}_{j=1}^{\infty} \subset \Theta=\left\{f_j\right\}_{j=1}^{\infty}$, we have, in particular
$$\lim _{\varepsilon \rightarrow 0}\left\langle\mu_{i, \varepsilon}^*-\nu_{i,\varepsilon}^*, x^j\right\rangle=0, \quad \text{for all } j\geq 1,\ i=1,2.$$
From Lemma \ref{contin} we know the mapping $\mu \mapsto\left\langle\mu, x^j\right\rangle$ is continuous on $\mathcal{M}_{Ne^{K^*T}}$. Hence, for $k=1,2,$
$$\left\langle\mu_k^*-\nu_k^*, x^j\right\rangle=\lim _{i \rightarrow \infty}\left\langle\mu_{k,\varepsilon_i}^*-\nu_{k,\varepsilon_i}^*, x^j\right\rangle=0,\ j\geq 1.$$
This implies $\left(\mu^*_1, \mu_2^*\right)=\left(\nu_1^*, \nu^*_2\right)$.
\smallskip\\
\textbf{(2)} We now show that $t^*<T$.
We assume on the contrary that $t^*=T$. Since by hypothesis $(u-v)(T, \cdot, \cdot) \leq 0$, and $u-v$ is continuous, we get
$$
\begin{aligned}
0 \geq (u-v)\left(T, \mu_1^*, \mu_2^*\right) &\geq \limsup_{i\rightarrow \infty}\left(u(t^*_{\varepsilon_i,}, \mu^*_{1,\varepsilon_i}, \mu^*_{2,\varepsilon_i})-v(s^*_{\varepsilon_i,}, \nu^*_{1,\varepsilon_i}, \nu^*_{2,\varepsilon_i})\right)\\
&\geq \limsup_{i\rightarrow \infty}\phi_{\varepsilon_i}(t^*_{\varepsilon_i,}, \mu^*_{1,\varepsilon_i}, \mu^*_{2,\varepsilon_i},s^*_{\varepsilon_i,}, \nu^*_{1,\varepsilon_i}, \nu^*_{2,\varepsilon_i})\geq\ell>0,
\end{aligned}$$
i.e., we have a contradiction. Hence, $t^*<T$.
\smallskip\\
\textbf{(3)} We prove $\displaystyle\limsup _{i \rightarrow \infty} \frac{\zeta_{\varepsilon_i}}{\varepsilon_i}=0.$ In fact, from the continuity of  $u-v$ and the inequalities (\ref{>0}) and (\ref{jing}) we obtain
$$\begin{aligned}
\ell & \geq \phi_\varepsilon\left(t^*, \mu_1^*,\mu_2^*, t^*,\mu_1^*,\mu_2^*\right) =u\left(t^*, \mu_1^*,\mu_2^*\right)-v\left(t^*,\mu_1^*,\mu_2^*\right)-2 \eta\left(T-t^*\right) \\
& \geq \limsup _{i \rightarrow \infty}\left(u\left(t_{\varepsilon_i}^*,  \mu^*_{1,\varepsilon_i}, \mu^*_{2,\varepsilon_i}\right)-v\left(s_{\varepsilon_i}^*, \nu^*_{1,\varepsilon_i}, \nu^*_{2,\varepsilon_i}\right)-\eta\left(T-t_{\varepsilon_i}^*+T-s_{\varepsilon_i}^*\right)\right) \\
& \geq \ell+\limsup _{i \rightarrow \infty} \frac{1}{\varepsilon_i}\left(d( \mu^*_{1,\varepsilon_i},\nu^*_{1,\varepsilon_i})+d(\mu^*_{2,\varepsilon_i},\nu^*_{2,\varepsilon_i})+\left(t_{\varepsilon_i}^*-s_{\varepsilon_i}^*\right)^2\right) =\ell+\limsup _{i \rightarrow \infty} \frac{\zeta_{\varepsilon_i}}{\varepsilon_i} .
\end{aligned}$$
As $\displaystyle\frac{\zeta_{\varepsilon_i}}{\varepsilon_i}\geq 0$, for all $i\geq 1$, this yields $\displaystyle\limsup _{i \rightarrow \infty} \frac{\zeta_{\varepsilon_i}}{\varepsilon_i}=0.$\smallskip\\
\textbf{Step 3.\ Initial estimate.}

For $s,t\in [0,T]$, $y\in\mathbb{R}$, we introduce the notations
\begin{equation}\label{pi}
\begin{gathered}
\beta_{\eta,\varepsilon}(t,s):=\frac{1}{\varepsilon}(t-s)^2+\eta(T-t+T-s),\\
\pi^*_{1,\varepsilon}(y):=\frac{2}{\varepsilon}\sum^{\infty}_{j=1}c_j\langle\mu^*_{1,\varepsilon}-\nu^*_{1,\varepsilon},f_j\rangle f'_j(y), \ \ \pi^*_{2,\varepsilon}(y):=\frac{2}{\varepsilon}\sum^{\infty}_{j=1}c_j\langle\mu^*_{2,\varepsilon}-\nu^*_{2,\varepsilon},f_j\rangle f'_j(y).
\end{gathered}\end{equation}
Then recalling (\ref{varphi_u}) and (\ref{varphi_v}), we see that, for $s,t\in [0,T]$, $\mu_1,\mu_2,\nu_1,\nu_2\in\mathcal{M}_{Ne^{K^*T}}$, $y\in\mathbb{R}$,
$$
\begin{gathered}
\partial_t \varphi_u(t,\mu_1,\mu_2)=\partial_t\beta_{\eta,\varepsilon}(t,s_\varepsilon^*),\ \ \partial_s \varphi_v(s,\nu_1,\nu_2)=-\partial_s\beta_{\eta,\varepsilon}(t_\varepsilon^*,s),\\
\partial_{\mu_1}\varphi_u(t,\mu_{1,\varepsilon}^*,\mu_{2,\varepsilon}^*;y)=\pi^*_{1,\varepsilon}(y),\ \
\partial_{\mu_2}\varphi_u(t,\mu_{1,\varepsilon}^*,\mu_{2,\varepsilon}^*;y)=\pi^*_{2,\varepsilon}(y),\\
\partial_{\nu_1}\varphi_v(s,\nu_{1,\varepsilon}^*,\nu_{2,\varepsilon}^*;y)=\pi^*_{1,\varepsilon}(y),\ \
\partial_{\nu_2}\varphi_v(s,\nu_{1,\varepsilon}^*,\nu_{2,\varepsilon}^*;y)=\pi^*_{2,\varepsilon}(y).
\end{gathered}
$$
From Step 1 we have
$$
\begin{aligned}
&(\partial_t \beta_{\eta,\varepsilon}(t_\varepsilon^*,s_\varepsilon^*),\pi^*_{1,\varepsilon}(\cdot),\pi^*_{2,\varepsilon}(\cdot))
\\&=(\partial_t\varphi_u(t_\varepsilon^*,\mu_{1,\varepsilon}^*,\mu_{2,\varepsilon}^*),\partial_{\mu_1}\varphi_u(t_\varepsilon^*,\mu_{1,\varepsilon}^*,\mu_{2,\varepsilon}^*;\cdot),\partial_{\mu_2}\varphi_u(t_\varepsilon^*,\mu_{1,\varepsilon}^*,\mu_{2,\varepsilon}^*;\cdot))
\in J^{+}_{\overline{\mathcal{O}}_N}u(t_\varepsilon^*,\mu_{1,\varepsilon}^*,\mu_{2,\varepsilon}^*),\\
&(-\partial_s \beta_{\eta,\varepsilon}(t_\varepsilon^*,s_\varepsilon^*),\pi^*_{1,\varepsilon}(\cdot),\pi^*_{2,\varepsilon}(\cdot))
\\&=(\partial_s\varphi_v(s_\varepsilon^*,\nu_{1,\varepsilon}^*,\nu_{2,\varepsilon}^*),\partial_{\nu_1}\varphi_v(s_\varepsilon^*,\nu_{1,\varepsilon}^*,\nu_{2,\varepsilon}^*;\cdot),\partial_{\nu_2}\varphi_v(s_\varepsilon^*,\nu_{1,\varepsilon}^*,\nu_{2,\varepsilon}^*;\cdot))
\in J^{-}_{\overline{\mathcal{O}}_N}v(t_\varepsilon^*,\mu_{1,\varepsilon}^*,\mu_{2,\varepsilon}^*).
\end{aligned}
$$
Then by the assumption that $u$ is a viscosity sub- and $v$ a supersolution on $\mathcal{O}_N$, and as $t^*<T$, we have
$$
\begin{aligned}
&-\partial_t \beta_{\eta,\varepsilon}(t_\varepsilon^*,s_\varepsilon^*)-\mathcal{H}(t_\varepsilon^*,\mu_{1,\varepsilon}^*,\mu_{2,\varepsilon}^*,\pi^*_{1,\varepsilon}(\cdot),\pi^*_{2,\varepsilon}(\cdot))\\
&\leq 0\leq\ \partial_s \beta_{\eta,\varepsilon}(t_\varepsilon^*,s_\varepsilon^*)-\mathcal{H}(s_\varepsilon^*,\nu_{1,\varepsilon}^*,\nu_{2,\varepsilon}^*,\pi^*_{1,\varepsilon}(\cdot),\pi^*_{2,\varepsilon}(\cdot)).
\end{aligned}
$$
By the definition of $\beta_{\eta,\varepsilon}$, this yields
\begin{equation}\label{maodun}
\begin{aligned}
&\mathcal{H}(s_\varepsilon^*,\nu_{1,\varepsilon}^*,\nu_{2,\varepsilon}^*,\pi^*_{1,\varepsilon}(\cdot),\pi^*_{2,\varepsilon}(\cdot))-\mathcal{H}(t_\varepsilon^*,\mu_{1,\varepsilon}^*,\mu_{2,\varepsilon}^*,\pi^*_{1,\varepsilon}(\cdot),\pi^*_{2,\varepsilon}(\cdot))
\\&\leq \partial_t \beta_{\eta,\varepsilon}(t_\varepsilon^*,s_\varepsilon^*)+\partial_s \beta_{\eta,\varepsilon}(t_\varepsilon^*,s_\varepsilon^*)=-2\eta<0.
\end{aligned}\end{equation}
On the other hand, from the definition of Hamiltonian $\mathcal{H}$, with the notation
$$
\begin{aligned}\mathcal{G}:=\big\{(\widetilde{\gamma}_1,\widetilde{\gamma}_2,\overline{\gamma}_1,\overline{\gamma}_2):
\widetilde{\gamma}_i\in\Pi_{\nu_{i,\varepsilon}^*},
\overline{\gamma}_i\in\Pi_{\mu_{i,\varepsilon}^*},\ \text{s.t. } \widetilde{\gamma}_i(\mathbb{R}\times \cdot)=\overline{\gamma}_i(\mathbb{R}\times \cdot),\ i=1,2 \big\},
\end{aligned}$$
 we have
\begin{equation}\label{end}
\begin{aligned}
&\mathcal{H}(s_\varepsilon^*,\nu_{1,\varepsilon}^*,\nu_{2,\varepsilon}^*,\pi^*_{1,\varepsilon}(\cdot),\pi^*_{2,\varepsilon}(\cdot))
-\mathcal{H}(t_\varepsilon^*,\mu_{1,\varepsilon}^*,\mu_{2,\varepsilon}^*,\pi^*_{1,\varepsilon}(\cdot),\pi^*_{2,\varepsilon}(\cdot))
\\&\geq \inf_{(\widetilde{\gamma}_1,\widetilde{\gamma}_2,\overline{\gamma}_1,\overline{\gamma}_2)\in\mathcal{G}}\Big\{
\big\langle \widetilde{\gamma}_1,\mathcal{L}^{\nu_{1,\varepsilon}^*,\nu_{2,\varepsilon}^*,\widetilde{\gamma}_2}_{s_{\varepsilon}^*}[\pi^*_{1,\varepsilon}(\cdot)]\big\rangle+
\big\langle \widetilde{\gamma}_2,\overline{\mathcal{L}}^{\nu_{1,\varepsilon}^*,\nu_{2,\varepsilon}^*,\widetilde{\gamma}_2}_{s_{\varepsilon}^*}[\pi^*_{2,\varepsilon}(\cdot)]\big\rangle
 \\&\quad \ \qquad\qquad\qquad-\big\langle \overline{\gamma}_1,\mathcal{L}^{\mu_{1,\varepsilon}^*,\mu_{2,\varepsilon}^*,\overline{\gamma}_2}_{t_{\varepsilon}^*}[\pi^*_{1,\varepsilon}(\cdot)]\big\rangle
-\big\langle \overline{\gamma}_2,\overline{\mathcal{L}}^{\mu_{1,\varepsilon}^*,\mu_{2,\varepsilon}^*,\overline{\gamma}_2}_{t_{\varepsilon}^*}[\pi^*_{2,\varepsilon}(\cdot)]\big\rangle\Big\}\\
&=-\sup_{(\widetilde{\gamma}_1,\widetilde{\gamma}_2,\overline{\gamma}_1,\overline{\gamma}_2)\in\mathcal{G}}I^{\widetilde{\gamma}_1,\widetilde{\gamma}_2,\overline{\gamma}_1,\overline{\gamma}_2}.
\end{aligned}
\end{equation}
Above we have used the abbreviating notation
$$
I^{\widetilde{\gamma}_1,\widetilde{\gamma}_2,\overline{\gamma}_1,\overline{\gamma}_2}:=
\sum_{i=1}^4 I_i(\widetilde{\gamma}_1,\widetilde{\gamma}_2,\overline{\gamma}_1,\overline{\gamma}_2)$$
with
$$\begin{aligned}
&I_1(\widetilde{\gamma}_1,\widetilde{\gamma}_2,\overline{\gamma}_1,\overline{\gamma}_2):=\big\langle \overline{\gamma}_1-\widetilde{\gamma}_1,\mathcal{L}^{\mu_{1,\varepsilon}^*,\mu_{2,\varepsilon}^*,\overline{\gamma}_2}_{t_{\varepsilon}^*}[\pi^*_{1,\varepsilon}(\cdot)]\big\rangle;\\
& I_2(\widetilde{\gamma}_1,\widetilde{\gamma}_2,\overline{\gamma}_1,\overline{\gamma}_2):=\big\langle \widetilde{\gamma}_1,{\mathcal{L}}^{\mu_{1,\varepsilon}^*,\mu_{2,\varepsilon}^*,\overline{\gamma}_2}_{t_{\varepsilon}^*}[\pi^*_{1,\varepsilon}(\cdot)]
-\mathcal{L}^{\nu_{1,\varepsilon}^*,\nu_{2,\varepsilon}^*,\widetilde{\gamma}_2}_{s_{\varepsilon}^*}[\pi^*_{1,\varepsilon}(\cdot)]\big\rangle;\\
&I_3(\widetilde{\gamma}_1,\widetilde{\gamma}_2,\overline{\gamma}_1,\overline{\gamma}_2):=\big\langle \overline{\gamma}_2-\widetilde{\gamma}_2,\mathcal{L}^{\mu_{1,\varepsilon}^*,\mu_{2,\varepsilon}^*,\overline{\gamma}_2}_{t_{\varepsilon}^*}[\pi^*_{2,\varepsilon}(\cdot)]\big\rangle;\\
& I_4(\widetilde{\gamma}_1,\widetilde{\gamma}_2,\overline{\gamma}_1,\overline{\gamma}_2):=
\big\langle \widetilde{\gamma}_2,\overline{\mathcal{L}}^{\mu_{1,\varepsilon}^*,\mu_{2,\varepsilon}^*,\overline{\gamma}_2}_{t_{\varepsilon}^*}[\pi^*_{2,\varepsilon}(\cdot)]
-\overline{\mathcal{L}}^{\nu_{1,\varepsilon}^*,\nu_{2,\varepsilon}^*,\widetilde{\gamma}_2}_{s_{\varepsilon}^*}[\pi^*_{2,\varepsilon}(\cdot)]\big\rangle.
\end{aligned}$$
Furthermore, we also need
$$\begin{gathered}
I_1^*:=\sup_\mathcal{G}I_1(\widetilde{\gamma}_1,\widetilde{\gamma}_2,\overline{\gamma}_1,\overline{\gamma}_2),\
I_2^*:=\sup_\mathcal{G}I_2(\widetilde{\gamma}_1,\widetilde{\gamma}_2,\overline{\gamma}_1,\overline{\gamma}_2),\\
I_3^*:=\sup_\mathcal{G}I_3(\widetilde{\gamma}_1,\widetilde{\gamma}_2,\overline{\gamma}_1,\overline{\gamma}_2),\
I_4^*:=\sup_\mathcal{G}I_4(\widetilde{\gamma}_1,\widetilde{\gamma}_2,\overline{\gamma}_1,\overline{\gamma}_2).
\end{gathered}$$
Then, obviously,
$$
\begin{aligned}
\mathcal{H}(s_\varepsilon^*,\nu_{1,\varepsilon}^*,\nu_{2,\varepsilon}^*,\pi^*_{1,\varepsilon}(\cdot),\pi^*_{2,\varepsilon}(\cdot))
-\mathcal{H}(t_\varepsilon^*,\mu_{1,\varepsilon}^*,\mu_{2,\varepsilon}^*,\pi^*_{1,\varepsilon}(\cdot),\pi^*_{2,\varepsilon}(\cdot))\geq -(I_1^*+I_2^*+I_3^*+I_4^*)
.
\end{aligned}$$
In the next two steps, we study the estimates of $I^*_1,I^*_2,I^*_3,I^*_4$, in order to obtain a contradiction to equation (\ref{maodun}).\\
\textbf{Step 4.\ Estimates of $I_1^*$ and $I_3^*$.}

From (\ref{L}) and (\ref{pi}), we have
$$\begin{aligned}
I_1^*\leq\sup_{(\widetilde{\gamma}_1,\widetilde{\gamma}_2,\overline{\gamma}_1,\overline{\gamma}_2)\in\mathcal{G}}\frac{2}{\varepsilon}
\sum_{j=1}^{\infty}c_j\big|\langle\mu^*_{1,\varepsilon}-\nu^*_{1,\varepsilon},f_j\rangle
\big\langle \overline{\gamma}_1-\widetilde{\gamma}_1,\mathcal{L}^{\mu_{1,\varepsilon}^*,\mu_{2,\varepsilon}^*,\overline{\gamma}_2}_{t_{\varepsilon}^*}[f_j']\big\rangle
\big|\leq I_{1,b}^*+I_{1,\sigma}^*,
\end{aligned}
$$
where due to Assumption \ref{A1}-(i) that $b_2$ is bounded by $C_0$, it follows that
$$\begin{aligned}
I_{1,b}^*:=&\sup_{(\widetilde{\gamma}_1,\widetilde{\gamma}_2,\overline{\gamma}_1,\overline{\gamma}_2)\in\mathcal{G}}\frac{2}{\varepsilon}
\sum_{j=1}^{\infty}c_j\big|\langle\mu^*_{1,\varepsilon}-\nu^*_{1,\varepsilon},f_j\rangle
\int_{\mathbb{R}\times U}f'_j(y)b_2(t_{\varepsilon}^*,\overline{\gamma}_2)( \overline{\gamma}_1-\widetilde{\gamma}_1)(dvdy)\big|
\\&\leq C_0\frac{2}{\varepsilon}
\sum_{j=1}^{\infty}c_j\big|\langle\mu^*_{1,\varepsilon}-\nu^*_{1,\varepsilon},f_j\rangle
\langle \mu_{1,\varepsilon}^*-\nu_{1,\varepsilon}^*,f'_j\rangle\big|.
\end{aligned}$$
By the construction of $\Theta=\left\{f_j\right\}_{j=1}^{\infty}$, for every $j\in\mathbb{N}$, there
exists an index $k_1(j)$ such that $f_j^{\prime}=$ $f_{k_1(j)}$, and as $f_j^{\prime}=f_{k_1(j)} \in \chi\left(f_j\right)$,  we have $c_j \leq c_{k_1(j)}$. Consequently,
$$\begin{aligned}
I_{1,b}^* &\leq C_0\frac{1}{\varepsilon}
\big(\sum_{j=1}^{\infty}c_j\langle\mu^*_{1,\varepsilon}-\nu^*_{1,\varepsilon},f_j\rangle^2+
\sum_{j=1}^{\infty}c_{k_1(j)}\langle \mu_{1,\varepsilon}^*-\nu_{1,\varepsilon}^*,f_{k_1(j)}\rangle^2\big)
\\&\leq C_0\frac{2}{\varepsilon}
\sum_{j=1}^{\infty}c_j\langle\mu^*_{1,\varepsilon}-\nu^*_{1,\varepsilon},f_j\rangle^2\rightarrow 0,
\end{aligned}$$
as $\varepsilon\rightarrow 0$.
The convergence result is due to Step 2-(3).

We similarly define and estimate $I_{1,\sigma}^*$. For any $j\in\mathbb{N}$, there exists an index $k_2(j)$ such that $f_j''=f_{k_2(j)}$ and $c_j\leq c_{k_2(j)}$. Also from Assumption \ref{A1}-(i) of the boundedness of $\sigma_2(t,\gamma)$ by $C_0$, we have
$$\begin{aligned}
I_{1,\sigma}^*:=&\sup_{(\widetilde{\gamma}_1,\widetilde{\gamma}_2,\overline{\gamma}_1,\overline{\gamma}_2)\in\mathcal{G}}\frac{1}{\varepsilon}
\sum_{j=1}^{\infty}c_j\big|\langle\mu^*_{1,\varepsilon}-\nu^*_{1,\varepsilon},f_j\rangle
\int_{\mathbb{R}\times U}f''_j(y)(\sigma_2(t_{\varepsilon}^*,\overline{\gamma}_2))^2( \overline{\gamma}_1-\widetilde{\gamma}_1)(dvdy)\big|
\\&\leq C^2_0\frac{1}{\varepsilon}
\sum_{j=1}^{\infty}c_j\big|\langle\mu^*_{1,\varepsilon}-\nu^*_{1,\varepsilon},f_j\rangle
\langle \mu_{1,\varepsilon}^*-\nu_{1,\varepsilon}^*,f''_j\rangle\big|
\\&\leq C^2_0\frac{1}{2\varepsilon}
\big(\sum_{j=1}^{\infty}c_j\langle\mu^*_{1,\varepsilon}-\nu^*_{1,\varepsilon},f_j\rangle^2+
\sum_{j=1}^{\infty}c_{k_2(j)}\langle \mu_{1,\varepsilon}^*-\nu_{1,\varepsilon}^*,f_{k_2(j)}\rangle^2\big)
\\&\leq C^2_0\frac{1}{\varepsilon}
\sum_{j=1}^{\infty}c_j\langle\mu^*_{1,\varepsilon}-\nu^*_{1,\varepsilon},f_j\rangle^2\rightarrow 0,
\end{aligned}$$
as $\varepsilon\rightarrow 0$.

Repeating the same argument for $I_3^*$, we conclude that there exists a function $g_3$ of $\varepsilon$ such that
 $$I_3^*\leq I_{3,b}^*+I_{3,\sigma}^*\leq g_3(\varepsilon)\rightarrow 0,\quad \text{as } \varepsilon\rightarrow 0,$$
 where $I_{3,b}^*$ and $I_{3,\sigma}^*$ are defined similarly to $I_{1,b}^*$ and $I_{1,\sigma}^*$, as we have discussed above.
 \\
\textbf{Step 5.\ Estimates of $I_2^*$ and $I_4^*$.}

From (\ref{L}) and (\ref{pi}), we have
$$\begin{aligned}
I_2^*\leq\sup_{(\widetilde{\gamma}_1,\widetilde{\gamma}_2,\overline{\gamma}_1,\overline{\gamma}_2)\in\mathcal{G}}\frac{2}{\varepsilon}
\sum_{j=1}^{\infty}c_j\big|\langle\mu^*_{1,\varepsilon}-\nu^*_{1,\varepsilon},f_j\rangle
\big\langle \widetilde{\gamma}_1,{\mathcal{L}}^{\mu_{1,\varepsilon}^*,\mu_{2,\varepsilon}^*,\overline{\gamma}_2}_{t_{\varepsilon}^*}[f_j']
-\mathcal{L}^{\nu_{1,\varepsilon}^*,\nu_{2,\varepsilon}^*,\widetilde{\gamma}_2}_{s_{\varepsilon}^*}[f_j']\big\rangle
\big|\leq I_{2,b}^*+I_{2,\sigma}^*,
\end{aligned}
$$
where $I_{2,b}^*$ and $I_{2,\sigma}^*$ are defined and estimated as follows. With the help of Assumption \ref{a2}-(i), (\ref{c_j}) and Step 2-(3), we deduce
$$\begin{aligned}
I_{2,b}^*:&=\sup_{(\widetilde{\gamma}_1,\widetilde{\gamma}_2,\overline{\gamma}_1,\overline{\gamma}_2)\in\mathcal{G}}\frac{2}{\varepsilon}
\sum_{j=1}^{\infty}c_j\Big|\langle\mu^*_{1,\varepsilon}-\nu^*_{1,\varepsilon},f_j\rangle
\int_{\mathbb{R}\times U}f'_j(y)|b_2(t_{\varepsilon}^*,\overline{\gamma}_2)-b_2(s_{\varepsilon}^*,\widetilde{\gamma}_2)|\widetilde{\gamma}_1(dvdy)\Big|
\\&\leq \frac{1}{\varepsilon}
\sum_{j=1}^{\infty}c_j\Big(\langle\mu^*_{1,\varepsilon}-\nu^*_{1,\varepsilon},f_j\rangle^2+
\big(\int_{\mathbb{R}\times U}f'_j(y)|b_2(t_{\varepsilon}^*,\overline{\gamma}_2)-b_2(s_{\varepsilon}^*,\widetilde{\gamma}_2)|\widetilde{\gamma}_1(dvdy)\big)^2\Big)
\\&\leq \frac{1}{\varepsilon}
\sum_{j=1}^{\infty}c_j\langle\mu^*_{1,\varepsilon}-\nu^*_{1,\varepsilon},f_j\rangle^2+
\frac{2}{\varepsilon}\kappa_0^2
\Big((t_\varepsilon^*-s_\varepsilon^*)^2+\sum_{j=1}^{\infty}c_j\langle\mu^*_{2,\varepsilon}-\nu^*_{2,\varepsilon},f_j\rangle^2\Big)
\sum_{j=1}^{\infty}c_{k_1(j)}\langle \nu_{1,\varepsilon}^*,f_{k_1(j)}\rangle^2
\\&\leq C\frac{1}{\varepsilon}\Big((t_\varepsilon^*-s_\varepsilon^*)^2+
\sum_{j=1}^{\infty}c_j\langle\mu^*_{1,\varepsilon}-\nu^*_{1,\varepsilon},f_j\rangle^2+
\sum_{j=1}^{\infty}c_j\langle\mu^*_{2,\varepsilon}-\nu^*_{2,\varepsilon},f_j\rangle^2\Big)
\rightarrow 0,
\end{aligned}$$
as $\varepsilon\rightarrow 0$, where $C>0$ is a constant. Similarly, from Assumption \ref{A1}-(i) and Assumption \ref{a2}-(i),
$$\begin{aligned}
&I_{2,\sigma}^*:=\sup_{(\widetilde{\gamma}_1,\widetilde{\gamma}_2,\overline{\gamma}_1,\overline{\gamma}_2)\in\mathcal{G}}\frac{1}{\varepsilon}
\sum_{j=1}^{\infty}c_j\Big|\langle\mu^*_{1,\varepsilon}-\nu^*_{1,\varepsilon},f_j\rangle
\int_{\mathbb{R}\times U}f''_j(y)\big|(\sigma_2(t_{\varepsilon}^*,\overline{\gamma}_2))^2
-(\sigma_2(s_{\varepsilon}^*,\widetilde{\gamma}_2))^2\big|\widetilde{\gamma}_1(dvdy)\Big|
\\&\leq \frac{1}{2\varepsilon}
\sum_{j=1}^{\infty}c_j\langle\mu^*_{1,\varepsilon}-\nu^*_{1,\varepsilon},f_j\rangle^2
+\frac{4}{\varepsilon}C^2_0\kappa_0^2
\Big((t_\varepsilon^*-s_\varepsilon^*)^2+\sum_{j=1}^{\infty}c_j\langle\mu^*_{2,\varepsilon}-\nu^*_{2,\varepsilon},f_j\rangle^2\Big)
\sum_{j=1}^{\infty}c_{k_2(j)}\langle \nu_{1,\varepsilon}^*,f_{k_2(j)}\rangle^2
\\&\leq C\frac{1}{\varepsilon}\Big((t_\varepsilon^*-s_\varepsilon^*)^2+
\sum_{j=1}^{\infty}c_j\langle\mu^*_{1,\varepsilon}-\nu^*_{1,\varepsilon},f_j\rangle^2+
\sum_{j=1}^{\infty}c_j\langle\mu^*_{2,\varepsilon}-\nu^*_{2,\varepsilon},f_j\rangle^2\Big)
\rightarrow 0,
\end{aligned}$$
as $\varepsilon\rightarrow 0$, where $C>0$ is a constant.

Repeating the same argument for $I_4^*$, we conclude that there exist a function $g_4$, such that
 $$I_4^*\leq I_{4,b}^*+I_{4,\sigma}^*\leq g_4(\varepsilon)\rightarrow 0,\quad \text{as } \varepsilon\rightarrow 0,$$
 where $I_{4,b}^*$ and $I_{4,\sigma}^*$ can be defined similarly to $I_{2,b}^*$ and $I_{2,\sigma}^*$ as we discussed above.\\
\textbf{Step 6.\ Conclusion.}

In the Steps 4 and 5, we have proved that there exists a function $G$, such that as $\varepsilon\rightarrow 0$,
\begin{equation}\label{relation}
-\sup_{(\widetilde{\gamma}_1,\widetilde{\gamma}_2,\overline{\gamma}_1,\overline{\gamma}_2)\in\mathcal{G}}I^{\widetilde{\gamma}_1,\widetilde{\gamma}_2,\overline{\gamma}_1,\overline{\gamma}_2}
\geq -(I_1^*+I_2^*+I_3^*+I_4^*)\geq G(\varepsilon)\rightarrow 0.
\end{equation}
However, by (\ref{maodun}) and (\ref{end}) we have
$$
-\sup_{(\widetilde{\gamma}_1,\widetilde{\gamma}_2,\overline{\gamma}_1,\overline{\gamma}_2)\in\mathcal{G}}I^{\widetilde{\gamma}_1,\widetilde{\gamma}_2,\overline{\gamma}_1,\overline{\gamma}_2}
\leq -2\eta <0,$$
which leads to a contradiction with (\ref{relation}). The proof is complete.
\endpf
\vspace{3mm}
The following corollary is a straightforward conclusion of Theorems \ref{verify} and \ref{comparison}.
\begin{corollary}
Under the assumption that \begin{equation}
 b_j(t,(y,v),\gamma)=b_j(t,\gamma),\quad \sigma_j(t,(y,v),\gamma)=\sigma(t,\gamma),
 \end{equation}
$(t,(y,v),\gamma)\in [0,T]\times(\mathbb{R}\times U)\times \mathcal{P}_2(\mathbb{R}\times U),\ j=1,2,$ are independent of $(y,v)$, the value function $\vartheta$ is the unique viscosity solution to HJB equation (\ref{HJB}) on $\mathcal{O}$.
\end{corollary}

\end{document}